\documentclass[11pt,a4paper]{article}
\usepackage[T1]{fontenc}
\usepackage[utf8]{inputenc}
\usepackage[top=22mm,bottom=24mm,left=24mm,right=24mm]{geometry}
\usepackage{amsmath,amssymb,amsthm,newpxtext,newpxmath,microtype,xcolor}
\definecolor{ink}{HTML}{202928}
\definecolor{pine}{HTML}{17584F}
\definecolor{rust}{HTML}{A2462B}
\definecolor{sand}{HTML}{E9E2CF}
\definecolor{mist}{HTML}{DCE9E6}
\definecolor{plum}{HTML}{5B3F78}
\definecolor{chg}{HTML}{1F5FA8}
\usepackage{tikz,pgfplots}
\pgfplotsset{compat=1.18}
\usetikzlibrary{arrows.meta,positioning,calc,decorations.pathreplacing,patterns,shapes.geometric,fit,backgrounds}
\usepackage{booktabs,tabularx,array}
\usepackage{placeins}
\usepackage{caption}
\usepackage[colorlinks=true,urlcolor=pine,linkcolor=pine,citecolor=pine,pdftitle={The Thurston algorithm for general entire maps on the complex plane: Existence and uniqueness},pdfauthor={Sebastian van Strien},pdfsubject={MSC 2020: 37F20; 37F10; 37F34; 30D05; 30F60},pdfkeywords={Thurston equivalence, entire maps, Levy cycles, Teichmuller space, extremal length, quadratic differentials}]{hyperref}
\newtheorem{lemma}{Lemma}[section]
\newtheorem{prop}[lemma]{Proposition}
\newtheorem{theorem}[lemma]{Theorem}
\theoremstyle{remark}
\newtheorem{remark}[lemma]{Remark}
\theoremstyle{plain}
\newtheorem*{thmA}{Theorem A}
\newtheorem*{thmB}{Theorem B}
\newtheorem{pushforwardlemma}{Lemma}[section]

\newtheorem{matrixlemma}[pushforwardlemma]{Lemma}
\newtheorem{analyticprep}[pushforwardlemma]{Lemma}
\numberwithin{equation}{section}
\newcommand{\CC}{\mathbb C}\newcommand{\RR}{\mathbb R}
\newcommand{\T}{\mathcal T}\newcommand{\Ext}{\operatorname{Ext}}\newcommand{\Mod}{\operatorname{Mod}}
\newcommand{\dd}{\mathrm d}\newcommand{\area}{\operatorname{area}}
\newcommand{\mpr}{Mukundan--Prochorov--Reinke~\cite{MPR}}
\title{The Thurston algorithm for general entire maps on the complex plane:\\
Existence and uniqueness}
\author{Sebastian van Strien}
\date{20 November 2026}
\makeatletter
\renewcommand{\maketitle}{%
  {\small\hfill\@date}\par\vspace{3pt}
  {\bfseries\@title}\par\medskip
  {\large\@author}\par
  {\small Department of Mathematics, Imperial College London\quad
  \href{mailto:s.van-strien@imperial.ac.uk}{\texttt{s.van-strien@imperial.ac.uk}}}\par\medskip}
\makeatother
\begin{document}\color{ink}
\maketitle
\begin{abstract}
We prove that a marked topological entire Thurston map $(f,A)$ of stable parabolic type is Thurston equivalent relative to $A$ to a postsingularly finite entire function if and only if it has no Levy cycle relative to $A\cup\{\infty\}$. Here $A$ is a finite forward-invariant set containing the postsingular set. In the realizable case, the entire function is unique up to complex affine conjugacy, and every orbit of the Thurston pullback converges to its fixed point.

We establish a Wolff--Denjoy type estimate for nonexpanding self-maps $T$ of the Teichm\"uller space of any connected orientable surface of finite type without boundary. If $T$ has an unbounded orbit, there are a nonzero measured foliation $G$ and a constant $C$ such that $\operatorname{Ext}_{T^k x}(G)\le C\operatorname{Ext}_x(G)$ for every $x$ and $k\ge0$. For a fixed-point-free transcendental Thurston pullback, this supplies the foliation used in the existence proof.

A central tool is a fixed finite measure obtained as the total-variation limit of renormalized pushforward measures transported back along a Teichm\"uller ray. It allows us to select a cylinder at a late time and compare the mass of the same marked region at different times. If $G$ has a cylinder, a four-point extremal-length estimate and the Thurston matrix force a Levy cycle. Otherwise, in the absence of a Levy cycle, cluster rescaling shows that cylinders of large modulus carry almost all the fixed measure at all late times. Following a cylinder of definite mass until its modulus decreases gives a contradiction. We compare this geometric argument with the proofs for real entire maps and exponential maps.
\end{abstract}
{\small
\par}

\section{Introduction}\label{sec:intro}
\subsection*{Singular values, forward orbits and the marked plane}\label{bg:singular}
\addcontentsline{toc}{subsection}{Singular values, forward orbits and the marked plane}
For an entire function, local inverses can fail at critical values and at finite asymptotic values, reached along paths escaping to infinity. For example, $e^x\to0$ as $x\to-\infty$, although the exponential has no critical point. The covering definition below treats both phenomena and applies before a complex structure is chosen. A forward-invariant marking containing the singular values ensures that inverse branches over an unmarked disk can be continued through successive pullbacks \cite[\S\S2.3--2.4]{MPR}.

Let $f:\RR^2\to\RR^2$ be noninjective and topologically holomorphic. Locally there are orientation-preserving source and target coordinates in which $f$ is $z\mapsto z^e$ for an integer $e\ge1$. A finite singular value is a point with no evenly covered neighborhood. A disk is evenly covered if its preimage is a nonempty disjoint union of disks mapped homeomorphically onto it. Write $S_f$ for the set of these values and $P_f=\bigcup_{j\ge0}f^j(S_f)$ for their forward orbits. We assume $P_f$ is finite. We also assume \emph{stable parabolic type}: for every target complex structure conformally equivalent to $\CC$, its pullback by $f$ is again conformally equivalent to $\CC$. This hypothesis allows us to uniformize both planes at each step of the argument. Following the terminology of \mpr, we call such an $f$ an entire Thurston map. The degree of the covering over $\RR^2\setminus S_f$ distinguishes topological polynomials (finite degree) from transcendental maps (infinite degree), in plane coordinates. See \mpr, Definitions 2.2, 2.7, 2.9 and 2.14--2.15.

A concrete sufficient condition is that $f$ can be written as $f=\phi\circ g\circ\psi$, where $g:\CC\to\CC$ is a nonconstant entire function and $\phi,\psi$ are orientation-preserving quasiconformal plane homeomorphisms. Since $S_f$ is finite, one may even allow arbitrary orientation-preserving plane homeomorphisms. Such a factorization supplies a target complex structure whose pullback is a plane, so Eremenko's stable-type dichotomy excludes stable hyperbolic type. Conversely, stable parabolic type gives this factorization by uniformization. Thus the equivalence follows from the dichotomy and factorization stated in \cite[Remark~2.10]{MPR}. The two coordinate changes are independent; no isotopy between $\phi^{-1}$ and $\psi$ relative to the marked set is assumed.

Infinity is treated as an end: polynomials extend across it, whereas transcendental maps have an essential singularity there. Curves are considered on the marked sphere, and coverings on the punctured plane.

Fix a finite forward-invariant marked set $A\supseteq P_f$. We call $(f,A)$ a \emph{marked entire Thurston map}. We use the same marked set in the domain and range, and put $A^+=A\cup\{\infty\}$ on the sphere. The point $\infty$ marks the end of the plane. We say that $(f,A)$ is \emph{Thurston equivalent relative to $A$} to an entire function $h$ if there are orientation-preserving homeomorphisms $\phi,\psi:\RR^2\to\CC$ with $\phi f=h\psi$, $\phi|_A=\psi|_A$, and $\phi$ isotopic to $\psi$ relative to $A$. Thus each marked point has a constant image throughout the isotopy. In this case $h$, with marked set $\phi(A)$, is a \emph{realization} of $(f,A)$.

A simple closed curve in the complement of $A^+$ is \emph{essential} if each complementary disk contains at least two marks. A Levy cycle consists of distinct essential isotopy classes $\gamma_0,\ldots,\gamma_{p-1}$ with disjoint representatives such that $f^{-1}(\gamma_{i+1})$ has a circle component in class $\gamma_i$ mapping by degree one, with indices taken modulo $p$. Here and below, pullback of a class means pullback of a representative followed by isotopy relative to the full source marking. Homotopy lifting makes the resulting classes and degrees independent of the representative. The case $p=1$ is included.
\phantomsection\label{bg:statements}
\begin{thmA}[Levy criterion]
A marked entire Thurston map $(f,A)$ is Thurston equivalent relative to $A$ to a postsingularly finite entire function if and only if it has no Levy cycle relative to $A^+$.
\end{thmA}
\begin{thmB}[Uniqueness and convergence]
If $(f,A)$ has a realization, it is unique up to complex affine conjugacy. In that case, the Thurston pullback map of $(f,A)$ has a unique fixed point, and every orbit converges to it.
\end{thmB}
For transcendental maps the existence half of Theorem~A occupies Sections~\ref{sec:records}--\ref{sec:persistence}; polynomial existence, using \cite[\S1.3 and Proposition~B.1]{BMW}, necessity and Theorem~B are treated in Section~\ref{sec:necessity}.
\subsection*{Two examples of the marking}
\addcontentsline{toc}{subsection}{Two examples of the marking}
For $h(z)=z^2-1$, the singular orbit is $-1\mapsto0\mapsto-1$. Add the positive fixed point $a=(1+\sqrt5)/2$ and take $A=\{-1,0,a\}$. The circle $\gamma=\{|w+1|=3/2\}$ and its full preimage $\{|z|=\sqrt{3/2}\}$ both enclose $-1,0$ and omit $a$. Moving the center from $-1$ to $0$ and interpolating the radii gives an isotopy relative to the marks. The lift returns to the same class with degree two, so it gives the Thurston matrix $[1/2]$, not a Levy cycle.

For $h(z)=2\pi i\,e^z$, the omitted singular value has orbit $0\mapsto2\pi i\mapsto2\pi i$. A circle centered at $0$ lifts to a vertical line: each turn changes a logarithm by $2\pi i$, so no finite return closes the lift. A disk avoiding $0$ instead has a closed boundary lift on each inverse branch. This explains the requirement of a \emph{closed preimage component} in a Levy cycle \cite[\S2]{HSS}.

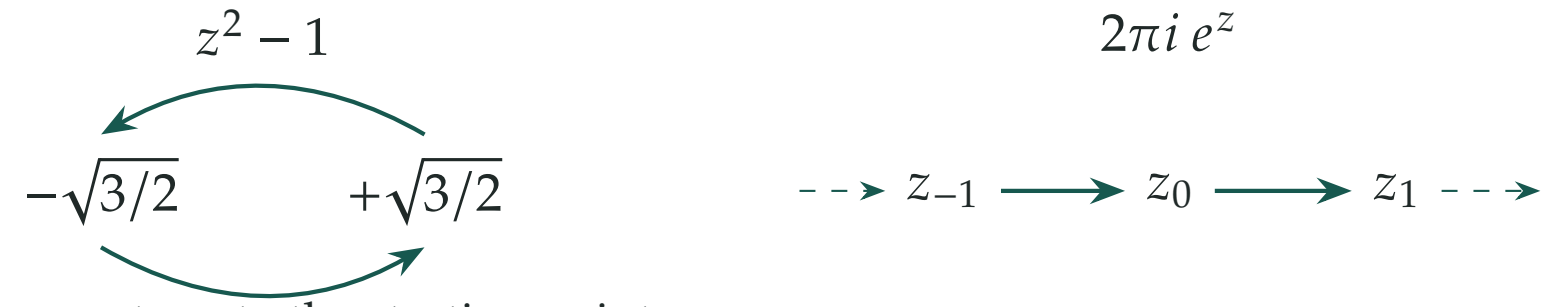
\begin{center}
\begin{tikzpicture}[>=Stealth,font=\small]
\begin{scope}
\node at (0,1.05) {$z^2-1$};
\node (plus) at (1.07,0) {$+\sqrt{3/2}$};
\node (minus) at (-1.07,0) {$-\sqrt{3/2}$};
\draw[->,pine,thick] (plus.north) to[bend right=30] (minus.north);
\draw[->,pine,thick] (minus.south) to[bend right=30] (plus.south);
\node at (0,-.89) {\footnotesize two turns return to the starting point};
\end{scope}
\begin{scope}[shift={(6,0)}]
\node at (0,1.05) {$2\pi i\,e^z$};
\foreach \k/\x in {-1/-1.5,0/0,1/1.5}{\node (p\k) at (\x,0) {$z_{\k}$};}
\draw[->,pine,thick] (p-1)--(p0);
\draw[->,pine,thick] (p0)--(p1);
\draw[->,pine,dashed] (-2.45,0)--(p-1);
\draw[->,pine,dashed] (p1)--(2.45,0);
\node at (0,-.89) {\footnotesize each turn adds $2\pi i$};
\end{scope}
\end{tikzpicture}
\end{center}
Each arrow shows the effect on a preimage of continuing once around the target loop. On the left the loop is $\gamma$ based at $1/2$. On the right it is a small positively oriented circle about $0$, and $z_j=z_0+2\pi i j$ are the preimages of its base point.

The extra fixed point in the polynomial example makes essential curves possible: a sphere with only three marks has none. Four marks already allow infinitely many essential isotopy classes; their partitions of the marks determine a class only within a disjoint family. The torus model on p.~\pageref{bg:four} explains this distinction.

Douady--Hubbard \cite{DH} give the rational characterization; see \cite[\S1.3 and Proposition~B.1]{BMW} for the marked polynomial criterion. Hubbard--Schleicher--Shishikura \cite[Theorem~2.4]{HSS} treat postsingularly finite topological exponential maps with their postsingular marking, including $\infty$; our $P_f$ contains only finite points. Prochorov \cite[Main Theorem~A]{ProSmall} treats four postsingular points in his sphere convention under hypotheses (A)--(C). With additional marks, \cite[Main Theorem~A]{ProExtra} assumes a realization on a smaller marked set and excludes a degenerate Levy multicurve becoming inessential when the additional marks are forgotten. His multicurves use the same disjoint, pairwise nonisotopic convention as our Levy cycles.

\emph{Related work.} The general transcendental realization problem is described as open in \cite{MPR,ProSmall,ProExtra}. Shemyakov \cite{Shem} proves a Levy criterion for postsingularly finite models of $\int_0^z e^{p(t)}\,\dd t$ and their compositions. His Conjecture~C asks for realization in the structurally finite family $\int u\,e^p$, with $u,p$ polynomials. Since a fixed point gives a realization in that family \cite[Theorem~1.1.8 and Proposition~1.2.2]{Shem}, Theorems~A and~B imply the conjecture. D'Souza \cite{DSou} gives a partial characterization for cosine maps. Schleicher has announced joint work with Shemyakov for postsingularly finite maps with hyperbolic orbifold modeled on finite compositions of structurally finite transcendental maps and rational maps \cite[Theorem~1, slide~3]{SchTalk}.

For real entire maps with finitely many singular values, all real, \cite{vS} proves realization of finite labelled real models in a fixed natural family, assuming a nonmonotone real restriction and the interval no-Levy-cycle condition. Its displacement result \cite[Theorem~5.1]{vS}, following \cite[Section~II.4]{dMvS}, motivates the present proof.

A \emph{Thurston obstruction} is a multicurve invariant under essential preimages whose Thurston matrix has spectral radius at least one. Its entries sum reciprocal covering degrees. For entire maps every obstruction contains a Levy cycle \cite[Proposition~1.1.16]{Shem}; the Matrix Lemma proves this for the full marking $A^+$. Conversely, iterating essential preimages of a Levy cycle $\Gamma$ gives nested multicurves, starting with $\Gamma\subseteq P(\Gamma)$ and using preservation of disjointness under pullback. Since each has at most $|A|-2$ classes, they stabilize at an invariant multicurve. Its matrix dominates the cycle's permutation matrix extended by zeros, so its spectral radius is at least one. A Levy cycle is \emph{degenerate} if each prescribed source disk is isotopic to a lift mapped homeomorphically onto the corresponding target disk. Every entire-map Levy cycle is degenerate: the argument principle makes each chosen degree-one boundary lift bound a proper disk of degree one, hence a disk mapped homeomorphically onto the target disk (Section~\ref{sec:closed}). Here disks are the bounded sides, and isotopies fix $A^+$. This explains the formulation in \cite{HSS}. The new implication in Theorem~A is from no Levy cycle to realization, including for the postsingular marking.

The following example shows where the stable parabolic type hypothesis is needed. Let $u:\mathbb D\to\CC\setminus\{0,1\}$ be the universal covering and choose an orientation-preserving homeomorphism $j:\CC\to\mathbb D$ taking the three points $0,1,2$ to three distinct points of $u^{-1}(2)$. Then $f=u\circ j$ has singular values $0,1$ and sends each of $0,1,2$ to $2$. Take $A=\{0,1,2\}$. The bounded side of every essential curve contains two of these points, including at least one omitted value. A closed lift would bound a Jordan disk in the disk uniformization on which the argument principle makes $u$ proper onto that bounded side. This is impossible for a disk containing an omitted value. Thus $f$ has no Levy cycle, but every entire realization would omit two values, contrary to Picard's theorem. The standard target structure pulls back to the disk.

Section~\ref{sec:comparison} compares our geometric estimates with the real displacement argument of \cite{vS} and the near-equality method of \cite{HSS}.

There must be a finite singular value since a connected cover of the plane is one-to-one. With a single finite forward-invariant mark the map is Thurston equivalent to $z\mapsto z^d$ and its marked Teichmüller space is a point (\mpr, Proposition 2.17). With two finite marks the sphere has three marks, its Teichmüller space is a point, and there are no essential curve classes. The fixed-point and rigidity results below give the conclusion in this case. In Sections \ref{sec:records}--\ref{sec:persistence} the map is transcendental and $s=|A|\ge3$.

\subsection*{Idea of the proof}
\addcontentsline{toc}{subsection}{Idea of the proof}
Assume that $(f,A)$ is transcendental, has no Levy cycle relative to $A^+$, and that its pullback $\sigma$ has no fixed point. Contraction on compact subsets of marked Teichmüller space makes its orbit unbounded. The proof has the following steps, summarized in Figure~\ref{fig:flow}.

\emph{One foliation at every depth.} Theorem~\ref{l:record} constructs a nonzero measured foliation $G$ and one constant $C$ such that
\[
 \Ext_{\sigma^j x}G\le C\Ext_xG\qquad(x\in\T,\ j\ge0).
\]
After normalization, choose its Teichmüller ray $g(t)$, along which $\Ext_{g(t)}G=e^{-2t}$. The same bound $Ce^{-2t}$ then holds at every pullback depth. This replaces the displacement estimate in the real proof \cite[Theorem~5.1]{vS}, \cite[Lemmas~II.4.1--II.4.2]{dMvS}.

\emph{Exclude a cylinder of $G$.} Such a cylinder would supply a curve that stays short at every depth. Retaining four source marks detects its class. Proposition~\ref{p:quartet} bounds its leading reciprocal extremal length by finite lifts of vertical cylinders, a degree-$d$ lift contributing its modulus divided by $d$. The resulting Thurston-matrix powers have a uniform positive lower bound. The Matrix Lemma makes those powers tend to zero in the absence of a Levy cycle, a contradiction.

\emph{Construct a fixed measure.} If $G$ has no cylinder, fix one sufficiently large depth $k$. Proposition~\ref{l:transport} uses the differentials realizing the same $G$ on $\sigma^k g(t)$ to construct a finite measure $\mu$, of mass $E>0$, on the initial target. For the renormalized pushforward $P_t$ and ray marking $\phi_t$,
\[
 \int_{\phi_t(U)}|P_t|=\mu(U)+o(1)
\]
uniformly over measurable sets $U$. Thus a cylinder selected on a late surface can be compared at two different times through the mass of its \emph{same} inverse image. Near-equality of the pushforward, used in \cite[Sections~5--6]{HSS}, is a further consequence.

\emph{Locate and follow its mass.} Rescale clusters of marks along any sequence $t_n\to\infty$. A nonzero component limit would, by Proposition~\ref{p:degree}, give a proper disk of degree $d$ containing at least $d+2$ marked poles. The bounded degree and compatible lifted boundary classes force a degree-one return, hence a Levy cycle (Lemma~\ref{l:components}). All component limits therefore vanish. Proposition~\ref{l:cylinder-mass} locates almost all the mass on long flat cylinders between clusters.

There are at most $r=|A|-2$ maximal cylinders of modulus greater than any fixed $K>1$. Consequently one late cylinder carries definite mass. Its core is nonvertical, so its modulus eventually falls to $K$, while its mass for the fixed measure is unchanged. At that time it is disjoint from every cylinder of modulus greater than $K$, contradicting concentration on those cylinders. Section~\ref{sec:persistence} makes this argument uniform over all sufficiently late times.

\clearpage
\subsection*{Main new ingredients}
\begin{figure}[!htp]\centering
\begin{tikzpicture}[
  frame/.style={rounded corners=2pt,inner xsep=7pt,inner ysep=5pt,font=\small},
  shared/.style={frame,draw=pine,fill=mist!65,align=flush center,text width=13.9cm},
  case head/.style={frame,fill=pine,text=white,font=\small\bfseries,
                    align=center,text width=6.65cm},
  proof step/.style={frame,draw=pine!50,fill=white,align=flush left,text width=6.65cm},
  end step/.style={proof step,draw=rust,fill=sand!35},
  line/.style={draw=pine,line width=.65pt},
  arrow/.style={line,-{Stealth[length=1.8mm]}},
  node distance=3mm]
\node[shared] (assume) {Let $f$ be transcendental with $|A|\ge3$ and no Levy cycle.\\
  Suppose its pullback $\sigma$ has no fixed point.};
\node[shared,below=of assume] (foliation) {\textbf{One measured foliation $G\ne0$ works at every depth}
  \quad (Theorem~\ref{l:record})\\[2pt]
  $\Ext_{\sigma^k x}G\le C\,\Ext_xG$ for all $x\in\T$ and $k\ge0$.\\
  On the ray $g(t)$ of $G$, this gives $\Ext_{\sigma^k g(t)}G\le Ce^{-2t}$ by \eqref{eq:ray}.};
\draw[arrow] (assume)--(foliation);

\node[case head,anchor=north] (vertical) at ([xshift=-3.8cm,yshift=-7mm]foliation.south)
  {$G$ has a cylinder of closed leaves};
\node[case head,anchor=north] (nonvertical) at ([xshift=3.8cm,yshift=-7mm]foliation.south)
  {$G$ has no cylinder of closed leaves};
\coordinate (fork) at ([yshift=-3.5mm]foliation.south);
\draw[line] (foliation.south)--(fork);
\draw[arrow] (fork)-|(vertical.north);
\draw[arrow] (fork)-|(nonvertical.north);

\node[proof step,below=of vertical] (core) {\textbf{Follow the core $\gamma$.}
  At every depth $k$, the pullback admits annuli in class $\gamma$ with modulus
  at least $c_\gamma e^{2t}$, where $c_\gamma>0$ is independent of $k$.};
\node[proof step,below=of core] (lifts) {\textbf{Collect the essential lifts.}
  The cylinder cores and their essential lifts at all depths can be drawn
  disjointly. Their classes form a finite invariant multicurve $\Lambda$
  with Thurston matrix $M_\Lambda$.};
\node[proof step,below=of lifts] (matrix) {\textbf{Keep four marks that isolate $\gamma$.}
  A lift of degree $d$ contributes its cylinder modulus divided by $d$.
  For the initial moduli $m$, the four-point bound gives
  $(M_\Lambda^k m)_\gamma\ge c_\gamma$ for every $k\ge1$
  (Proposition~\ref{p:quartet}, \eqref{eq:matrix}).};
\node[end step,below=of matrix] (decay) {\textbf{The matrix powers must decay.}
  In the absence of a Levy cycle, $M_\Lambda^k\to0$. This contradicts the
  positive lower bound for $\gamma$ (Proposition~\ref{l:cylinder-levy}).};

\node[proof step,below=of nonvertical] (measure) {\textbf{Fix one sufficiently large depth.}
  Near equality of areas gives a measure $\mu$ of mass $E>0$ on the initial
  surface, the same at every time (Proposition~\ref{l:transport}).};
\node[proof step,below=of measure] (clusters) {\textbf{Rescale the clusters of marks.}
  A nonzero limit gives a degree-$d$ disk with at least $d+2$ marked poles.
  This bounds $d$ and forces a repeated class among degree-one lifts.
  That would be a Levy cycle, so every component limit vanishes
  (Proposition~\ref{p:degree}, Lemma~\ref{l:components}).};
\node[proof step,below=of clusters] (coverage) {\textbf{Locate the mass between clusters.}
  For fixed $K>1$, at most $|A|-2$ maximal flat cylinders have \mbox{modulus $>K$}.
  Together they carry almost all of $\mu$ at every sufficiently late time
  (Proposition~\ref{l:cylinder-mass}, \eqref{eq:coverage}).};
\node[end step,below=of coverage] (exit) {\textbf{Follow one cylinder of positive mass.}
  One carries a fixed fraction of $E$. Its core is nonvertical, so its
  modulus falls to $K$. It then lies outside all cylinders of \mbox{modulus $>K$},
  contradicting the late-time mass bound (\S\ref{sec:persistence}).};

\foreach \a/\b in {vertical/core,core/lifts,lifts/matrix,matrix/decay,
                    nonvertical/measure,measure/clusters,clusters/coverage,coverage/exit}
  \draw[arrow] (\a.south)--(\b.north);
\coordinate (join) at ([yshift=-3.5mm]foliation.center |- exit.south);
\node[shared,anchor=north] (fixed) at ([yshift=-3.5mm]join)
  {\textbf{The pullback $\sigma$ has a fixed point.}};
\draw[line] (decay.south)|-(join);
\draw[line] (exit.south)|-(join);
\draw[arrow] (join)--(fixed.north);
\foreach \s in {foliation,matrix,measure,clusters}
  \node[anchor=center,text=rust,font=\large] at ([xshift=-7pt,yshift=-8pt]\s.north west) {$\star$};
\end{tikzpicture}
\caption{The existence proof; stars mark the main new ingredients listed above. The absence of a Levy cycle enters in the matrix lemma and in Lemma~\ref{l:components}. Transcendence gives contraction on compact sets, excludes closed lifts of a curve enclosing all finite marks, and excludes a nonzero limit when all finite target marks remain in the same bounded chart in Proposition~\ref{p:degree}.}
\label{fig:flow}
\end{figure}
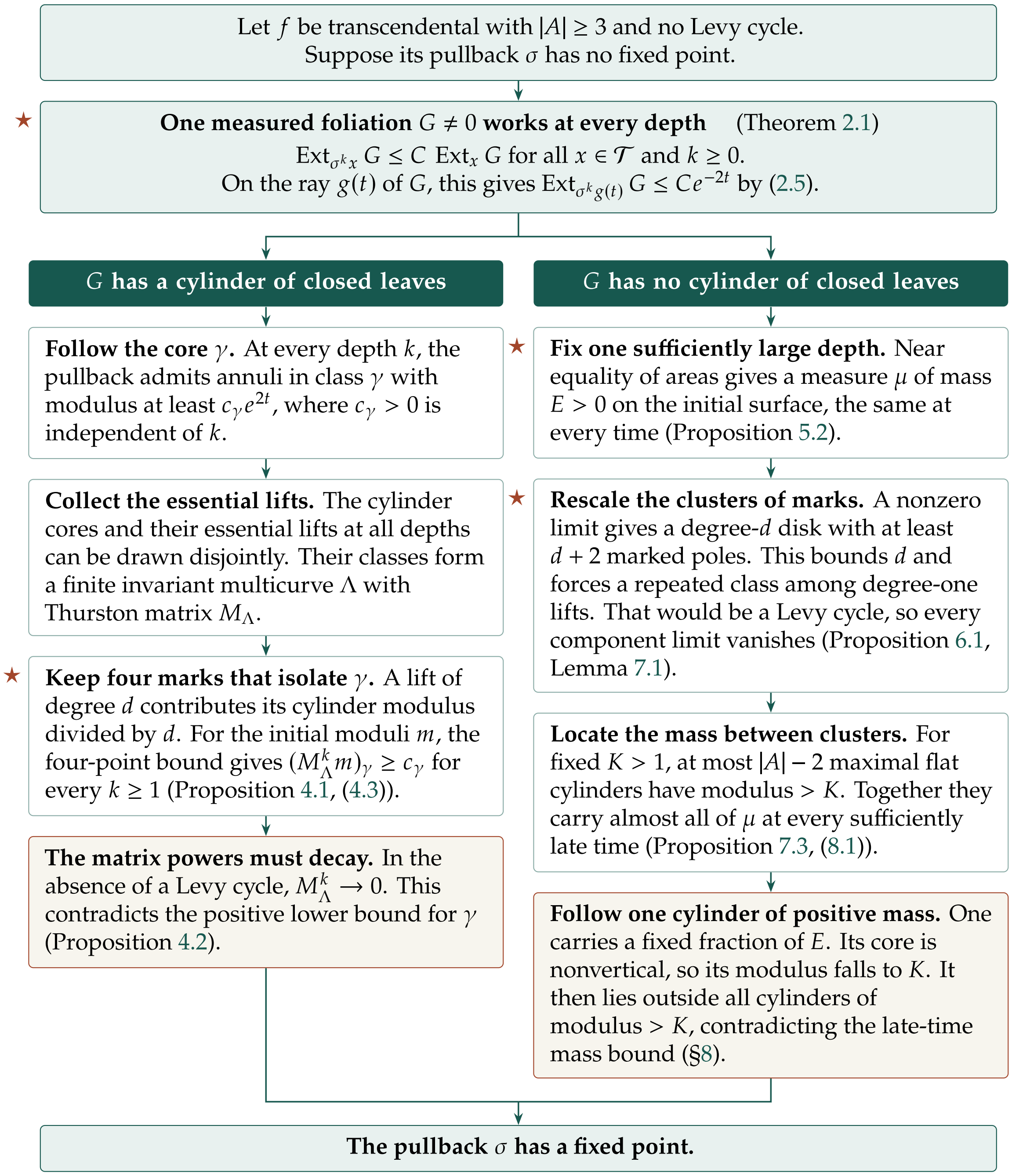

\addcontentsline{toc}{subsection}{Main new ingredients}
Figure~\ref{fig:flow} summarizes the strategy; stars mark the new tools.
\begin{itemize}
\setlength{\itemsep}{0pt}\setlength{\parsep}{0pt}\setlength{\parskip}{0pt}
\item The proof strategy, originating in the real arguments of \cite[Section~II.4]{dMvS} and \cite{vS}.
\item One foliation with an extremal-length bound uniform in point and depth (Theorem~\ref{l:record}).
\item A four-point energy bound for infinite-degree pullbacks (Proposition~\ref{p:quartet}).
\item A fixed measure comparing the same cylinder at different times (Proposition~\ref{l:transport}).
\item Bounded-degree disks from nonzero component limits (Proposition~\ref{p:degree}).
\end{itemize}

\begin{table}[tp]
\subsection*{Notation}
\addcontentsline{toc}{subsection}{Notation}
Table~\ref{tab:notation} lists the symbols used throughout the existence proof.
\par\smallskip
\centering\small
\renewcommand{\arraystretch}{1.12}
\begin{tabularx}{\textwidth}{@{}>{\raggedright\arraybackslash}p{2.7cm}X>{\raggedright\arraybackslash}p{1.9cm}@{}}
\toprule
\textbf{Symbol} & \textbf{Meaning} & \textbf{Introduced}\\
\midrule
$A$, $A^+$, $s$, $r$ & finite marked set; $A^+=A\cup\{\infty\}$; $s=|A|$; $r=s-2$ for $s\ge3$, the largest number of disjoint, pairwise nonisotopic essential curves & \S\ref{sec:intro}, \S\ref{sec:flat}\\
$\T$, $d_T$, $\sigma$ & $\mathrm{Teich}(\widehat\CC,A^+)$, Teichmüller distance, Thurston pullback & \S\ref{sec:records}\\
$o$, $x_j$ & base point, $x_j=\sigma^jo$ & \S\ref{sec:records}\\
$\Ext_x$, $i(\cdot,\cdot)$ & extremal length on $x$, intersection number & \S\ref{sec:records}\\
$G$, $C$ & the foliation and constant of Theorem~\ref{l:record} & \S\ref{sec:records}\\
$\Sigma$, $q$, $g(t)$, $q_t$, $\phi_t$ & initial punctured target at $o$, unit-area differential of $G$, its Teichmüller ray, the stretched differential and ray marking & \S\ref{sec:records}\\
$m_C(t)$ & modulus at time $t$ of a maximal $q$-cylinder $C$ & \S\ref{sec:flat}\\
$M_\Lambda$ & matrix of reciprocal degrees on $\Lambda$ & \S\ref{sec:closed}\\
$\mathcal E_t$, $H_a$ & double cover branched at four marks; least energy of a closed form with period one on $a$ & \S\ref{sec:fourlimit}\\
$k$, $H$, $H_t$, $Q$ & fixed depth; entire representative of $f^k$; its time-$t$ version; $Q=H^*q$ & \S\ref{sec:transport}\\
$r_t$, $E_t$, $E$ & Hubbard--Masur differential of $G$ on $\sigma^kg(t)$; $E_t=e^{2t}\|r_t\|_1$; $E=\lim E_t$ & \S\ref{sec:transport}\\
$v$, $S$, $\mu$, $P_t$ & $v^2$ is the limiting source density; $S=\sum_jv^2\circ\xi_j$; $\mu=S|q|$; $P_t=e^{2t}(H_t)_*r_t$ & \S\ref{sec:transport}\\
$t_n$, $q_n$, $P_n$, $\mu_n$ & a sequence of times and the corresponding objects & \S\ref{sec:mass}\\
$\mathcal U(t,K)$ & union of the maximal cylinders of modulus $>K$ at time $t$ & \S\ref{sec:persistence}\\
\bottomrule
\end{tabularx}
\caption{Notation for the existence proof.}\label{tab:notation}
\end{table}

\newpage
\section{A measured foliation with uniformly bounded extremal-length growth}\label{sec:records}
Theorem~\ref{l:record} gives one measured foliation whose extremal length grows by a bounded factor under every iterate of a nonexpanding map with an unbounded forward orbit. It is inspired by the displacement criterion of \cite[Theorem~5.1]{vS}, following \cite[Section~II.4]{dMvS}: a continuous fixed-point-free self-map of a bounded open convex subset of a finite-dimensional normed space has points approaching the boundary with bounded displacement-to-boundary-distance ratio. Remark~\ref{r:wolff} explains the relation to Wolff--Denjoy theory.

\subsection*{Marked complex structures and the Thurston pullback}\label{bg:mark}
\addcontentsline{toc}{subsection}{Marked complex structures and the Thurston pullback}
Put $\T=\mathrm{Teich}(\widehat\CC,A^+)$. A point is a conformal structure together with a marking from the reference sphere, considered up to conformal changes and isotopy relative to the labels. Thus the same curve class is defined on every marked structure. For four labels normalized to $0,1,\infty,\lambda$, moving $\lambda$ around $0$ returns to the same configuration but changes the marking; Teichmüller space is the universal cover of $\CC\setminus\{0,1\}$ in this example.

We use $d_T=\tfrac12\inf\log K$, minimizing quasiconformal dilatation in the marked class.\footnote{The metric in \cite[\S2.5]{MPR} is twice ours; nonexpansion and contraction factors are unchanged.} With $s\ge2$ finite labels, $\dim_\CC\T=s-2$.

In the pullback diagram we choose the bottom marking first and solve for the top marking.
\[
\begin{array}{ccc}
(\RR^2,A)&\xrightarrow{\ \psi\ }&(\CC,\psi(A))\\[-2pt]
{\scriptstyle f}\downarrow\phantom{\scriptstyle f} &&\phantom{\scriptstyle h}\downarrow{\scriptstyle h}\\[-2pt]
(\RR^2,A)&\xrightarrow{\ \phi\ }&(\CC,\phi(A)).
\end{array}
\qquad h=\phi f\psi^{-1},\qquad \sigma[\phi]=[\psi].
\]
Stable parabolic type makes $h$ entire. Pulling back a quasiconformal target change preserves its dilatation bound on regular charts, so $\sigma$ is nonexpanding. A fixed point is equivalent to a realization in the prescribed marking \cite[Proposition~2.24]{MPR}; for transcendental maps, $\sigma$ is uniformly contracting on compact subsets \cite[Proposition~2.22(1)--(2)]{MPR}. Pullback for $f^k$ is $\sigma^k$, represented by compositions of one-step entire maps between successive marked planes.

Analytically, a target change $\Phi$ with Beltrami coefficient $\nu=\Phi_{\bar w}/\Phi_w$ pulls back through an entire representative $H$ to
\[
 \widetilde\nu(z)=\nu(H(z))\frac{\overline{H'(z)}}{H'(z)}
\]
off the isolated critical points, which do not affect the equation almost everywhere. Solving $\Psi_{\bar z}=\widetilde\nu\,\Psi_z$ gives a quasiconformal marking, unique up to an affine map, for which $\Phi H\Psi^{-1}$ is holomorphic \cite{AIM}, \cite[Propositions~2.6 and~2.21]{MPR}.

For a composition of one-step representatives, a disk avoiding the target marks again has inverse branches obtained by successive continuation through the factors. Indeed, by forward invariance an intermediate preimage over that disk cannot be a mark. Thus the finite singular values of the composition are among its target marks.

\subsection*{Pushforward and the cotangent map}
\addcontentsline{toc}{subsection}{Pushforward and the cotangent map}
The following statement supplies the analytic input for coverings with infinitely many sheets as well as for finite-degree maps. For a finite set $B\subset\CC$, write $Q(B^+)$ for the integrable holomorphic quadratic differentials on $\widehat\CC\setminus B^+$, where $B^+=B\cup\{\infty\}$; their poles at the marks are at most simple.

\begin{pushforwardlemma}[Pushforward]\label{l:pushforward}
Let $H:\CC\to\CC$ be nonconstant and entire, and let $A_s,A_t\subset\CC$ be finite sets such that $H(A_s)\subset A_t$ and
$H:\CC\setminus H^{-1}(A_t)\to\CC\setminus A_t$ is a covering. For every integrable measurable quadratic differential $\eta=u(z)\,\dd z^2$, the inverse-branch sum
\[
 p=H_*\eta,\qquad
 p(w)=\left(\sum_{H(z)=w}\frac{u(z)}{H'(z)^2}\right)\dd w^2
\]
converges absolutely for almost every $w$, defines an integrable measurable quadratic differential, and satisfies $\|p\|_1\le\|\eta\|_1$. For every measurable target set $W$,
\[
 \int_W\sum_{H(z)=w}\frac{|u(z)|}{|H'(z)|^2}\,\dd A(w)
       =\int_{H^{-1}(W)}|u(z)|\,\dd A(z).
\]
This area identity also holds for arbitrary measurable coefficients, with nonnegative integrals allowed to be infinite.

If in addition $\eta\in Q(A_s^+)$, the sum converges absolutely locally uniformly off $A_t$ and $p\in Q(A_t^+)$. For a nonzero $\eta\in Q(A_s^+)$, equality of the norms forces $\eta=cH^*p$ for a constant $c>0$. Consequently $\|H_*\eta\|_1<\|\eta\|_1$ for every nonzero $\eta\in Q(A_s^+)$ when $H$ is transcendental; when $H$ is a polynomial of degree $D$, equality for such $\eta$ forces $c=1/D$.
\end{pushforwardlemma}
\begin{proof}
Partition $\CC\setminus A_t$ into countably many disjoint Borel subsets of disks with closure disjoint from $A_t$. Each disk has at most countably many inverse branches with disjoint images. Change of variables on every sheet and Tonelli's theorem give the displayed area identity for measurable $u$. The exceptional sets and their inverse images have area zero. For integrable $u$ the nonnegative fiber sum has finite integral, so it is finite almost everywhere. This proves absolute convergence of the complex sum almost everywhere, and the triangle inequality gives the norm bound.

Now suppose $\eta\in Q(A_s^+)$. Over a disk $D_0$ with $\overline D_0\cap A_t=\varnothing$, write $\xi_j$ for the inverse branches. The holomorphic coefficients $b_j=(u\circ\xi_j)(\xi_j')^2$ satisfy
\[
 \sum_j\int_{D_0}|b_j|\,\dd A
       =\int_{H^{-1}(D_0)}|\eta|\le\|\eta\|_1.
\]
For each compact $K\subset D_0$, the area mean-value inequality gives
$\sup_K|b_j|\le C_{K,D_0}\int_{D_0}|b_j|\,\dd A$.
The Weierstrass test proves absolute locally uniform convergence and holomorphy of the sum.

An $L^1$ holomorphic coefficient on a punctured disk has at most a simple pole. Indeed, integrability gives radii $r_j\downarrow0$ with
$r_j^2\int_0^{2\pi}|v(r_je^{i\theta})|\,\dd\theta\to0$; otherwise the area integral would dominate a divergent multiple of $\int\dd r/r$. Cauchy's formula for the Laurent coefficients then makes every coefficient of $z^{-m}$, $m\ge2$, zero. Apply this at each target mark and, in the coordinate $1/w$, at infinity. Thus $p$ extends meromorphically to the sphere with at most simple marked poles.

If equality holds and $\eta\ne0$, then $p\ne0$. The triangle-inequality defect is nonnegative and has integral $\|\eta\|_1-\|p\|_1=0$. Thus, on almost every regular target fiber, every nonzero inverse-branch contribution has the same argument as $p$. Hence the meromorphic function $\eta/H^*p$ on the plane is nonnegative real almost everywhere where defined. On a regular disk where this quotient is holomorphic, continuity and the open mapping theorem make it constant; meromorphic continuation on the connected source gives the same constant everywhere. Thus $\eta=cH^*p$ with $c>0$, since $\eta\ne0$. For transcendental $H$, the great Picard theorem gives infinitely many preimages of every finite value except possibly one. The area identity therefore yields
\[
 \int_{\CC}|H^*p|=\int_{\CC}N_H(w)|p(w)|=\infty,
\]
where $N_H$ counts regular preimages, whereas $\eta=cH^*p$ is integrable. This contradiction proves strictness. For a polynomial of degree $D$, the identity $H_*H^*p=Dp$ on regular target disks gives $p=cDp$, so $c=1/D$.
\end{proof}

For a finitely marked sphere, $T_x^*\T=Q(x)$ with its $L^1$ norm; see \cite[\S1.1, p.~5 of the author version]{BCT}. Poles are allowed at every mark, including additional finite marks; the dimension depends on the number of marks. For the entire representative $h_x$ from the marked structure $\sigma x$ to $x$, the cotangent map is
\[
             (D\sigma_x)^*=(h_x)_*:Q(\sigma x)\longrightarrow Q(x).
\]
The marked entire setting, allowing additional marks, is given in \cite[Definition~2.14 and Proposition~2.21]{MPR}; holomorphy of pullback is stated just before Proposition~2.22 there. To prove the displayed identity also for infinitely many sheets, let $h=h_x$. Differentiating Beltrami pullback gives $\nu\mapsto h^*\nu$ on tangent representatives, where
\[
 (h^*\nu)(z)=\nu(h(z))\frac{\overline{h'(z)}}{h'(z)}
\]
off the isolated critical points. For a bounded target Beltrami differential $\nu$ and $q=u(z)\,\dd z^2\in Q(\sigma x)$, absolute integrability follows from
\[
 \int_{\CC}|u(z)|\,|\nu(h(z))|\,\dd A(z)
       \le\|\nu\|_\infty\|q\|_1<\infty.
\]
Partition the regular target into countably many Borel pieces in inverse-branch disks, as in the lemma, and change variables on each sheet. The Jacobian cancels the conjugate derivative factor. Absolute integrability permits the sums and integrals to be exchanged, giving
\[
 \int_{\CC}u(z)\nu(h(z))\frac{\overline{h'(z)}}{h'(z)}\,\dd A(z)
 =\int_{\CC}\nu(w)\sum_{h(z)=w}\frac{u(z)}{h'(z)^2}\,\dd A(w).
\]
This is precisely $\langle q,h^*\nu\rangle=\langle h_*q,\nu\rangle$. The lemma ensures that $h_*q\in Q(x)$: target marks contain all finite singular values and all images of source marks. If $\nu$ is infinitesimally trivial on the target, it annihilates $Q(x)$, so the identity shows that $h^*\nu$ annihilates every source cotangent vector. The formula therefore descends to tangent classes and identifies the transpose of $D\sigma_x$. Both cotangent spaces are those of the finitely marked spheres.

The finite-degree model is \cite[Proposition~3.2]{DH}, \cite[\S1.3]{BCT}; the exponential case is treated in \cite[\S3.2]{HSS}. For transcendental $h$, combine the cotangent identity $(D\sigma_x)^*=h_*$ with the strict norm inequality in Lemma~\ref{l:pushforward}. Compactness of the unit sphere of $Q(\sigma x)$ gives $\|(D\sigma_x)^*\|<1$ at each point. The $L^1$ cotangent norms vary continuously in local bundle coordinates: on a marked sphere, this follows from the rational coefficients and the uniform integrability of moving simple poles. Since $\sigma$ is holomorphic, $x\mapsto\|(D\sigma_x)^*\|$ is therefore continuous. The tangent norm is dual to this $L^1$ norm, and $d_T$ is its associated length metric \cite[\S\S3, 5]{Ant}. A compact Teichmüller geodesic segment consequently has a common derivative bound below one; integration along it gives $d_T(\sigma x,\sigma y)<d_T(x,y)$ for $x\ne y$. This proves the strict contraction asserted in \cite[Proposition~2.22(2)]{MPR}.

\emph{Scope of the extremal-length preliminaries.} The facts below, and Theorem~\ref{l:record}, apply to connected orientable surfaces $S$ of finite type without boundary, with ends treated as punctures. The preceding Thurston-pullback discussion concerns the marked sphere; we return to that setting after the theorem. When using geodesic laminations we assume $\chi(S)<0$; the proof treats the closed torus by passing to the once-punctured torus, and the spheres with at most three punctures separately.

\subsection*{Extremal length from a cylinder}\label{bg:length}
\addcontentsline{toc}{subsection}{Extremal length from a cylinder}
For an essential simple class $\alpha$, extremal length and its annular characterization are
\[
 \Ext_x\alpha=\sup_\rho\frac{L_\rho(\alpha)^2}{\area_\rho(x)},
 \qquad
 \frac1{\Ext_x\alpha}=\sup_{D:\,\mathrm{core}(D)=\alpha}\Mod D .
\]
Here $\rho|\dd z|$ ranges over conformal metrics of positive finite area, $L_\rho(\alpha)$ is the infimum of lengths in the marked class, and $D$ ranges over embedded annuli \cite[Definition~2.1]{FM}.

A Euclidean cylinder of circumference $h$ and height $w$ has modulus $w/h$ and core extremal length $h/w$. Indeed, for any density, Cauchy--Schwarz on the horizontal circles and integration in height give $wL_\rho^2\le h\area_\rho$; the Euclidean metric gives equality. The coordinate $\log z$ gives
\[
 \Mod\{r<|z|<R\}=\frac1{2\pi}\log(R/r).
\]
A connected degree-$d$ cover of an annulus $D$ has modulus $\Mod D/d$. This accounts for the reciprocal-degree weights below and explains why an annulus of large modulus makes its core short.

\subsection*{Quadratic differentials and measured foliations}\label{bg:fol}
\addcontentsline{toc}{subsection}{Quadratic differentials and measured foliations}
A quadratic differential $q=u(z)\,\dd z^2$ defines the coordinate-independent length element $|u|^{1/2}|\dd z|$, area measure $|q|=|u|\,\dd x\,\dd y$, and norm $\|q\|_1=\int|q|$. Away from zeros and poles, a primitive of a local square root gives $q=\dd\zeta^2$; transitions are $\zeta\mapsto\pm\zeta+c$. Integrability allows at most simple poles at punctures, since an order-$m$ pole has radial area integral comparable to $\int_0^\epsilon r^{1-m}\,\dd r$.

The vertical leaves $\Re\zeta=\mathrm{constant}$ have transverse measure $|\dd\Re\zeta|$; horizontal leaves use $|\dd\Im\zeta|$. A measured foliation is considered up to isotopy and Whitehead moves, which collapse or expand leaf segments between singularities. A cylinder of transverse width $a$ represents the weighted core $a\alpha$. For a test curve $\beta$, $i(U,\beta)$ is the infimum of its transverse measure in the class; for $U=a\alpha$ this is $a$ times the geometric intersection number. The pairing is jointly continuous \cite[Corollary~1.11, p.~468]{Rees}; see \cite[\S5, p.~5 of the arXiv version]{Ant} for punctured surfaces. We include $0$, with $i(0,U)=0$; continuity there also follows from Minsky's inequality below.

Weighted simple curves are dense in measured foliation space, and multiplying a transverse measure by $c>0$ multiplies extremal length by $c^2$. On each marked structure $x$, the Hubbard--Masur theorem gives the unique integrable holomorphic differential whose \emph{vertical} measured foliation represents a prescribed $U$. We define $\Ext_xU$ to be its area and set $\Ext_x0=0$. This agrees with the simple-curve definition and its continuous, degree-two homogeneous extension on surfaces of finite type, including punctured surfaces \cite[Theorem~5.1]{Ant}; the original closed-surface statement is in \cite[remark after Proposition~3, p.~34]{Ker}. Changing $x$ keeps $U$ as a marked topological object but generally changes its representing differential and its area. In particular, the differentials realizing $G$ on $g(t)$ and on $\sigma^k g(t)$ generally differ.

We use the standard measured foliation space, or, when $\chi(S)<0$, equivalently measured geodesic laminations on one fixed complete hyperbolic surface. Denote geometric intersection number by $i$. The square root of extremal length is continuous and homogeneous of degree one, its unit ball is compact, and the duality formula is
\begin{equation}\label{eq:duality}
 \sqrt{\Ext_x U}=\sup_{V\ne0}\frac{i(U,V)}{\sqrt{\Ext_x V}} .
\end{equation}
We use the topology of measured foliations determined by intersection with simple classes and the density of weighted simple curves \cite{FLP}; see also \cite[Section~2 of the arXiv version]{Wal}. Two measured foliations $U,V$ are called \emph{compatible} if $i(U,V)=0$. When $\chi(S)<0$, their geodesic representatives have no transverse crossings, although their supports may overlap. Their sum is defined by adding the transverse measures on the union of these supports, and intersection with a test class is additive on such sums. A pointwise bound on all intersection numbers gives the same bound on square roots of extremal lengths by the duality formula \eqref{eq:duality}.

Projective measured foliation space identifies positive rescalings. The normalization $\Ext_oU=1$ selects one representative of each nonzero class, giving a compact space in which short foliations can converge even as their curve classes vary \cite{FLP,Wal}.

The duality formula is due to Gardiner--Masur \cite{GM}. To recall the argument, let $q_U$ be the Hubbard--Masur differential for $U$ on $x$, allowing at most simple poles at the punctures \cite{HM}, \cite[\S5, Theorem~5.1]{Ant}. For any simple class $\beta$ its $q_U$-length dominates its transverse intersection with $U$. Testing extremal length in the metric $|q_U|^{1/2}$ therefore gives $i(U,\beta)^2\le\|q_U\|_1\Ext_x\beta$. Density gives the inequality $i(U,V)^2\le\Ext_xU\,\Ext_xV$ for every pair of measured foliations. This is Minsky's inequality \cite[Lemma~5.1]{Min}. The horizontal foliation of $q_U$ has intersection $\|q_U\|_1$ with $U$ and extremal length $\|q_U\|_1=\Ext_xU$, giving equality in the supremum when $U\ne0$. We shall also use Kerckhoff's formula \cite[Theorem~4]{Ker} in its finite-type formulation \cite[Theorem~2.2 and the following remark]{FM}, together with its distortion consequence:
\begin{equation}\label{eq:kerckhoff}
 e^{2d_T(x,y)}=\max_{U\ne0}\frac{\Ext_y U}{\Ext_x U},
 \qquad
 e^{-2d_T(x,y)}\Ext_xU\le\Ext_y U\le e^{2d_T(x,y)}\Ext_xU .
\end{equation}
The ratio depends only on the projective class of $U$. The factor two in \eqref{eq:kerckhoff}, and hence in the record-time normalization, comes from $d_T=\tfrac12\log K$.

At a \emph{record time} $n$, the distance $R=d_T(o,\sigma^n o)$ dominates all earlier distances. Kerckhoff's formula supplies a foliation of extremal length $1$ at $o$ and $e^{-2R}$ at $\sigma^n o$. The proof retains compact sets of such short foliations, since a single extremizer may change with $n$. Their compatible limits and finitely many test curves will produce one $G$ and one bound for every starting point and depth.

\begin{samepage}
\begin{theorem}\label{l:record}
Let $S$ be a connected orientable surface of finite type without boundary, and let $\T=\mathrm{Teich}(S)$ be its Teichmüller space, with the ends of $S$ treated as punctures.
Suppose that $\sigma:\T\to\T$ is $1$-Lipschitz for the Teichmüller metric and has an unbounded forward orbit.
Then there are a nonzero measured foliation $G$ on $S$ and a constant $C\ge1$ with
\begin{equation}\label{eq:displacement}
                  \Ext_{\sigma^k x}G\le C\Ext_xG
                  \qquad(x\in\T,\ k\ge0).
\end{equation}
\end{theorem}
\end{samepage}
\begin{proof}
A sphere with at most three punctures has a one-point Teichmüller space, so every orbit is bounded.
Marking one point on a closed torus identifies its Teichmüller space isometrically with that of the once-punctured torus, since a translation can align the marked points without changing the dilatation of a comparison map.
On both surfaces a measured foliation is specified by a slope and its transverse weight. In the flat model $\CC/(\mathbb Z+\tau\mathbb Z)$ with $\Im\tau>0$, the class of a primitive vector $(a,b)$ has extremal length $|a+b\tau|^2/\Im\tau$ on either surface \cite[pp.~566--568]{Min99}.
Homogeneity and density give the same identification for all measured foliations. Thus the estimate for the once-punctured torus gives the estimate on the closed torus with the same $G$ and $C$.

We may now assume that $S$ has negative Euler characteristic and realize its measured foliations as geodesic laminations with transverse measures in a complete finite-area hyperbolic metric.

\textbf{1. Record times and normalization.}
Fix $o\in\T$ and put $x_j=\sigma^j o$. If $y$ has an unbounded orbit, then $d_T(\sigma^j o,\sigma^j y)\le d_T(o,y)$ for every $j$, so the orbit of $o$ is also unbounded.

Choose $n_i\uparrow\infty$ with $R_i=d_T(o,x_{n_i})\to\infty$ and $d_T(o,x_j)\le R_i$ for $0\le j\le n_i$. For fixed $k$ and $n_i\ge k$ put
\[
     K_{i,k}=\{U:\Ext_{x_{n_i-k}}U\le e^{-2R_i}\},
     \qquad \mathcal U_o=\{U:\Ext_o U\le1\}.
\]
The closed unit ball $\mathcal U_o$ consists of all measured foliations whose extremal length on the fixed marked surface $o$ is at most one. If $F$ is scaled so that $\Ext_o F=1$, the ball contains exactly the multiples $tF$ with $0\le t\le1$ along that ray, since $\Ext_o(tF)=t^2$.
Extremal-length distortion and the record property put every $K_{i,k}$ in the same compact metrizable ball $\mathcal U_o$. Each is nonempty and compact. Compactness of the Hausdorff space of nonempty compact subsets of $\mathcal U_o$ and a diagonal extraction give limits $K_k$ for all nonnegative integers $k$. Kerckhoff's formula in the order $\Ext_o/\Ext_{x_{n_i}}$ has a maximizer on compact projective foliation space. Normalize it at $o$ to length $\Ext_o U_i=1$. Then $\Ext_{x_{n_i}}U_i=e^{-2R_i}$. Any cluster point of the $U_i$ lies in $K_0$ and still has length one at $o$.

\textbf{2. A single finite family of test curves.}
Let $D=d_T(o,\sigma o)$. Nonexpansion gives $d_T(x_{n_i-k},x_{n_i})\le kD$. Distortion and quadratic homogeneity, first before taking limits and then by continuity of rescaling, give
\begin{equation}\label{eq:incl}
                  e^{-kD}K_0\subset K_k\subset e^{kD}K_0 .
\end{equation}
For example, if $U_i\in K_{i,0}$ tends to $U$, then
$\Ext_{x_{n_i-k}}(e^{-kD}U_i)\le e^{-2kD}e^{2kD}\Ext_{x_{n_i}}U_i\le e^{-2R_i}$.
Thus $e^{-kD}U\in K_k$. Starting with approximants in $K_{i,k}$ and interchanging the two times proves $e^{-kD}K_k\subset K_0$. Although these inclusions identify a common linear span, their factors grow with $k$. The bound on coefficients below uses the separate fact $K_k\subset\mathcal U_o$.

For $U,V\in K_0$ take approximants from the same $K_{i,0}$. Minsky's inequality there bounds their intersection by $e^{-2R_i}$, so $i(U,V)=0$. Any finite family in $K_0$ has noncrossing geodesic supports. Its nonnegative sums are measured laminations, and intersection with every simple class is additive on these sums.

Write $I_U(\alpha)=i(U,\alpha)$ on the set of essential simple classes. The functions $I_U$ for $U\in K_0$ span a finite-dimensional real vector space. If $I_{F_1},\ldots,I_{F_m}$ are linearly independent, the map $(c_1,\ldots,c_m)\mapsto\sum c_\ell F_\ell$ for $c_\ell>0$ is continuous and injective, by intersection additivity. Measured foliation space is homeomorphic by the finite-type Hubbard--Masur theorem at $o$ \cite[Theorem~5.1]{Ant} to the space of integrable holomorphic quadratic differentials on $S$. If $S$ has genus $g$ and $n$ punctures, this space has complex dimension $3g-3+n$ by Riemann--Roch, allowing at most simple poles at the punctures. Invariance of dimension therefore gives $m\le 6g-6+2n$.

Choose a basis $I_{F_1},\ldots,I_{F_b}$ from $K_0$ and put $G=\sum_{\ell=1}^b F_\ell$. Step~1 gives a nonzero element of $K_0$ with extremal length one at $o$, so $b\ge1$ and $G\ne0$. There are simple classes $\beta_1,\ldots,\beta_b$ for which the matrix $J=(i(F_\ell,\beta_j))_{j,\ell}$ is invertible. Otherwise the evaluation rows would span a proper subspace of $\RR^b$ and a nonzero vector orthogonal to all of them would give a dependence. For $U\in K_k$ its intersection function is in this span by the second inclusion in \eqref{eq:incl}. Its coordinate vector is
\[
                 c(U)=J^{-1}(i(U,\beta_j))_{j=1}^b .
\]
Since $U\in\mathcal U_o$, each entry on the right before applying $J^{-1}$ is at most $\sqrt{\Ext_o\beta_j}$ by Minsky. Thus $|c_\ell(U)|\le B$ for every $k$, with the single finite constant $B=\|J^{-1}\|_\infty\max_j\sqrt{\Ext_o\beta_j}$, where $\|\cdot\|_\infty$ is the induced matrix norm. Nonnegativity of $I_{F_\ell}$ implies $i(U,\alpha)\le B\,i(G,\alpha)$ for every simple $\alpha$. By homogeneity and density this holds against every measured foliation. Extremal-length duality \eqref{eq:duality} gives $\Ext_xU\le B^2\Ext_xG$ for all $x$.

\textbf{3. The comparison at arbitrary $x$ and $k$.}
For every fixed $y$ and $k$, quadratic homogeneity and Kerckhoff's formula give the exact identity
\[
     \max_{U\in K_{i,k}}\Ext_y U
        =\exp\!\left(2d_T(y,x_{n_i-k})-2R_i\right).
\]
Indeed, put $z=x_{n_i-k}$. For every $U\in K_{i,k}$, Kerckhoff's formula \eqref{eq:kerckhoff} and the defining bound $\Ext_z U\le e^{-2R_i}$ give
\[
    \Ext_y U\le e^{2d_T(z,y)}\Ext_z U\le e^{2d_T(z,y)-2R_i}.
\]
To attain equality, choose $F$ maximising $\Ext_y F/\Ext_z F$ and set $U=tF$ with $t=e^{-R_i}/\sqrt{\Ext_z F}$. Both extremal lengths scale by $t^2$, so $U$ still attains the maximal ratio $e^{2d_T(z,y)}$, while $\Ext_z U=e^{-2R_i}$. Thus $U\in K_{i,k}$ attains both inequalities, and $d_T(z,y)=d_T(y,z)$ gives the claimed identity.

Compare $y=\sigma^k x$ with index zero to $y=x$ with index $k$. Since $\sigma^k x_{n_i-k}=x_{n_i}$, nonexpansion gives $d_T(\sigma^k x,x_{n_i})\le d_T(x,x_{n_i-k})$. Applying the identity above with the same $R_i$ on both sides gives $\max_{U\in K_{i,0}}\Ext_{\sigma^k x}U\le\max_{U\in K_{i,k}}\Ext_xU$. As $K_{i,0}\to K_0$ and $K_{i,k}\to K_k$ in Hausdorff distance, continuity of extremal length gives convergence of both maxima and hence the first inequality below.
\[
      \max_{U\in K_0}\Ext_{\sigma^k x}U
         \le\max_{U\in K_k}\Ext_xU\le B^2\Ext_xG .
\]
Indeed, if compact sets $L_i$ converge in Hausdorff distance to $L$ inside $\mathcal U_o$, each maximizer of a fixed continuous function on $L$ has approximants from $L_i$. Conversely any sequence of maximizers on $L_i$ has a subsequence tending to a point of $L$. These two comparisons give convergence of the maxima. We apply this separately for each fixed pair $(x,k)$ to the same previously chosen sets.

Put $y=\sigma^k x$. Each $F_\ell$ belongs to $K_0$, so the preceding estimate gives
\[
   \Ext_y F_\ell\le\max_{U\in K_0}\Ext_y U\le B^2\Ext_xG
   \qquad(1\le\ell\le b).
\]
To pass from the individual foliations to their sum $G$, fix a test foliation $V$ with $\Ext_y V=1$.
The $F_\ell$ are compatible, so their transverse measures add. Thus $i(G,V)=\sum_{\ell=1}^b i(F_\ell,V)$, first for weighted simple curves and then for every measured foliation by density and continuity.
By duality \eqref{eq:duality}, $\sqrt{\Ext_y F_\ell}$ is the supremum of $i(F_\ell,W)$ over all $W$ with $\Ext_y W=1$. In particular it bounds $i(F_\ell,V)$ for this same $V$. Hence
\[
   i(G,V)=\sum_{\ell=1}^b i(F_\ell,V)
      \le\sum_{\ell=1}^b\sqrt{\Ext_y F_\ell}
      \le bB\sqrt{\Ext_xG}.
\]
Taking the supremum over all such $V$ gives $\sqrt{\Ext_yG}\le bB\sqrt{\Ext_xG}$ by the same duality formula. Squaring proves \eqref{eq:displacement} with $C=(bB)^2$. The family of $b$ foliations was fixed along the record subsequence, and Step~2 chose one $B$ for every $K_k$. Thus the same $C$ works for all $x$ and $k$.
\end{proof}

We return to the marked sphere $S=\widehat{\CC}\setminus A^+$. For the Thurston pullback of a transcendental map, the unbounded-orbit hypothesis of Theorem~\ref{l:record} holds whenever the pullback has no fixed point. Indeed, \cite[Proposition~2.22(2)]{MPR} gives $d_T(\sigma x,\sigma y)<d_T(x,y)$ for every $x\ne y$ by taking the compact set $\{x,y\}$. A bounded orbit would have compact closure $\Omega$ by properness of $\T$, with $\sigma(\Omega)\subset\Omega$. The continuous function $x\mapsto d_T(x,\sigma x)$ would attain its minimum on $\Omega$ at some $y$. If $\sigma y\ne y$, then $d_T(\sigma y,\sigma^2y)<d_T(y,\sigma y)$ contradicts minimality. Thus $y$ would be a fixed point.

\begin{figure}[ht]\centering
\begin{tikzpicture}[scale=0.85]
\foreach \r in {1.1,2.0,3.0}{\draw[pine!45,dashed] (0,0) circle (\r);}
\fill (0,0) circle (1.6pt) node[below left]{$o$};
\draw[rust,thick,-{Stealth}] plot[smooth,tension=0.6] coordinates {(0,0) (0.5,0.3) (0.8,-0.5) (1.1,0.0) (0.6,0.9) (1.4,1.1) (2.0,0.0) (1.5,-1.0) (1.8,-1.6) (2.5,-0.9) (3.0,0.0)};
\foreach \p in {(0.5,0.3),(0.8,-0.5),(0.6,0.9),(1.4,1.1),(1.5,-1.0),(1.8,-1.6),(2.5,-0.9)}{\fill[rust] \p circle (1.2pt);}
\foreach \p/\l in {(1.1,0.0)/x_{n_1},(2.0,0.0)/x_{n_2},(3.0,0.0)/x_{n_3}}{\fill[pine] \p circle (2pt); \node[pine,font=\scriptsize,above right] at \p {$\l$};}
\node[font=\scriptsize,pine,align=left] at (-3.4,2.3) {records: $d_T(o,x_j)\le R_i$\\ for $j\le n_i$};
\draw[pine,-{Stealth}] (3.3,2.1) node[right,font=\scriptsize,align=left]{at $x_{n_i}$ a foliation $U_i$ with\\ $\Ext_oU_i=1$ has $\Ext_{x_{n_i}}U_i=e^{-2R_i}$} -- (3.08,0.25);
\begin{scope}[shift={(8.6,-0.2)}]
\draw[thick,pine] (0,0) circle (1.5);
\fill[mist] (0,0) -- (40:1.5) arc (40:70:1.5) -- cycle;
\draw[very thick,rust] (40:1.5) arc (40:70:1.5);
\node[font=\scriptsize,rust] at (55:2.02) {$\mathbb P(K_0\setminus\{0\})$};
\node[font=\scriptsize,align=center] at (0,-1.95) {projective classes of nonzero\\ limits with pairwise $i=0$};

\end{scope}
\end{tikzpicture}
\caption{Record times (left) give a Hausdorff limit $K_0$ of normalized short-foliation sets. The arc (right) schematically represents its nonzero projective classes, whose pairwise intersections vanish. The construction of $G$ and the uniform bound over all $K_k$ are in Steps~2--3.}
\label{fig:record}
\end{figure}
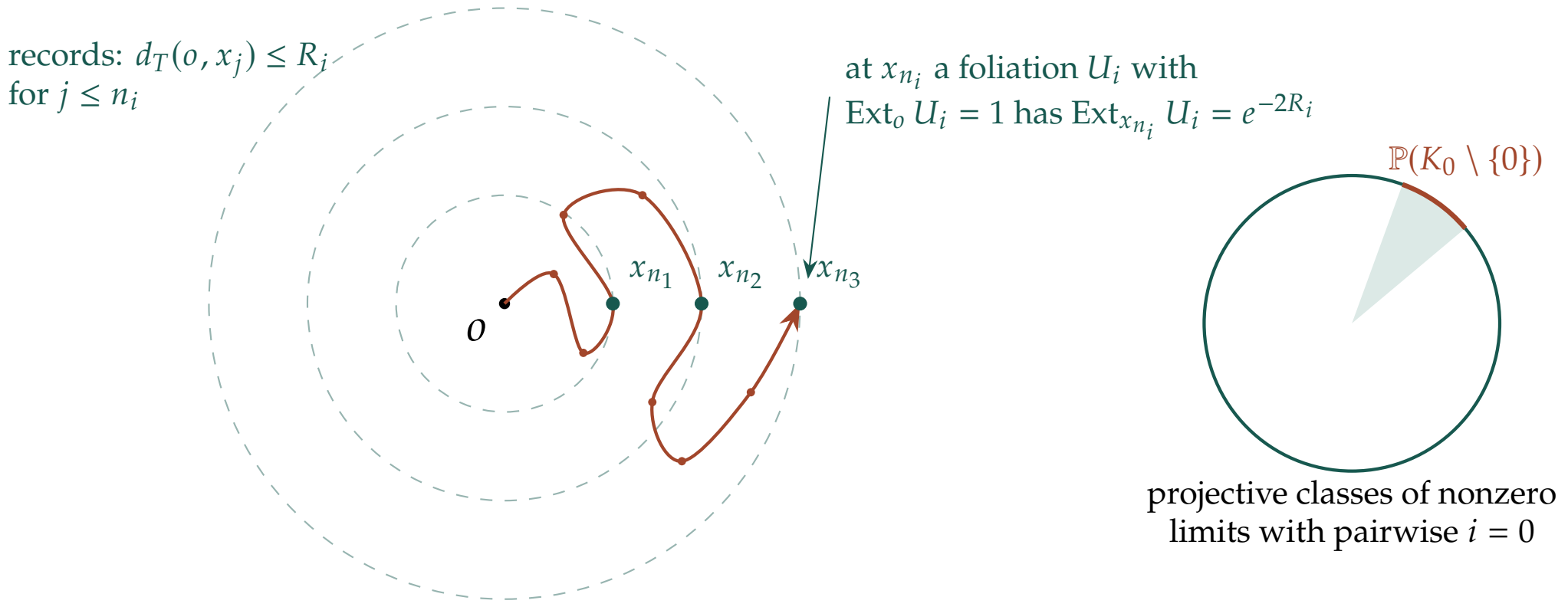

\emph{A boundary interpretation.} A horofunction is a limit of normalized distance profiles $x\mapsto d(x,z)-d(o,z)$ as $z$ escapes to infinity. The following remark relates the record-time construction to this boundary language.

\begin{remark}\label{r:wolff}
For the record subsequence, put $h_k(x)=\lim_i(d_T(x,x_{n_i-k})-R_i)$. Step~3 and Hausdorff convergence give existence of the limit and $\max_{K_k}\Ext_x=e^{2h_k(x)}$. This is the horofunction normalization for $k=0$, and differs from it by a constant for fixed $k$. Liu--Su \cite{LS} identify the horofunction and Gardiner--Masur compactifications; see also \cite[Theorem~1.3 and \S6]{Wal} and \cite{Miy}. Nonexpansion gives $h_0\circ\sigma^k\le h_k$, as in Beardon--Karlsson's Wolff--Denjoy approach \cite{Bea,Kar}, applied to Teichmüller space in \cite{Kar14}. The theorem adds $h_k\le h_0+\log(bB)$ uniformly in $k$, by representing the compatible limit sets through one foliation $G$.

For a mapping class $\varphi$, $\Ext_{\varphi^kx}G=\Ext_x(\varphi^{-k}G)$. A pseudo-Anosov gives the expanding lamination, and a Dehn twist $T_\gamma$ gives $G=\gamma$, both with $C=1$. On $\mathrm{Teich}(S_{0,4})\cong\mathbb H$, the theorem is a coarse analogue of Wolff's invariant-horoball lemma: $\log\Ext_xG$ increases by at most $\log C$ under any iterate. For a Levy cycle of $p$ curves, the degree-one annular lifts and duality give \eqref{eq:displacement} for their sum with $C=p^2$.
\end{remark}

Write $\Sigma$ for the punctured target sphere representing $o$, with its marking. Rescale $G$ so its Hubbard--Masur differential $q$ on $\Sigma$ has area one. Let $g(t)$ be its Teichmüller ray. In regular coordinates write $\zeta_t=e^tx+ie^{-t}y$ and $q_t=\dd\zeta_t^2$; the affine map $\phi_t:x+iy\mapsto\zeta_t$ from the initial flat surface to the time-$t$ surface is the \emph{ray marking}. The marked vertical measure of $q_t$ is $e^tG$, and the area stays one. Thus $e^{-2t}q_t$ realizes $G$ and the area identity gives $\Ext_{g(t)}G=e^{-2t}$ \cite[Theorem~5.1]{Ant}. In particular
\begin{equation}\label{eq:ray}
                  \Ext_{\sigma^k g(t)}G\le C e^{-2t}
                       \qquad(t\ge0,\ k\ge0).
\end{equation}

\emph{Standing input.} Sections~\ref{sec:flat}--\ref{sec:persistence} use only a normalized nonzero $G$ satisfying \eqref{eq:displacement}, its Hubbard--Masur differential and its ray, together with the stated hypotheses on the marked entire maps. No further property of the record-time construction is needed. Nonexpansion of $\sigma$ supplies the energy monotonicity in Section~\ref{sec:transport}; absence of Levy cycles is invoked where stated.

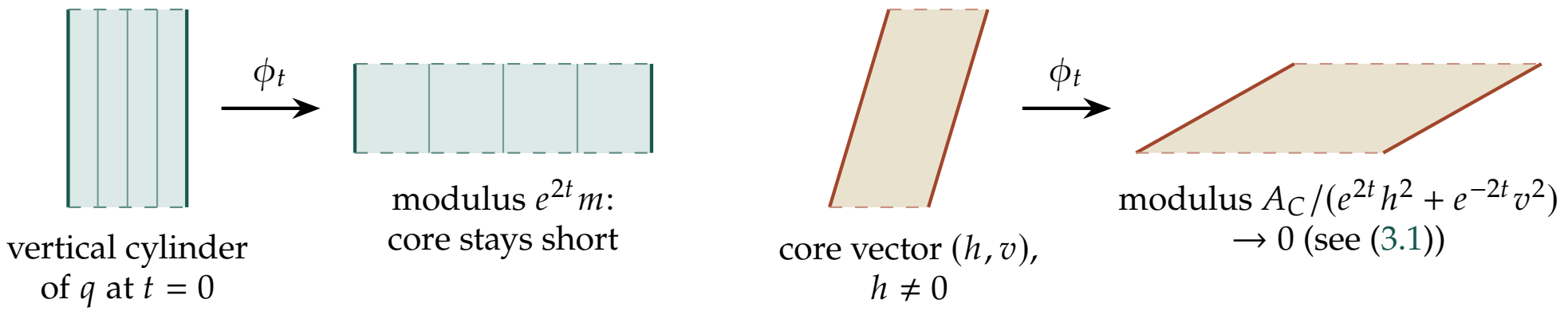
\begin{figure}[ht]\centering
\begin{tikzpicture}[scale=0.82]
\begin{scope}
\fill[mist] (0,0) rectangle (1.2,2.0); \draw[pine,thick] (0,0)--(0,2) (1.2,0)--(1.2,2);
\draw[pine!70,dashed] (0,0)--(1.2,0) (0,2)--(1.2,2);
\foreach \x in {0.3,0.6,0.9}{\draw[pine!60] (\x,0)--(\x,2);}
\node[font=\scriptsize,below,align=center] at (0.6,-0.1) {vertical cylinder\\ of $q$ at $t=0$};
\end{scope}
\draw[-{Stealth},thick] (1.55,1)--(2.55,1) node[midway,above,font=\scriptsize]{$\phi_t$};
\begin{scope}[shift={(2.9,0.55)}]
\fill[mist] (0,0) rectangle (3.0,0.9); \draw[pine,thick] (0,0)--(0,0.9) (3.0,0)--(3.0,0.9);
\draw[pine!70,dashed] (0,0)--(3,0) (0,0.9)--(3,0.9);
\foreach \x in {0.75,1.5,2.25}{\draw[pine!60] (\x,0)--(\x,0.9);}
\node[font=\scriptsize,below,align=center] at (1.5,-0.1) {modulus $e^{2t}m$:\\ core stays short};
\end{scope}
\begin{scope}[shift={(7.7,0)}]
\fill[sand] (0,0)--(1.0,0)--(1.6,2.0)--(0.6,2.0)--cycle;
\draw[rust,thick] (0,0)--(0.6,2.0) (1.0,0)--(1.6,2.0);
\draw[rust!60,dashed] (0,0)--(1,0) (0.6,2)--(1.6,2);
\node[font=\scriptsize,below,align=center] at (0.8,-0.1) {core vector $(h,v)$,\\ $h\neq0$};
\end{scope}
\draw[-{Stealth},thick] (9.65,1)--(10.55,1) node[midway,above,font=\scriptsize]{$\phi_t$};
\begin{scope}[shift={(10.8,0.55)}]
\fill[sand] (0,0)--(2.5,0)--(4.1,0.9)--(1.6,0.9)--cycle;
\draw[rust,thick] (0,0)--(1.6,0.9) (2.5,0)--(4.1,0.9);
\draw[rust!60,dashed] (0,0)--(2.5,0) (1.6,0.9)--(4.1,0.9);
\node[font=\scriptsize,below,align=center] at (2.05,-0.1) {modulus $A_C/(e^{2t}h^2+e^{-2t}v^2)$\\ $\to0$ (see \eqref{eq:modulus})};
\end{scope}
\end{tikzpicture}
\caption{The ray stretches natural coordinates by $(x,y)\mapsto(e^tx,e^{-t}y)$. Vertical cylinders gain modulus by $e^{2t}$; every nonvertical cylinder eventually loses modulus. These are the alternatives used in Sections~\ref{sec:closed} and~\ref{sec:persistence}.}
\label{fig:stretch}
\end{figure}

\section{Flat cylinders}\label{sec:flat}
We prove two facts used in Section~\ref{sec:persistence}: at most $r=s-2$ maximal flat cylinders on the $(s+1)$-punctured sphere have modulus greater than one, and the modulus of each nonvertical cylinder tends to zero along the ray. Ray markings identify each cylinder with one fixed open subset of the initial surface, allowing its mass for a fixed measure to be followed as its modulus changes.

We work with a nonzero integrable quadratic differential $\Phi$ on $\widehat\CC$ with marked set $A^+$. Marked regular points are removed along with the other punctures. A \emph{flat cylinder} is an embedded open annulus foliated by parallel regular closed geodesics of equal length. For the flat geometry used here see Strebel \cite{Str}.

\subsection*{The flat metric at a zero and along a cylinder}
\addcontentsline{toc}{subsection}{The flat metric at a zero and along a cylinder}
For $q=z^m\dd z^2$ ($m\ge-1$), a natural coordinate is $w=2z^{(m+2)/2}/(m+2)$, giving cone angle $(m+2)\pi$. Zeros have negative curvature, regular marks have none, and simple poles have curvature $\pi$. Gauss--Bonnet requires $2\pi$ inside a disk bounded by a regular closed geodesic, so such a disk contains at least two simple poles. In particular, cylinder cores are essential.

\emph{Crossing geodesics intersect essentially.} Two regular closed geodesics that cross transversely meet in finitely many points. An empty bigon between them would be a disk of nonpositive interior curvature whose two exterior turns sum to less than $2\pi$, which Gauss--Bonnet excludes. By the bigon criterion \cite[Proposition~1.7]{FMar} their classes have positive intersection. Tangent closed geodesics coincide.

\emph{Maximal cylinders.} Two distinct regular closed geodesics in the same class are therefore disjoint, and they cobound an unmarked annulus. By Gauss--Bonnet this annulus contains no zero, so it is a flat cylinder foliated by parallel closed geodesics. Consequently the union of all regular closed geodesics in a given class is a single open flat cylinder, the \emph{maximal cylinder} of that class, and it contains every cylinder with that core \cite{Str}.

\emph{At most $r$ cylinders above modulus one.} Let $\alpha$ and $\beta$ be the core classes of two distinct maximal cylinders of moduli $M,M'>1$. By Minsky's inequality and the annulus bound,
\[
 i(\alpha,\beta)^2\le\Ext(\alpha)\Ext(\beta)\le\frac{1}{MM'}<1.
\]
Since $i(\alpha,\beta)$ is a nonnegative integer, it is zero. If the cylinders overlapped in different directions, the closed geodesics through a point of overlap would cross transversely and hence have positive intersection. If they overlapped in the same direction, they would share a closed geodesic and thus be the same maximal cylinder.

We can therefore choose one core curve from each cylinder to obtain disjoint essential curves in distinct isotopy classes. Cutting a sphere with $s+1$ punctures along $p$ such curves gives $p+1$ components. Each has at least three ends because the curves are essential and pairwise nonisotopic. Counting each puncture once and each cut twice gives $s+1+2p$ ends in total, so $3(p+1)\le s+1+2p$ and hence $p\le s-2=r$. In any fixed direction the same bound holds at every modulus, since parallel closed geodesics cannot cross.

\emph{Moduli along the ray.} Let $\Phi=q$ and let $\phi_t$ be the affine ray marking of Section~\ref{sec:records}. It maps straight closed geodesics of $q$ to straight closed geodesics of $q_t$, and its inverse does the same, so each maximal $q$-cylinder $C$ is carried onto a maximal $q_t$-cylinder. We regard $C$ as one fixed open subset of the initial surface; its area $A_C$ is constant. Its core has translational holonomy, since its tangent returns to itself; let $(h,v)$ be its translation vector in natural coordinates of $q$, up to sign. At time $t$ the circumference is the length of $(e^th,e^{-t}v)$, and the width is $A_C$ divided by the circumference, even when the image of an initially perpendicular segment becomes oblique. Therefore
\begin{equation}\label{eq:modulus}
 m_C(t)=\frac{A_C}{e^{2t}h^2+e^{-2t}v^2}.
\end{equation}
If $h=0$, that is, if $C$ is a \emph{vertical} cylinder, $m_C(t)$ grows like $e^{2t}$. If $h\ne0$, then $m_C$ is continuous and tends to zero as $t\to\infty$.

\section{Vertical cylinders and the Thurston matrix}\label{sec:closed}
In this section we show that if the measured foliation $G$ from Theorem~\ref{l:record} contains a vertical cylinder of positive width, then $f$ has a Levy cycle. The argument leads to the following incompatible bounds if there is no Levy cycle:
\[
  0<c_\gamma\ \le\ \frac{L_k}{2}\ \le\ (M_\Lambda^k m)_\gamma
        \ \longrightarrow\ 0\qquad(k\to\infty).
\]
Here $\gamma$ is a cylinder core of $G$, $m$ records the initial cylinder moduli, and $\Lambda$ is the invariant multicurve formed from their essential lifts. The quantity $L_k/2$ is the limiting renormalized reciprocal extremal length of $\gamma$ after retaining four suitable source marks at depth $k$. Theorem~\ref{l:record} supplies the positive lower bound $c_\gamma$, Proposition~\ref{p:quartet} supplies the upper bound by the matrix contribution, and the Matrix Lemma below supplies its decay. We first let ray time $t\to\infty$ with $k$ fixed to obtain $L_k$; only afterwards do we let $k\to\infty$. Proposition~\ref{l:cylinder-levy} carries out these steps and gives~\eqref{eq:matrix}.

For an entire representative $F$, let $A_s,A_t$ be its finite source and target marking sets. We have $F(A_s)\subset A_t$, and $A_t$ contains every finite singular value. Thus $F:\CC\setminus F^{-1}(A_t)\to\CC\setminus A_t$ is a covering. For a Jordan curve $\beta$ avoiding $A_t$, restriction gives a covering of a circle. Each connected component maps onto $\beta$ by path lifting. A compact component has finite fibers and is therefore a finite covering circle; a noncompact component is a line \cite[\S1.3, Theorem~1.38]{Hat}. Since $\infty$ is marked, the bounded disk of an essential plane circle contains between $2$ and $s-1$ finite labels. This asymmetry between the bounded and unbounded sides lets us control an infinite-degree map with a finite matrix.

We use the bounded side of each simple closed plane curve. Representatives may be chosen smooth away from the marks. If $\alpha$ is a compact component of $F^{-1}(\beta)$, the Jordan--Sch\"onflies theorem gives a bounded open disk $W$ with frontier $\alpha$; see \cite[Theorems~2.12 and~3.1, pp.~122--123]{Tho}. Let $D$ be the bounded target disk and $d$ the degree of $F|_\alpha$. With $\alpha$ positively oriented, the argument principle gives
\[
 \sum_{z\in W:\ F(z)=w}\deg_z F
       =\operatorname{wind}(F\circ\alpha,w),\qquad w\notin\beta.
\]
The winding number is zero outside $\overline D$, so no point of $W$ maps there. Openness also excludes an interior point mapping to $\beta$. For $w\in D$ the winding number is either $d$ or $-d$; since the nonempty open image $F(W)$ meets $D$, the nonnegative preimage count forces the sign to be positive. Thus every point of $D$ has $d$ preimages counted with multiplicity. The inverse image in $W$ of a compact subset of $D$ is closed in the compact set $\overline W$ and misses its boundary. Hence $F:W\to D$ is proper of degree $d$. Distinct compact components over the same curve bound disjoint disks, since nesting would put a point over $\beta$ inside a disk whose image lies in $D$. Every source mark in $W$ consequently maps to a target mark in $D$. The same conclusions hold for topological representatives by marked equivalence. Only the map on the selected disk has finite degree; the entire map may be transcendental.

If the target disk contains no finite mark, its proper preimage disks contain none by forward invariance. If it contains just one mark $p$, remove the full finite nonempty fiber over $p$ from a proper preimage disk $W$. If that fiber has $m$ points, the resulting connected surface is a disk minus $m$ points and is a finite cover of the once-punctured target disk, since all finite singular values are marked. Its fundamental group is therefore isomorphic to a finite-index subgroup of $\mathbb Z$ \cite[\S1.3, Theorem~1.38]{Hat}. A disk minus $m$ points has free fundamental group of rank $m$, so $m=1$. Every source mark in $W$ lies in this fiber, and hence there is at most one. At the other extreme, suppose the bounded disk $\Delta$ contains all finite marks. The exterior $U=\CC\setminus\overline\Delta$ is a punctured disk and contains no singular value of an entire representative $F$. If a component over $U$ had finite degree $d$, classification of connected covers would give an injective holomorphic parametrization $\varphi:\mathbb D^*\to\CC$ of that component in which $F\circ\varphi$ tends to infinity at the puncture (in a punctured-disk coordinate on $U$ the covering is $\zeta\mapsto\zeta^d$). Injectivity excludes an essential singularity of $\varphi$ by Picard's theorem. A finite limit contradicts continuity of $F$. Thus $\varphi$ has a simple pole, its image contains a neighborhood of infinity, and $F(z)\to\infty$ there. This would make $F$ a polynomial. For transcendental $F$ every component over $U$ is consequently a universal cover. A closed lift of $\partial\Delta$ would have a finite-degree exterior collar. Inclusion of the target collar in $U$ induces an isomorphism on fundamental groups, so the containing component would have that same finite degree. This excludes a closed lift.

It follows that an essential closed preimage component has essential image. In particular the successive images of an essential compact $k$-step lift are essential closed components at every intermediate depth. For a polynomial the exterior containing infinity is connected when the target disk contains all finite singular values, since every component over the exterior disk on the sphere maps onto it and contains the unique preimage of infinity. The exterior curve therefore has a single lift enclosing all finite marks, peripheral relative to $A^+$.

Choose one-step representatives $F_j:X_j\to X_{j-1}$, where $X_j$ is the marked plane representing $\sigma^jo$. Starting with $m_0$, choose the source marking $m_j$ in the pullback construction so that $m_{j-1}\circ f=F_j\circ m_j$ exactly (\mpr, Proposition~2.21). Target isotopies lift over $\CC\setminus A=\widehat\CC\setminus A^+$ in reference coordinates \cite[Proposition~1.30]{Hat}. These identities identify compact components and their degrees at every depth. More explicitly, put $H_j=F_1\circ\cdots\circ F_j$. If $b$ is a compact component of $H_{j+1}^{-1}(c)$ for a target curve $c$ avoiding the marks, then $F_{j+1}(b)$ lies in one component $K$ of $H_j^{-1}(c)$. It is relatively open in $K$ by covering charts and closed by compactness. Hence it equals $K$, which is itself a compact covering circle. The restricted covering degrees multiply. The $i$ factors adjacent to the source of a $j$-step lift are $F_{j-i+1}\circ\cdots\circ F_j$. Through the markings this composition of successive entire representatives corresponds to the topological iterate $f^i$ on the reference plane.

\subsection*{Reading a Thurston matrix}
\addcontentsline{toc}{subsection}{Reading a Thurston matrix}
A \emph{multicurve} is a finite set of pairwise disjoint, nonisotopic essential curve classes. A degree-$d$ closed lift of an annulus of modulus $m_\beta$ contributes $m_\beta/d$ to its source class $\alpha$. The entry $M_{\alpha\beta}$ sums these reciprocal degrees: columns correspond to target classes, rows to source classes. Products of degrees along lifted chains explain matrix powers. For $z^2-1$ in the introduction, the matrix is $[1/2]$, so its $k$th power is $[2^{-k}]$. In general, the spectral radius $\rho(M)$ is the largest modulus of an eigenvalue; $\rho(M)<1$ implies $M^k\to0$.

Draw $\alpha\to\beta$ when $M_{\alpha\beta}>0$. A \emph{block} consists of a maximal group of vertices mutually reachable by directed paths; it is \emph{recurrent} if it contains a cycle. A \emph{return of length $p$} is a path of $p$ arrows back to its starting class. The lemma shows that each recurrent block is a single weighted cycle. Its key input is that a periodic marked point has a unique periodic predecessor, which disjoint lifted disks cannot both contain. This restricts internal arrows of a recurrent block; a column of the full matrix may still have several nonzero entries.

The plane hypothesis in the Matrix Lemma is essential. Shishikura and Tan Lei give a Thurston map of the sphere with a Thurston obstruction and no Levy cycle \cite{ShiTan}. For the Levy criterion for topological polynomials, see Hubbard \cite[\S10.3]{Hub} and the marked formulation obtained from \cite[\S1.3 and Proposition~B.1]{BMW}.

\begin{matrixlemma}[Matrix lemma]\label{l:matrix}
Let $(f,A)$ be a marked entire Thurston map of the plane and let $\Lambda$ be a finite multicurve relative to $A^+=A\cup\{\infty\}$. Form the matrix $M=(M_{\alpha\beta})_{\alpha,\beta\in\Lambda}$ with entries
\[
  M_{\alpha\beta}=\sum_{\delta}\frac{1}{\deg(f:\delta\to\beta)},
\]
where $\delta$ ranges over the closed components of the preimage of a representative of $\beta$ that are isotopic to $\alpha$ relative to $A^+$. If $f$ has no Levy cycle relative to $A^+$, then $\rho(M)<1$.
\end{matrixlemma}
\begin{proof}
For a fixed target representative the compact preimage disks are disjoint. Each essential one contains at least two finite marks, so there are at most $\lfloor s/2\rfloor$ essential compact components and the entries are finite sums. This count concerns actual components; there may be infinitely many inessential compact components. An isotopy of target curves relative to the marks lifts in the cover over their complement, and the inverse isotopy gives the inverse correspondence. The trace of each lifted compact circle is the continuous image of a circle times $[0,1]$, hence compact. The curves continue to avoid the source marks, whose images are target marks. Thus this is a bijection of compact components preserving covering degree and marked source class, and in particular preserving essentiality.

Draw $\alpha\to\beta$ when $M_{\alpha\beta}>0$. Let $S_\alpha$ be the finite labels in the bounded side and $P_\alpha$ the periodic labels among them. Both are well defined on the marked class. An arrow gives $f(S_\alpha)\subset S_\beta$. In a strongly connected component containing an edge on a cycle, each vertex lies on a cycle. If a cycle through $\alpha$ has length $p$, then $f^p(S_\alpha)\subset S_\alpha$. Iterating $f^p$ on this nonempty finite set eventually reaches a cycle, so $S_\alpha$ contains a periodic marked point $b$ of $f$. Its least period need not equal the return length $p$: for example, a one-step return of a disk class can preserve a pair of marks that $f$ exchanges. Thus the periodic point $b$ is one of the marked labels in the disk. Let $\mathcal P\subset A$ be the periodic labels. The restriction $f|_{\mathcal P}$ is a permutation since $\mathcal P$ is a disjoint union of finite cycles. Thus along internal arrows the cardinality of $P_\alpha$ cannot decrease, and by strong connectivity it is constant. In particular $f(P_\alpha)=P_\beta$ along each internal arrow.

Fix a representative of $\beta$ and choose $b\in P_\beta$. Since $f(P_\alpha)=P_\beta$ for every internal arrow $\alpha\to\beta$, each corresponding preimage disk must contain the same point $a$, the unique periodic predecessor of $b$ on its orbit. Distinct preimage disks are disjoint, so only one such lift is possible, even when counting lifts in the same class. Thus every vertex has at most one incoming internal edge, and strong connectivity gives at least one. Hence each has exactly one. There are as many internal edges as vertices, and every vertex has at least one outgoing internal edge, so each has exactly one. Following these edges gives a single cycle through the whole component by strong connectivity. Its weights are $1/d_1,\ldots,1/d_p$. The $p$th power of its matrix is $(\prod d_i)^{-1}$ times the identity, so the spectral radius of this block is $(\prod d_i)^{-1/p}$. Some $d_i>1$ because otherwise these multicurve classes form a Levy cycle. Acyclic diagonal blocks have radius zero. Ordering the strongly connected components makes $M$ block triangular and proves the claim.
\end{proof}

\subsection{A four-point bound}\label{sec:fourlimit}
At fixed depth $k$, a degree-$d$ lift of a vertical cylinder of initial modulus $m$ has modulus $e^{2t}m/d$. Proposition~\ref{p:quartet} gives the corresponding \emph{upper} bound after retaining four source marks: finite lifts in the selected class account for the leading reciprocal extremal length. Its proof uses finite-energy one-forms on a torus, even though the pulled-back flat area near the point over infinity can be infinite. Proposition~\ref{l:cylinder-levy} then combines this bound with Theorem~\ref{l:record} and lets $k\to\infty$.

\subsection*{Four marked points and the torus}\label{bg:four}
\addcontentsline{toc}{subsection}{Four marked points and the torus}
Join four distinct sphere marks in two pairs by disjoint arcs. Cut two copies along these arcs and glue opposite banks across the copies. The result is the double cover branched at the four marks: a loop about an endpoint changes sheets, and the local projection is $z\mapsto z^2$. Riemann--Hurwitz gives Euler characteristic $2\cdot2-4=0$, hence a torus. Conversely, the quotient of $\CC/(\mathbb Z+\tau\mathbb Z)$ by $z\mapsto-z$ is a sphere with the four half-lattice points as branch points. Each essential curve downstairs has two degree-one lifts whose common unoriented primitive slope specifies its class \cite[\S4]{CFPP}.

For instance take the square torus and the horizontal direction $(1,0)$. The integer linear transformation $(x,y)\mapsto(x,2x+y)$ commutes with the involution and fixes each half-lattice point modulo the lattice. On the quotient it preserves every label and carries the chosen class to the slope $(1,2)$. These two slopes are distinct. Thus the same partition of four labels can occur for different curve classes. When Section~\ref{sec:excluding} chooses four marks to detect a member of a multicurve, the prior disjointness of the family is essential.

For a rectangular torus of circumference one and height $T$, the quotient strip between branch heights $0$ and $T/2$ has modulus $T/2$. The form $\dd x$ has horizontal period one and energy $T$. Thus the torus energy is twice the reciprocal extremal length downstairs, as proved in Step~1 of Proposition~\ref{p:quartet}. Each counted lift also gives two torus annuli, explaining the cancelling factors of two in \eqref{eq:matrix}.

Fix a depth $k$, its entire representative $h$ at $g(0)$, and four source marks $B$ including infinity. Let $X_{t,B}$ be the sphere representing $\sigma^k g(t)$ with only $B$ retained, put $\Ext_{B,t}=\Ext_{X_{t,B}}$, and let $\mathcal E_t$ be its branched double cover. An essential class $\beta$ downstairs determines the common unoriented slope $[a]$ of its two primitive torus lifts. Choose either orientation $a$, a homology class upstairs, and define
\[
 H_a(\mathcal E_t)=\min_{\substack{\dd\eta=0\\ \int_a\eta=1}}
       \int_{\mathcal E_t}|\eta|^2\,\dd A,
\]
where $\eta$ ranges over real $L^2$ one-forms and the period $\int_a\eta$ is understood on cohomology. In flat coordinates $\mathcal E_t=\CC/(\mathbb Z+\tau_t\mathbb Z)$ in which $a$ is horizontal and has length one, this minimum equals the height $\operatorname{Im}\tau_t$ of the fundamental parallelogram.

\emph{Surfaces and coordinates.} The differential $q$ lives on the time-zero \emph{target}; its deformation $q_t$ defines the vertical cylinders along the target ray $g(t)$. The torus $\mathcal E_t$ instead covers the four-marked \emph{source} $X_{t,B}$. It is a topological branched double cover equipped with the lifted conformal structure, not the canonical double cover on which $q$ becomes a square. There are three coordinate systems in the proof:
\begin{itemize}
\item $(X,Y)$ are uniformizing Euclidean coordinates on $\mathcal E_t$, normalized so that the chosen class $a$ is horizontal of length one. They give $H_a(\mathcal E_t)=\Im\tau_t$.
\item $(x,y)$ are regular natural coordinates of $(h\circ\pi)^*q$ on $\mathcal E_0$, and similarly at time $t$. They define the projections $P_x,P_y$ and the stretch $(x,y)\mapsto(e^t x,e^{-t}y)$. Their vertical direction need not be the $Y$-direction of the uniformizing coordinates.
\item $(u,v)\bmod\mathbb Z^2$ are real marked lattice coordinates, with the $u$-circle representing $a$. They are used to define periods and Fourier series, without a conformality requirement.
\end{itemize}
Here $\pi$ denotes the projection to the source sphere. The composition $h\circ\pi$ is used away from the point over $\infty$; no holomorphic extension of a transcendental $h$ across that point is assumed.

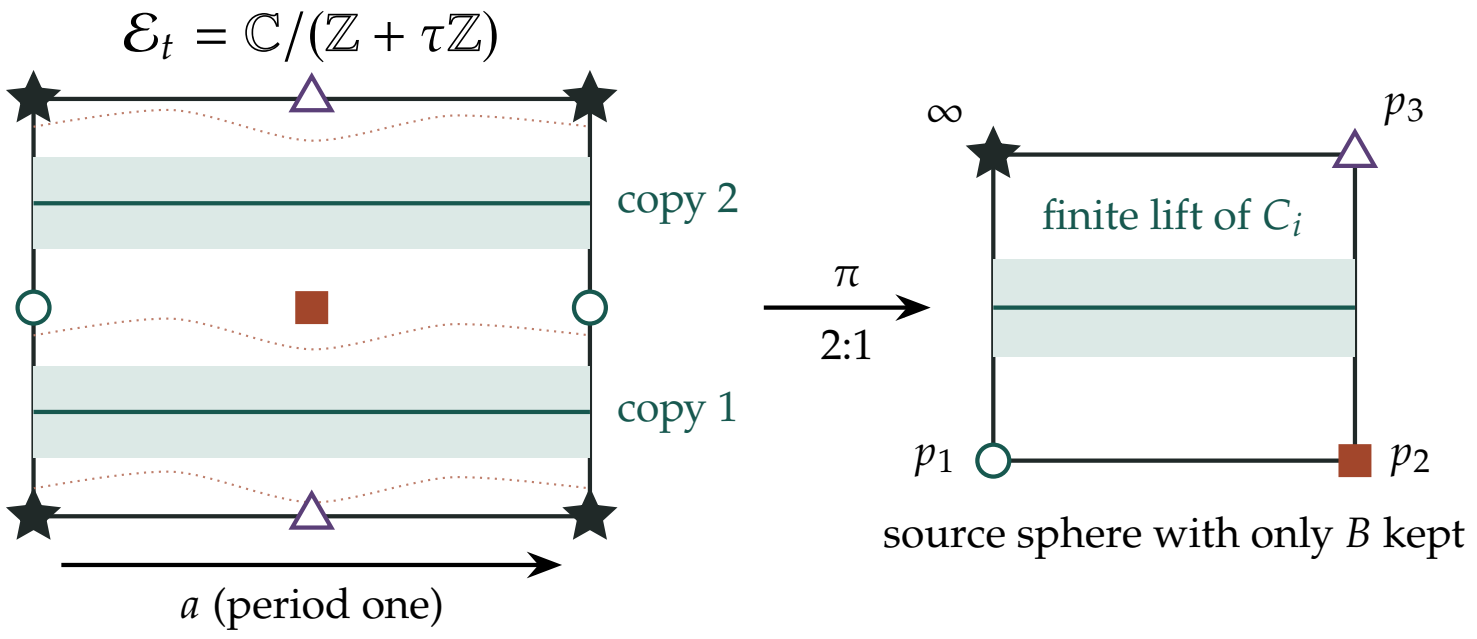
\begin{figure}[ht]\centering
\begin{tikzpicture}[scale=0.92,
  branch mark/.style={inner sep=0pt,outer sep=0pt,line width=0.8pt},
  p1 mark/.style={branch mark,circle,minimum size=6pt,draw=pine,fill=white},
  p2 mark/.style={branch mark,rectangle,minimum size=5.4pt,draw=rust,fill=rust},
  p3 mark/.style={branch mark,regular polygon,regular polygon sides=3,minimum size=8.5pt,draw=plum,fill=white},
  infinity mark/.style={branch mark,star,star points=5,star point ratio=2.1,minimum size=9pt,draw=ink,fill=ink}]
\begin{scope}
\draw[thick,ink] (0,0) rectangle (4,3);
\fill[mist] (0,0.42) rectangle (4,1.08);
\fill[mist] (0,1.92) rectangle (4,2.58);
\draw[pine,thick] (0,0.75)--(4,0.75); \draw[pine,thick] (0,2.25)--(4,2.25);
\foreach \y in {0.2,1.3,2.8}{\draw[rust!70,densely dotted] plot[smooth] coordinates {(0,\y) (1,\y+0.12) (2,\y-0.1) (3,\y+0.08) (4,\y)};}
\foreach \p in {(0,1.5),(4,1.5)}{\node[p1 mark] at \p {};}
\node[p2 mark] at (2,1.5) {};
\foreach \p in {(2,0),(2,3)}{\node[p3 mark] at \p {};}
\foreach \p in {(0,0),(4,0),(0,3),(4,3)}{\node[infinity mark] at \p {};}
\node[font=\scriptsize,pine,right] at (4.05,0.75) {copy 1};
\node[font=\scriptsize,pine,right] at (4.05,2.25) {copy 2};
\draw[-{Stealth},thick] (0.2,-0.35)--(3.8,-0.35) node[midway,below,font=\scriptsize]{$a$ (period one)};
\node[font=\small] at (2,3.45) {$\mathcal E_t=\mathbb C/(\mathbb Z+\tau\mathbb Z)$};
\end{scope}
\draw[-{Stealth},thick] (5.25,1.5)--(6.45,1.5) node[midway,above,font=\scriptsize]{$\pi$} node[midway,below,font=\scriptsize]{$2{:}1$};
\begin{scope}[shift={(6.9,0.4)}]
\draw[thick,ink] (0,0) rectangle (2.6,2.2);
\fill[mist] (0,0.75) rectangle (2.6,1.45);
\draw[pine,thick] (0,1.1)--(2.6,1.1);
\node[p1 mark,label={[font=\scriptsize]left:$p_1$}] at (0,0) {};
\node[p2 mark,label={[font=\scriptsize]right:$p_2$}] at (2.6,0) {};
\node[infinity mark,label={[font=\scriptsize]above left:$\infty$}] at (0,2.2) {};
\node[p3 mark,label={[font=\scriptsize]above right:$p_3$}] at (2.6,2.2) {};
\node[font=\scriptsize,align=center] at (1.3,-0.55) {source sphere with only $B$ kept};
\node[font=\scriptsize,pine,align=center] at (1.3,1.75) {finite lift of $C_i$};
\end{scope}
\end{tikzpicture}
\caption{The torus double cover branched at the four retained marks; only the front face of the pillowcase is shown. Opposite sides of the rectangle are identified, so its four stars represent one point over $\infty$. The upper-left and lower-right quarters map to the front face, the other two to the back; each interior point of a face has two preimages exchanged by $z\mapsto-z$. The shaded annuli are the two torus lifts of a counted cylinder, giving the factor $2$ in \eqref{eq:quartet}. Their cores are drawn horizontally in the uniformizing torus coordinates; their vertical direction refers to the pulled-back differential. The dotted leaves represent infinite-degree lifts approaching the point over $\infty$. Step~4 proves that the dual form $\theta$ vanishes almost everywhere there.}
\label{fig:torus}
\end{figure}

\subsection*{What the energy minimization measures}
\addcontentsline{toc}{subsection}{What the energy minimization measures}
For $\eta=u\,\dd x+v\,\dd y$, Dirichlet energy is $\int(u^2+v^2)\,\dd x\,\dd y$ and is conformally invariant. The constraint $\int_a\eta=1$ fixes the period while allowing exact changes of the form. The resulting minimum measures reciprocal extremal length:
\[
 \frac1{\Ext_{B,t}(\beta)}=\frac{H_a(\mathcal E_t)}2.
\]
Step~1 below proves the identity. Along the ray, division by $e^{2t}$ makes the horizontal part of the energy inexpensive; the remaining vertical energy detects the counted cylinders. If the limiting energy is positive, the period functional is bounded in the vertical $L^2$ norm; taking the Hodge star of its Riesz representative produces the dual form $\theta$ in Step~3. The main issue is to show that nonperiodic leaves contribute nothing (Step~4). We first justify periods and approximation for closed $L^2$ forms, whose coefficients need not have traces on a prescribed curve.

\begin{analyticprep}[Closed forms and periods on a torus]\label{l:analyticprep}
On a compact marked conformal torus, every distributionally closed real $L^2$ one-form has the representation
\[
        \eta=c_1\,\dd u+c_2\,\dd v+\dd F,
        \qquad F\in W^{1,2}(\RR^2/\mathbb Z^2).
\]
The constants $c_1,c_2$ are its marked periods. They are continuous linear functionals of $\eta$, and there are smooth closed forms with exactly these periods converging to $\eta$ in $L^2$. Marking-preserving quasiconformal homeomorphisms identify these spaces of closed forms and preserve the periods.
\end{analyticprep}
\begin{proof}
A period of a smooth closed form is its integral around a marked loop. For an $L^2$ form, changing the coefficients on that loop leaves the form unchanged in $L^2$, since the loop has area zero. To define the period, use marked linear coordinates $(u,v)\bmod\mathbb Z^2$ on the uniformizing flat torus, with the $u$-circle representing $a$. Write $\eta=\eta_u\,\dd u+\eta_v\,\dd v$. For a smooth closed form, the horizontal circles at different heights are homologous, so $\int_0^1\eta_u(u,v)\,\dd u$ has the same value for every $v$. Averaging over $0\le v\le1$ shows that the period on $a$ equals the mean $c_1$ of $\eta_u$ on the unit square. Likewise the period on the $v$-circle equals the mean $c_2$ of $\eta_v$. Both means are defined for $L^2$ coefficients and depend continuously on them by Cauchy--Schwarz.

To see that these means account for both periods at this regularity, we construct a periodic potential for the remainder $\eta-c_1\dd u-c_2\dd v$. Closedness in the sense of distributions means $\partial_u\eta_v=\partial_v\eta_u$. On a Fourier mode $e^{2\pi i(m_1u+m_2v)}$, differentiation in $u$ or $v$ multiplies the coefficient by $2\pi i m_1$ or $2\pi i m_2$, respectively. Thus $m_1\widehat{\eta}_v(m)=m_2\widehat{\eta}_u(m)$. For each $m\ne0$ the coefficient vector is parallel to $m$, so it is the gradient coefficient $2\pi i m\,\widehat F(m)$ of a scalar function. Explicitly, take
\[
 \widehat F(m)=\frac{m_1\widehat{\eta}_u(m)+m_2\widehat{\eta}_v(m)}
                       {2\pi i\,|m|^2}\quad(m\ne0),\qquad \widehat F(0)=0 .
\]
By Parseval, the two $L^2$ coefficients of $\eta$ have square-summable Fourier coefficients. Hence $\sum_{m\ne0}|m|^2|\widehat F(m)|^2<\infty$. This gives a mean-zero periodic function $F\in W^{1,2}$, meaning that $F$ and its first weak derivatives are square-integrable. Matching Fourier coefficients gives $\eta=c_1\dd u+c_2\dd v+\dd F$.

Now let $F_N$ be the finite Fourier sum of $F$ over $|m|\le N$. The smooth closed forms $\eta_N=c_1\dd u+c_2\dd v+\dd F_N$ converge to $\eta$ in $L^2$. Each $\dd F_N$ integrates to zero around every closed loop because $F_N$ is single-valued on the torus. Every $\eta_N$ has exactly the same two periods $c_1,c_2$. Convergence of their energies therefore extends any lower bound proved for smooth closed forms with these periods to $\eta$.

A marking-preserving quasiconformal homeomorphism $\Psi$ between compact tori also preserves these periods and identifies their closed $L^2$ forms. For a $K$-quasiconformal map one has
\[
             \int|\Psi^*\eta|^2\le K\int|\eta|^2,
\]
and composition with $\Psi$ preserves $W^{1,2}$ and obeys the weak chain rule \cite[Theorem~2.7 and its proof, $p=n=2$]{OP}. These local critical-exponent estimates apply in finitely many charts on the compact tori. For a smooth periodic potential the differential pulls back to that of its composition. For the constant coordinate forms, a lift of $\Psi$ to the plane is the identity in marked lattice coordinates plus a periodic $W^{1,2}$ map, so pullback changes these forms only by exact terms. The period-preserving approximation above, applied in both directions, proves the claim for all closed $L^2$ forms.
\end{proof}

For a finite-area Strebel ray $X_t$, Farb and Masur show that if $\gamma$ is the core of a vertical cylinder of initial modulus $m$, then $e^{-2t}/\Ext_{X_t}(\gamma)\to m$ as $t\to\infty$ \cite[Lemma~2.3(1)]{FM}. After pullback by $h$, a degree-$d$ lift of such a cylinder has modulus $e^{2t}m/d$. The new point in Proposition~\ref{p:quartet} is the upper bound showing that, even when $h$ has infinite degree, the sum of $m/d$ over the finite lifts in the chosen four-point class gives the entire coefficient of $e^{2t}$ in its reciprocal extremal length.

Douady and Hubbard~\cite[Theorem~7.1(c)]{DH} estimate the reciprocal hyperbolic length of a curve after filling punctures by summing the reciprocal lengths of the short geodesics that become homotopic to it. Their additive error depends on the number of punctures filled. In infinite degree the full preimage of the target marks can be infinite, so their estimate gives no finite error bound when only the four marks in $B$ are retained.

\begin{samepage}
\begin{prop}\label{p:quartet}
Fix a depth $k$ and its entire representative $h$ at $g(0)$. Let $B$ be any set of four distinct source marks containing $\infty$, and fix an essential curve class $\beta$ on the reference sphere marked by $B$. Through the markings, $\beta$ denotes the same class on every $X_{t,B}$. Let $\mathcal E_t$ be its double cover branched over $B$, and $[a]=\{a,-a\}$ the common unoriented slope of the two torus lifts of $\beta$. Let $C_i$ be the maximal positive-width vertical cylinders of $q$ on the target with its marks removed. Write $\gamma_i$ for their cores and $m_i$ for their initial moduli. For each $i$ at most one finite lift $\alpha_i$ of $\gamma_i$ under $h$ satisfies $[\alpha_i]_B=\beta$, where $[\alpha_i]_B$ denotes its class after forgetting the marks outside $B$. Put $d_i=\deg(h:\alpha_i\to\gamma_i)$ for such a lift. Then the limit below exists and
\[
       \lim_{t\to\infty}\frac{e^{-2t}}{\Ext_{B,t}(\beta)}
           \ \le\ \sum_{i:\,[\alpha_i]_B=\beta}\frac{m_i}{d_i}.
\]
Equivalently, for the fixed $k,B,\beta$, the limit $L=\lim_{t\to\infty}e^{-2t}H_a(\mathcal E_t)$ exists and
\begin{equation}\label{eq:quartet}
 L\ \le\ 2\sum_{i:\ \alpha_i\ \mathrm{exists}}\frac{m_i}{\deg(h:\alpha_i\to\gamma_i)} .
\end{equation}
A finite lift is counted precisely when it separates the four retained marks into two pairs and its two torus lifts have slope $[a]$. A $2+2$ partition alone does not specify the class $\beta$.
\end{prop}
\end{samepage}
The requirement $\infty\in B$ is part of the statement: it is used in the assertion that at most one finite lift is counted. The bounded disk of an essential lift contains exactly two of the three finite marks in $B$, whereas distinct compact lifts of the same target core bound disjoint disks. Two such disks cannot each contain two of those three marks. If all four retained marks were finite, this counting argument would no longer exclude two disjoint disks containing two retained marks each.

In fact equality holds in \eqref{eq:quartet}: at time $t$ each counted lift supplies two disjoint torus annuli of modulus $e^{2t}m_i/d$, and a closed form with period one on $a$ has energy at least the modulus on each. Only the upper bound is used.
\begin{proof}
\textbf{1. Periods and counted lifts.}
In conformal coordinates $\mathcal{E}_t=\CC/(\mathbb Z+\tau\mathbb Z)$ with the origin at a branch point of the quotient and $a$ the unit horizontal generator, $H_a(\mathcal{E}_t)=\Im\tau$. Put $Y_\tau=\Im\tau$. For a smooth closed form $\eta=\eta_1\,\dd X+\eta_2\,\dd Y$ of period one, $\int_0^1 \eta_1(X,Y)\,\dd X=1$ at each height, using periodicity in $X$. By Cauchy--Schwarz the energy per unit of height is at least one. Integrating over the height of the torus gives
\[
 E(\eta)
 =\int_0^{Y_\tau}
   \underbrace{\left(\int_0^1(\eta_1^2+\eta_2^2)\,\dd X\right)}_{\ge 1\ \text{at each height }Y}
   \,\dd Y
 \;\ge\;\int_0^{Y_\tau}1\,\dd Y
 =Y_\tau.
\]
The same bound holds for closed $L^2$ forms $\eta=\eta_1\,\dd X+\eta_2\,\dd Y$ of period one by the period-preserving approximation in the analytic preparation. The constant form $\dd X$ attains this value.

For the corresponding length bound downstairs, put $\Lambda=\mathbb Z+\tau\mathbb Z$. The map from $\CC\setminus\tfrac12\Lambda$ to the four-punctured sphere is the regular covering with deck group $z\mapsto\pm z+\lambda$, $\lambda\in\Lambda$. A horizontal unit segment avoiding the half-lattice projects to a core in class $\beta$ by the slope correspondence \cite[\S4, p.~220]{CFPP}. Lift a free homotopy from any representative of $\beta$ to that core \cite[Proposition~1.30]{Hat}. The deck transformation relating the endpoints remains constant along the lifted homotopy. A terminal lift is a deck translate of the unit segment, in either orientation, and conjugating $z\mapsto z+1$ in this group gives $z\mapsto z\pm1$. Thus every such representative has a lift with endpoint displacement $\pm1$, and length at least one. The lifted homotopy fixes the marked curve class throughout this argument.

The quotient by the elliptic involution has a cylinder of circumference one and height $\Im\tau/2$, and total flat area $Y_\tau/2$. To verify $1/\Ext_{B,t}(\beta)=H_a(\mathcal{E}_t)/2$, note that an open horizontal strip between successive branch heights projects injectively to the quotient apart from the period-one translation, which gives $1/\Ext_{B,t}(\beta)\ge Y_\tau/2$. The preceding length bound in the descended flat metric gives $\Ext_{B,t}(\beta)\ge 2/Y_\tau$ in that same class.

Suppose a plane lift $\alpha$ of $\gamma_i$ separates the four marks of $B$ into two pairs. Each complementary disk contains two branch points and lifts to an annulus. Their common boundary on the torus consists of two parallel essential curves, each mapping with degree one to $\alpha$. If $\alpha$ covers $\gamma_i$ with degree $d$, each torus copy of the lifted cylinder has initial modulus $m_i/d$. A different pairing of the four marks gives a different slope, and even the same pairing can occur for several slopes. Thus \eqref{eq:quartet} counts a $2+2$ lift precisely when its torus cores have slope $[a]$. For a $1+3$ or $0+4$ split, the disk with at most one branch point lifts to one or two disks, so every torus lift of $\alpha$ is null-homologous.

\textbf{2. The fixed Hilbert space.}
Write $\mathcal{E}=\mathcal{E}_0$ and $\pi:\mathcal{E}\to\widehat{\CC}$ for its projection to the source sphere. Remove the four ramification points of $\pi$, which includes the point over $\infty$. Also remove the inverse images under $h\circ\pi$ of the finite target marks and zeros of $q$. The resulting exceptional set $Z$ is closed and countable because the inverse images are locally finite away from the point above infinity. On $\mathcal{E}\setminus Z$ the pullback of $q$ has regular coordinates $x+iy$ with linear transitions $\pm I$. Let $P_x$ and $P_y$ project the one-form $\eta$ onto its horizontal and vertical components, respectively. If $\eta=u\,\dd x+v\,\dd y$, then $P_x\eta=u\,\dd x$ and $P_y\eta=v\,\dd y$. The Hilbert space consists of real $L^2$ one-forms that are distributionally closed on the whole compact torus $\mathcal E$. We use $\mathcal E\setminus Z$ for regular flat coordinates and integration, since $Z$ has area zero.

To compare energies, let $g_{\mathrm{sm}}$ be a smooth conformal metric on $\mathcal E$. On a regular chart the pulled-back flat metric is $\rho^2g_{\mathrm{sm}}$, and
\[
       |\eta|_{\rho^2g_{\mathrm{sm}}}^2\,\dd A_{\rho^2g_{\mathrm{sm}}}
           =|\eta|_{g_{\mathrm{sm}}}^2\,\dd A_{g_{\mathrm{sm}}}.
\]
Thus the Dirichlet norm is the ordinary $L^2$ norm of a one-form on a compact conformal torus, even if the pulled-back flat metric has infinite area near the point over infinity. The analytic preparation defines its periods and supplies smooth approximation without any assumption on that flat area.

Let $\phi_t$ be the affine marking on the target ray and solve on the source for a quasiconformal map $\psi_t$ with the complex structure pulled back by $h$ from $\phi_t$. Then $h_t=\phi_t h\psi_t^{-1}$ is entire. The source marking $\psi_t$ takes the four branch labels to their time-$t$ positions. It lifts to a quasiconformal map $\Phi:\mathcal{E}\to \mathcal{E}_t$: on the punctured spheres it preserves the double cover defined by odd monodromy at each branch label, and the lift extends across the four points. In corresponding regular natural coordinates $(x,y)$ and $(x_t,y_t)$, either lift has the form
\[
 \Phi(x,y)=(\varepsilon e^t x+c_0,\varepsilon e^{-t}y+d_0),
                    \qquad \varepsilon\in\{1,-1\}.
\]
This follows by taking coordinates through $h\circ\pi$ and its time-$t$ counterpart. The same sign occurs in both components because the regular chart transitions have derivative $\pm I$. Its Jacobian in these flat coordinates is one, so $\Phi$ preserves the pulled-back flat area, even if that area is infinite.

By Lemma~\ref{l:analyticprep}, pullback by $\Phi$ identifies closed $L^2$ forms on the whole compact tori and preserves their marked periods. For $\zeta=\zeta_1\,\dd x_t+\zeta_2\,\dd y_t$ on $\mathcal E_t$, write $\eta=\Phi^*\zeta$. Then
\[
 \eta=\varepsilon e^t(\zeta_1\circ\Phi)\,\dd x
       +\varepsilon e^{-t}(\zeta_2\circ\Phi)\,\dd y,
 \qquad
 E_{\mathcal E_t}(\zeta)
       =e^{-2t}\|P_x\eta\|_2^2+e^{2t}\|P_y\eta\|_2^2.
\]
The energy identity follows by squaring these coefficients and changing variables with Jacobian one; conformal invariance permits the chartwise computation of the Dirichlet energy. Thus it is the energy on $\mathcal E_t$ expressed through the pulled-back form on the fixed torus. Writing $\ell(\eta)=\int_a\eta$ on cohomology gives
\[
 e^{-2t}H_a(\mathcal{E}_t)=\inf_{\ell(\eta)=1}
        \bigl(\|P_y\eta\|_2^2+e^{-4t}\|P_x\eta\|_2^2\bigr)
       \ \longrightarrow\
 L:=\inf_{\ell(\eta)=1}\|P_y\eta\|_2^2 .
\]
Each left side is at least $L$. Given $\epsilon>0$, choose one admissible form $\eta_\epsilon$ with $\|P_y\eta_\epsilon\|_2^2\le L+\epsilon$. Keeping this form fixed makes $e^{-4t}\|P_x\eta_\epsilon\|_2^2$ tend to zero. The limsup is then at most $L+\epsilon$ for every $\epsilon>0$, with no compactness assumption on the minimizing forms at time $t$.

\textbf{3. A dual form $\boldsymbol{\theta}$ with integral periods.}
If $L=0$ there is nothing to prove. For $L>0$ homogeneity gives $|\ell(\eta)|\le L^{-1/2}\|P_y\eta\|_2$. Consequently $P_y\eta\mapsto\ell(\eta)$ is well defined on the projections of closed forms and has norm exactly $L^{-1/2}$ by the definition of the infimum. Let $b$ be its Riesz representative in the closure of those projections. That is, $b$ is the form in this closure satisfying
\[
 \langle P_y\eta,b\rangle_{L^2}=\int_a\eta
         \quad\text{for every closed }\eta.
\]
It is vertical almost everywhere. The Hodge star gives a horizontal form $\theta=*b$ with
\[
          \|\theta\|_2^2=1/L,\qquad
          \int_{\mathcal{E}}\eta\wedge\theta=\ell(\eta)
          \quad\text{for every closed }L^2\text{ form }\eta.
\]
Testing $\eta=\dd\varphi$ for smooth $\varphi$ gives distributional closedness of $\theta$ on \emph{all} of $\mathcal{E}$. Testing $\dd u,\dd v$ in an oriented lattice basis starting with $a$ gives constant mode $\dd v$ for $\theta$. Indeed, for the constant mode $c_1\dd u+c_2\dd v$, pairing in the first argument with $\dd u$ gives $c_2=\ell(\dd u)=1$, and pairing with $\dd v$ gives $-c_1=\ell(\dd v)=0$. Exact terms have zero pairing with these smooth closed forms by periodic integration by parts. By Lemma~\ref{l:analyticprep}, $\theta=\dd v+\dd F$ with $F$ real, periodic and in $W^{1,2}$. Thus $U=(v+F)\bmod1$ is a globally defined circle-valued Sobolev potential for $\theta$.

In a regular rectangle, a real local lift $u$ satisfies $\dd u=\theta=c\,\dd x$, with $c\in L^2$, and hence $\partial_yu=0$ in distributions. Sobolev slicing makes almost every vertical restriction absolutely continuous with zero derivative, so $u(x,y)=f_I(x)$ almost everywhere. On a relatively compact subrectangle $I\times J$, one can take
\[
 f_I(x)=\frac1{|J|}\int_Ju(x,y)\,\dd y,
 \qquad
 f_I'(x)=\frac1{|J|}\int_J\partial_xu(x,y)\,\dd y
\]
in the weak sense. Fubini and Cauchy--Schwarz give $f_I\in W^{1,2}(I)$. It has a representative absolutely continuous on every compact subinterval, with $c(x,y)=f_I'(x)$ almost everywhere \cite[Theorem~2.36 and Remark~2.37(2), pp.~55--56]{Kin}. The resulting continuous circle-valued representatives agree almost everywhere on overlaps and therefore everywhere. They give a continuous version of $U$ on $\mathcal E\setminus Z$ that is constant on every regular vertical leaf. 

\textbf{4. The form $\theta$ vanishes almost everywhere off periodic leaves.}
The simplest model is an infinite-degree component $\mathcal S$ of $(h\circ\pi)^{-1}(C_i)$. It is the universal cover of a target cylinder, with a regular chart $(0,w)\times\RR$ and vertical leaves $\{x\}\times\RR$. By Step~3, $\theta=c(x)\,\dd x$ there. Hence for every $R>0$,
\[
       2R\int_0^w|c(x)|^2\,\dd x
         =\int_{(0,w)\times(-R,R)}|\theta|^2
         \le\|\theta\|_2^2=1/L.
\]
Letting $R\to\infty$ forces $c=0$ almost everywhere: positive energy cannot repeat indefinitely along a strip. This model does not cover preimages of minimal components or leaves accumulating at $Z$. For those, we first control the circle-valued potential near every point of $Z$, including the point over infinity. Two uses of the area formula then show, respectively, that almost every level consists of periodic leaves and that the exceptional levels carry zero energy.

Fix $p\in Z$ and a smooth conformal coordinate disk with $z(p)=0$. On that disk $U$ has a real $W^{1,2}$ lift $\widetilde U$ by its construction. For almost every $r$, its trace on the circle $C_r=\{|z|=r\}$ is absolutely continuous: apply Sobolev slicing in polar coordinates on annuli away from $0$, then exhaust the punctured disk \cite[Theorem~2.36 and Remark~2.37(2)]{Kin}. Let $V_p(r)$ be its total variation there and let $s$ denote Euclidean arc length in the $z$-coordinate. Cauchy--Schwarz on $C_r$ gives
\[
 V_p(r)^2
 =\left(\int_{C_r}\left|\frac{\dd\widetilde U}{\dd s}\right|\,\dd s\right)^2
 \le \underbrace{\left(\int_{C_r}1\,\dd s\right)}_{2\pi r}
     \int_{C_r}\left|\frac{\dd\widetilde U}{\dd s}\right|^2\,\dd s
 \le 2\pi r\int_{C_r}|\theta|^2\,\dd s .
\]
The last inequality uses that $\dd\widetilde U/\dd s$ is the tangential component of $\theta=\dd\widetilde U$. Dividing by $r$ and integrating over the radii, with $\dd A=\dd s\,\dd r$ in polar coordinates, gives
\[
      \int_0^\epsilon V_p(r)^2\,\frac{\dd r}{r}
           \le 2\pi\int_{|z|<\epsilon}|\theta|^2\,\dd A <\infty .
\]
Circles through $Z$ exclude at most countably many radii. For almost every remaining radius the circle-valued trace agrees everywhere with the continuous representative on that circle, by Fubini and continuity. Since $\int_0^\epsilon\dd r/r=\infty$, the energy bound gives such radii $r_j\downarrow0$ with $V_p(r_j)\to0$. After passing to a subsequence, the images of the circles shrink to a single value $\lambda_p\in\RR/\mathbb Z$. A leaf accumulating at $p$ crosses each sufficiently small chosen circle, since an arc of the leaf joins a point outside the circle to one inside. Its constant value under $U$ must therefore be $\lambda_p$. The same argument applies at the point over infinity even if other points of $Z$ accumulate there. Since $Z$ is countable, leaves accumulating at $Z$ account for at most countably many values.

\emph{First use of the area formula.} Cover $\mathcal E\setminus Z$ by countably many open flow rectangles whose closures lie in larger regular rectangles. By Step~3, on each rectangle $U(x,y)=F_I(x)$, where $F_I$ is absolutely continuous on a compact horizontal transversal $I$ spanning the rectangle. The one-dimensional area formula gives
\[
  \int_{\RR/\mathbb Z}\#\{x\in I:F_I(x)=\lambda\}\,\dd\lambda
       =\int_I |F_I'(x)|\,\dd x<\infty .
\]
For an absolutely continuous real lift, the formula follows by parametrizing its graph by $s(x)=x+\int_{\min I}^x|F_I'|$ and applying the Lipschitz area formula and monotone substitution. Summing over integer translates gives the circle-valued formula. Thus almost every level meets $I$ in finitely many points. Choose a single Borel null set $N\subset\RR/\mathbb Z$ containing the levels of infinite multiplicity on any of the countably many transversals, together with the values $\lambda_p$ for $p\in Z$.

For $\lambda\notin N$, the level $U^{-1}(\lambda)$ consists locally of finitely many vertical segments. It is a one-manifold without boundary, and its connected components are individual leaves. Each component is closed in the level, which is closed in $\mathcal E\setminus Z$ by continuity of $U$. Its closure also misses $Z$ by the circle argument. The component is therefore closed in the compact torus $\mathcal E$. As a compact connected one-manifold without boundary it is a circle, hence a periodic leaf.

\emph{Second use of the area formula.} Every nonperiodic regular leaf lies in $U^{-1}(N)$. Although $N$ has measure zero in the circle, its preimage can have positive area. On each transversal the area formula restricted to $N$ gives
\[
  \int_{F_I^{-1}(N)} |F_I'(x)|\,\dd x
    =\int_N\#\{x\in I:F_I(x)=\lambda\}\,\dd\lambda=0 .
\]
Thus $F_I'=0$ almost everywhere on $F_I^{-1}(N)$. In the rectangle $\theta=\dd U=F_I'(x)\,\dd x$, so integrating $|F_I'(x)|^2$ in $x$ and $y$ gives zero energy on the portion of $U^{-1}(N)$ in that rectangle. The rectangles form a countable cover of the regular surface, and $Z$ itself has area zero. Hence $\theta$ vanishes almost everywhere off the periodic leaves.

\textbf{5. Energy on the cylinders.}
Away from the exceptional set $h\circ\pi$ is a covering of the target with its marks and zeros removed and is locally flat in the pulled-back charts. A regular periodic leaf projects to a regular periodic target leaf. The latter belongs to a target cylinder $C_i$. A component over $C_i$ that contains a closed lift is a finite cover of the annulus, by the classification of connected annulus coverings. This describes all the cylinders needed to support $\theta$. The target has finitely many such maximal cylinders (Section~\ref{sec:flat}). Upstairs there are at most countably many.

On a torus cylinder $C=(0,w)\times(\RR/h_0\mathbb Z)$, write $x$ for the width coordinate and $y$ for the periodic coordinate. By Step~3, $\theta=v_C(x)\,\dd x$, with $v_C\in L^2(0,w)$. In this particular cylinder the potential also admits a single-valued real lift: the function $u_C(x,y)=\int_{x_0}^x v_C(s)\,\dd s$ is a real primitive of $\theta$, and $U-u_C$ is constant modulo one. Its period around the vertical core is zero. Orient the core $c_C$ in the direction of increasing $y$. For a smooth closed one-form $\eta$, every circle $x=\mathrm{const}$ has the same period $\int_{c_C}\eta$. Integrating first in $y$ and then in $x$ gives
\[
 \int_C\eta\wedge\theta=b_C\int_{c_C}\eta,\qquad
 b_C=-\int_0^w v_C(x)\,\dd x .
\]
Reversing the core orientation reverses the sign of $b_C$ and preserves the pairing.

\emph{Why the signed contributions sum to one.} Step~4 permits us to sum these pairings over the disjoint periodic cylinders. Even for countably many cylinders the sum converges absolutely, since $\sum_C|\int_C\eta\wedge\theta|\le\|\eta\|_2\|\theta\|_2$. By Step~3 the resulting sum is $\int_a\eta$ for every smooth closed $\eta$. Since periods distinguish real homology classes, this says
\[
                     a=\sum_C b_C[c_C]\quad\text{in }H_1(\mathcal E;\RR).
\]
A null-homologous core contributes zero to this identity, even if its cylinder carries energy. Since $a\ne0$, some essential core contributes. All essential cores are disjoint simple curves on the torus, so they represent one primitive homology class up to sign. The identity forces their common unoriented slope to be $[a]$. Orient each essential core to represent $a$ and take $b_C$ with the corresponding sign. Evaluating the identity on a closed form of period one on $a$ now gives $\sum_{C\ \mathrm{essential}} b_C=1$, allowing either sign for the individual terms.

\emph{The energy bound.} Put $M_C=\Mod(C)=w/h_0$ at time zero. Cauchy--Schwarz across the width gives
\[
 b_C^2\le w\int_0^w v_C(x)^2\,\dd x
       =M_C\,\|\theta\|_{2,C}^2,
 \qquad \|\theta\|_{2,C}^2=h_0\int_0^w v_C(x)^2\,\dd x .
\]
In the following sums $C$ ranges over the essential cylinders. Their combined energy is at most $\|\theta\|_2^2=1/L$ by Step~3. Applying Cauchy--Schwarz to the signed sum yields
\[
 1=\Bigl(\sum_C b_C\Bigr)^2
   \le \Bigl(\sum_C M_C\Bigr)\Bigl(\sum_C b_C^2/M_C\Bigr)
   \le \frac{1}{L}\sum_C M_C .
\]
Thus $L\le\sum_C M_C$ when $L>0$. To evaluate this sum, fix a target cylinder $C_i$. Distinct closed plane lifts of its core bound disjoint disks. An essential torus core requires the bounded disk to contain two of the three finite marks of $B$, so at most one such plane lift $\alpha_i$ exists. Its torus cores have slope $[a]$ by the homology argument. If $\alpha_i$ has degree $d_i$ over the target core, the lifted plane cylinder has modulus $m_i/d_i$. By Step~1 it has two torus copies of that same modulus. Hence $\sum_C M_C=2\sum_{i:\,\alpha_i\ \mathrm{exists}}m_i/d_i$, which gives the upper bound in \eqref{eq:quartet}. For $L=0$ the bound is immediate.
\end{proof}

\subsection{Excluding a cylinder of \texorpdfstring{$G$}{G}}\label{sec:excluding}
We now combine Theorem~\ref{l:record}, the four-point bound and the Matrix Lemma. The proof first constructs an invariant multicurve, then chooses four marks isolating one core, and finally obtains the incompatible bounds \eqref{eq:matrix}.

\begin{prop}\label{l:cylinder-levy}
Assume $s=|A|\ge3$. Let $G$ be any measured foliation with $\Ext_oG=1$ satisfying \eqref{eq:displacement}, and let $q$ be its Hubbard--Masur differential at $o$. If $q$ has a positive-width vertical cylinder, then $f$ has a Levy cycle relative to $A^+$.
\end{prop}

\begin{proof}
\textbf{1. Construct an invariant multicurve.}
Let $\Gamma$ be the nonempty set of core classes of the positive-width vertical cylinders of $q$. Fix $\gamma\in\Gamma$ and let $w_\gamma>0$ be the transverse width of its cylinder at time zero. The weighted curve $w_\gamma\gamma$ is a component of $G$, so $w_\gamma\,i(\gamma,U)\le i(G,U)$ for every test foliation $U$. On any marked structure $X$, divide by $\sqrt{\Ext_X U}$ and take the supremum over $U\ne0$. Squaring the duality formula \eqref{eq:duality} gives $w_\gamma^2\Ext_X\gamma\le\Ext_XG$.
Since $q$ has area one, $\Ext_{g(t)}G=e^{-2t}$. The assumed bound \eqref{eq:displacement}, supplied in the application by Theorem~\ref{l:record}, gives $\Ext_{\sigma^k g(t)}G\le C e^{-2t}$ with one $C$ for all $t$ and $k$, as in \eqref{eq:ray}. Taking $X=\sigma^k g(t)$ and inverting the bound on $\Ext_X\gamma$ gives
\begin{equation}\label{eq:cgamma}
             1/\Ext_{\sigma^k g(t)}\gamma\ \ge\
             c_\gamma e^{2t},\qquad c_\gamma=w_\gamma^2/C>0
                    \quad(\gamma\in\Gamma).
\end{equation}
Fix a depth $k$ and an essential finite lift $\alpha$ of a core $\beta\in\Gamma$ of initial modulus $m_\beta$. The lifted annulus on $\sigma^k g(t)$ has modulus $e^{2t}m_\beta/d$ for its fixed degree $d$. Hence the extremal lengths of both $\alpha$ and each $\gamma\in\Gamma$ decay to zero as $t\to\infty$ with $k$ fixed, on that same marked surface. Minsky gives $i(\alpha,\gamma)=0$. To compare lifts at depths $i\le j$, follow both forward by $f^i$. The first reaches its original cylinder core and the second reaches an essential lift of depth $j-i$ by the preceding covering argument. We have just shown that this core and this lift can be drawn disjointly, using parallel copies if their classes coincide. Pulling the two representatives back together by $f^i$ gives disjoint inverse images, and isotopy lifting ensures that their components include representatives of the two prescribed classes.

Thus all the distinct essential lifted classes, together with $\Gamma$, are pairwise compatible. Pairwise compatibility implies simultaneous disjointness: for a complete hyperbolic metric on $\widehat\CC\setminus A^+$, geodesic representatives of distinct classes with zero intersection are disjoint \cite[Prop.~1.3 and Cor.~1.9]{FMar}. Taking the union over all depths therefore gives a single multicurve $\Lambda$ with at most $r=s-2$ distinct classes, since any $s-1$ of them would already exceed the bound for disjoint essential curves on the $(s+1)$-punctured sphere (Section~\ref{sec:flat}). This multicurve is invariant because every essential one-step lift of a class occurring at depth $k$ occurs at depth $k+1$, with the original cores at depth zero. Its matrix $M_\Lambda$ is the Thurston matrix.

\textbf{2. Isolate one core using four marks.}
Fix $\gamma\in\Gamma$. Choose $B_\gamma$ so that $\gamma$ is the only class of $\Lambda$ remaining essential after the other marks are forgotten. The same four labels will be used for every depth $k$ and time $t$.

Choose disjoint representatives of $\Lambda$ and form their inclusion tree of bounded disk label sets, adding all singleton leaves and the root $A$. Disjoint representatives with the same label set cobound an unmarked annulus, so different vertices have different sets. Choose labels $p_1,p_2$ from two different children of $\gamma$ and a label $p_3$ from a sibling of $\gamma$ in its parent. Exactly the vertex $\gamma$ contains $p_1,p_2$ while excluding $p_3$. Indeed a proper descendant of $\gamma$ misses $p_1$ or $p_2$, a proper ancestor contains the label set of the parent of $\gamma$, by laminarity, and hence $p_3$, and a disjoint vertex misses both $p_1$ and $p_2$. Each nonsingleton vertex, including the root, has at least two children because singleton labels have been included; in particular the sibling used to choose $p_3$ exists. Every vertex other than $\gamma$ contains zero, one or three of $p_1,p_2,p_3$, so the other possible two-point splits are excluded as well. Thus, for $B_\gamma=\{p_1,p_2,p_3,\infty\}$, $\gamma$ is the unique class of $\Lambda$ that remains essential after all other marks are forgotten. Figure~\ref{fig:tree} illustrates this choice.

\begin{figure}[ht]\centering
\begin{tikzpicture}[scale=0.88]
\begin{scope}
\draw[thick,pine] (0,0) ellipse (3.0 and 1.7);
\node[font=\scriptsize,pine] at (2.45,1.5) {$\gamma'$};
\draw[very thick,rust] (-1.0,0) ellipse (1.5 and 1.05);
\node[font=\scriptsize,rust] at (-1.0,1.28) {$\gamma$};
\draw[thick,pine] (-1.55,0) ellipse (0.62 and 0.6);
\node[font=\scriptsize,pine] at (-1.55,0.8) {$\delta$};
\fill (-1.8,0) circle (1.6pt) node[font=\scriptsize,below,inner sep=1.5pt]{$p_1$};
\fill (-1.3,0) circle (1.6pt);
\fill (-0.3,0.1) circle (1.6pt) node[font=\scriptsize,above]{$p_2$};
\fill (1.6,0.3) circle (1.6pt) node[font=\scriptsize,above]{$p_3$};
\fill (1.9,-0.6) circle (1.6pt);
\fill (3.9,0.8) circle (1.6pt);
\node[font=\scriptsize] at (3.6,-1.5) {$\infty$};
\end{scope}
\begin{scope}[shift={(8.2,1.5)},level distance=9mm,
  every node/.style={font=\scriptsize,draw=pine,rounded corners=2pt,fill=mist,inner sep=2pt},
  level 1/.style={sibling distance=24mm},level 2/.style={sibling distance=13mm},level 3/.style={sibling distance=9mm}]
\node {root $A$}
  child {node {$\gamma'$}
     child {node[draw=rust,very thick,fill=sand] {$\gamma$}
        child {node {$\delta$} child{node{$p_1$}} child{node{$\cdot$}}}
        child {node {$p_2$}}}
     child {node {$p_3$}}
     child {node {$\cdot$}}}
  child {node {$\cdot$}};
\end{scope}
\end{tikzpicture}
\caption{Isolating $\gamma$ in a multicurve (left) using its inclusion tree (right). Take $p_1$ inside $\delta$, $p_2$ between $\delta$ and $\gamma$, and $p_3$ between $\gamma$ and $\gamma'$, together with $\infty$. After forgetting the other marks, $\delta$ encloses only $p_1$ and $\gamma'$ leaves only $\infty$ outside. Only $\gamma$ separates two retained marks from two.}
\label{fig:tree}
\end{figure}
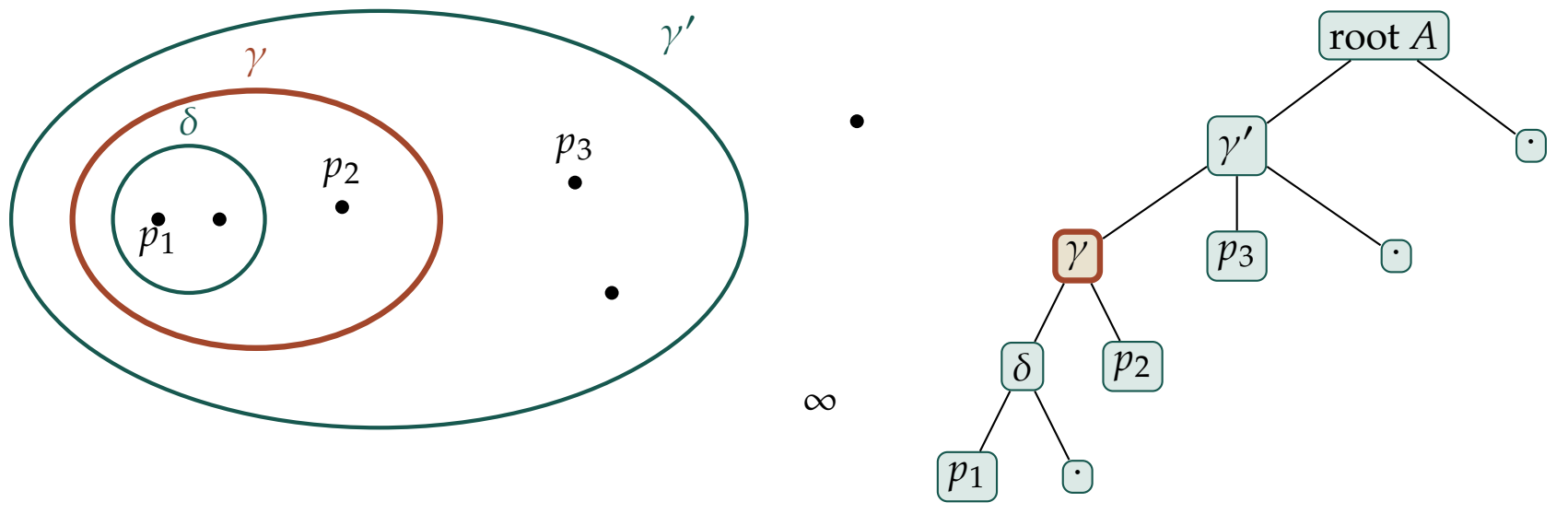

\textbf{3. Apply the four-point bound and the Matrix Lemma.}
Fix $k\ge1$ and write $\gamma_B$ for the class of $\gamma$ on the time-$t$ source sphere after forgetting the marks outside $B$. Each embedded annulus of full-marking core $\gamma$ remains embedded after forgetting the other marks and has core $\gamma_B$. The simple-curve annulus formula and \eqref{eq:cgamma} give
\[
 \frac{1}{\Ext_{B,t}(\gamma_B)}
 =\sup_{\operatorname{core}(D)=\gamma_B}\Mod D
 \ \ge\ \sup_{\operatorname{core}(D)=\gamma\ {\rm rel.}\ A^+}\Mod D
 =\frac{1}{\Ext_{\sigma^k g(t)}\gamma}\ \ge\ c_\gamma e^{2t}.
\]
On the four-marked sphere $1/\Ext_{B,t}(\gamma_B)=H_a(\mathcal E_t)/2$, so $e^{-2t}H_a(\mathcal E_t)/2\ge c_\gamma$ for every $t$. Write $L_k=\lim_{t\to\infty}e^{-2t}H_a(\mathcal E_t)$ for the limit at this fixed depth. Passing to the limit gives $L_k/2\ge c_\gamma$. The same $c_\gamma=w_\gamma^2/C$ works at every depth because $\gamma$ and its initial width were fixed before choosing $k$, and \eqref{eq:displacement} bounds all iterates with the same $C$; Theorem~\ref{l:record} supplies this bound in the application.

Now \eqref{eq:quartet} bounds $L_k/2$ from above by the sum of $m_i/d_i$ over the finite lifts of target cylinder cores in the four-marked class $\gamma_B$, where $m_i$ is the initial target modulus and $d_i$ the degree of the lift. Thus these lifts together must contribute at least $c_\gamma$ at every depth. Extending the vector $m$ of initial cylinder moduli by zero on $\Lambda\setminus\Gamma$, we obtain
\begin{equation}\label{eq:matrix}
                      0<c_\gamma\le \frac{L_k}{2}\le (M_\Lambda^k m)_\gamma
                          \qquad(k\ge1).
\end{equation}
Every counted lift is full-marking essential and belongs to $\Lambda$; the choice of $B$ selects precisely class $\gamma$. Its ancestors are essential at every intermediate depth. Conversely, a chain of actual lift components gives a lift of the iterate, with degree the product of its degrees and with a unique chain of ancestors. Isotopy lifting identifies the next components for any representative of an intermediate class. Thus multiplicities are preserved: two successive degrees $d_1,d_2$ contribute $1/(d_1d_2)$, and summing over components and intermediate classes gives $\sum_{\beta\in\Lambda}M_{\gamma\beta}M_{\beta\delta}$. Induction gives $M_\Lambda^k$. Having taken $t\to\infty$ at fixed $k$, we now let $k\to\infty$ in \eqref{eq:matrix}.

If $f$ has no Levy cycle, the Matrix Lemma gives $M_\Lambda^k\to0$, contradicting the uniform positive lower bound in \eqref{eq:matrix}. Thus $f$ has a Levy cycle relative to $A^+$.
\end{proof}

\section{A finite measure on the initial flat surface}\label{sec:transport}
At each fixed depth $k$, we construct a measure on the initial target that controls the mass on regions chosen at later times. The main input is convergence of near-minimizers on one Hilbert space, obtained by transporting the source differentials through their markings.

\subsection*{Pulling back and pushing forward a quadratic differential}\label{bg:transport}
\addcontentsline{toc}{subsection}{Pulling back and pushing forward a quadratic differential}
Markings identify the surfaces at different times. Holomorphic pullback instead transfers a differential through $H$, with coefficient $(p\circ H)(H')^2$, and pushforward sums the inverse-branch contributions as in Lemma~\ref{l:pushforward}. Their complex arguments can cancel; positive multiples of a common differential add with equality. The construction below makes this loss of mass tend to zero.

The pulled-back geometry $Q=H^*q$ has infinite area for transcendental $H$: Picard and the area formula give infinitely many sheets above almost every target point. The Hubbard--Masur differentials $r_t$ have finite area. Their comparison with $Q_t$ will produce an integrable density $v^2$ on the source and, by summing over fibers, a finite measure on the target.

Fix $k$ and an entire representative $H$ for the pullback by $f^k$ from the target $o$ to the source $\sigma^k o$. Set $Q=H^*q$. Take out the preimages of marked target points and of the zeros of $q$, and the source marks. Their union is locally finite in the plane. Write $Y$ for the regular complement. Use the affine ray marking $\phi_t$ on the target and solve its pulled-back Beltrami equation on the source to obtain $\psi_t$ in the induced marking. Then $H_t=\phi_t H\psi_t^{-1}$ is entire and represents the same marked topological cover at time $t$. In regular charts of $Q_t=H_t^*q_t$ the lift $\psi_t$ is the stretch $(x,y)\mapsto(e^tx,e^{-t}y)$ up to the regular chart transitions. This follows by using $H,H_t$ themselves as the local coordinates into the corresponding target charts. In particular $\psi_t^*|Q_t|=|Q|$ as area measures. For transcendental $H$, this source area is infinite. Throughout this section, curves and foliations are identified by the source markings, including the initial marking of $\sigma^k o$.

\emph{Order of choices.} The already fixed $o,G,C$ determine $q$ and the ray. The depth $k$, arbitrary here and prescribed in Section~\ref{sec:mass}, determines all source objects and the measure $\mu$ below. Only afterwards are sequences $t_n\to\infty$ chosen. This ensures that the same measure is available at every late time.

Let $r_t$ be the Hubbard--Masur differential for $G$ on $\sigma^k g(t)$. Its area norm is $\Ext_{\sigma^k g(t)}G$ by the punctured-surface Hubbard--Masur area identity \cite[Theorem~5.1]{Ant}. Since $d_T(\sigma^k g(t),\sigma^k g(t+h))\le h$ for $h\ge0$, the distortion bound~\eqref{eq:kerckhoff} gives $\|r_{t+h}\|_1\ge e^{-2h}\|r_t\|_1$. Thus
\begin{equation}\label{eq:energy}
 E_t:=e^{2t}\|r_t\|_1\uparrow E,\qquad
              0<E_0\le E\le C<\infty,
\end{equation}
where the upper bound is \eqref{eq:ray} and positivity follows from $G\ne0$. The limit and all the source objects in this section concern this one fixed depth.

Convergence of areas alone does not locate their mass. The following Hilbert-space lemma uses the parallelogram identity: if the midpoint remains admissible, two near-minimizers must be close. In the application, the competitors are absolute coefficients of square roots of $r_t$; the proof of Proposition~\ref{l:transport} verifies admissibility by a transverse-length estimate. The resulting total-variation convergence controls every measurable region, including regions chosen after seeing the late surface.

\begin{lemma}\label{l:hilbert}
Let $\mathcal H$ be a real Hilbert space, put $J_t(a,b)=\|a\|^2+e^{4t}\|b\|^2$ on $\mathcal H\times\mathcal H$, and let $(a_s,b_s)\in\mathcal H\times\mathcal H$ for $s\ge0$. Suppose that $E_t=J_t(a_t,b_t)$ increases to a finite limit $E$, and that $E_t\le J_t(a,b)$ for every $t\ge0$ and every finite convex combination $(a,b)$ of the pairs $(a_s,b_s)$ (each taken at an arbitrary $s\ge0$). Then $e^{2s}\|b_s\|\to0$, and $a_s$ converges to some $v\in\mathcal H$ with $\|v\|^2=E$.
\end{lemma}
\begin{proof}
For $s\ge1$, the hypothesis with $t=s-1$ and $(a,b)=(a_s,b_s)$ gives $E_{s-1}\le J_{s-1}(a_s,b_s)=E_s-(1-e^{-4})e^{4s}\|b_s\|^2$, so
\[
 (1-e^{-4})e^{4s}\|b_s\|^2\le E_s-E_{s-1}\longrightarrow0 .
\]
For $s,u\ge t$ apply the hypothesis to the midpoint of $(a_s,b_s)$ and $(a_u,b_u)$. The parallelogram identity for $J_t$, together with $J_t\le J_s$ and $J_t\le J_u$, gives
\begin{align*}
 \|a_s-a_u\|^2\le J_t(a_s-a_u,b_s-b_u)
 &=2J_t(a_s,b_s)+2J_t(a_u,b_u)-4J_t\Big(\frac{a_s+a_u}2,\frac{b_s+b_u}2\Big)\\
 &\le 2E_s+2E_u-4E_t\le 4(E-E_t).
\end{align*}
Thus $(a_s)$ is Cauchy. Its limit $v$ satisfies $\|v\|^2=\lim\big(E_s-e^{4s}\|b_s\|^2\big)=E$.
\end{proof}

\begin{prop}\label{l:transport}
There is $v\ge0$ in $L^2(Y,|Q|)$ with
\[
 \|e^{2t}r_t-(v^2\circ\psi_t^{-1})Q_t\|_1\longrightarrow0,\qquad
 \int_Yv^2|Q|=E .
\]
The fiber sum $S(w)=\sum_{H(z)=w}v(z)^2$, defined as a measurable nonnegative sum off the target marks, belongs to $L^1(|q|)$ and is finite for $|q|$-almost every $w$, even when the fibers are infinite. With $P_t=e^{2t}(H_t)_*r_t$, we have
\begin{equation}\label{eq:pushforward}
                 \|P_t-(S\circ\phi_t^{-1})q_t\|_1\to0,\qquad
                 \mu=S|q|,\qquad \mu(\widehat{\CC})=E.
\end{equation}
Here $\psi_t$ and $\phi_t$ are the Teichmüller markings from time zero on the source and target: they are affine in the natural charts of $Q$ and $q$, respectively, rather than in global plane coordinates. Set $v_t=v\circ\psi_t^{-1}$ and $S_t=S\circ\phi_t^{-1}$. The comparison differentials $v_t^2Q_t$ and $S_tq_t$ are measurable. Since $\phi_t^*(S_t|q_t|)=\mu$, for every measurable $U$ on the initial target the difference between $\int_{\phi_t(U)}|P_t|$ and $\mu(U)$ is bounded by $\|P_t-S_tq_t\|_1$, uniformly over $U$.
\end{prop}
\begin{proof}
All $L^2$ norms below use the single area element $|Q|$ on $Y$. In an initial regular chart $x+iy$, take the absolute values $a_s,b_s$ of the coefficients of $\Re\sqrt{r_s}$ in the basis $\dd x,\dd y$, using $\psi_s$. Changes of regular chart have linear part $\pm I$, so both functions are single-valued. In the conformal chart $\zeta_s=e^sx+ie^{-s}y$ the coefficient of $\sqrt{r_s}$ has real and imaginary parts of absolute values $e^{-s}a_s$ and $e^sb_s$. Consequently
\begin{equation}\label{eq:Es}
 E_s=e^{2s}\|r_s\|_1=\|a_s\|_2^2+e^{4s}\|b_s\|_2^2 .
\end{equation}
Equation \eqref{eq:energy} gives $E_s\uparrow E$.

For a finite convex combination $(a,b)=\sum_j\lambda_j(a_{s_j},b_{s_j})$ and \emph{every} time $t$, we claim
\begin{equation}\label{eq:Jt}
               E_t\le J_t(a,b):=\|a\|_2^2+e^{4t}\|b\|_2^2.
\end{equation}
If $A_s,B_s$ are the signed coefficients in the local expression $\psi_s^*(\Re\sqrt{r_s})=A_s\,\dd x+B_s\,\dd y$, then $a_s=|A_s|$ and $b_s=|B_s|$, so the triangle inequality bounds the transverse length element $|A_s\,\dd x+B_s\,\dd y|$ by $a_s|\dd x|+b_s|\dd y|$. Thus on a piecewise smooth regular representative $\delta$ of a simple class $\beta$, each $\int_\delta(a_{s_j}|\dd x|+b_{s_j}|\dd y|)$ bounds from above the transverse length of $\delta$ for the vertical foliation of $r_{s_j}$. That length is at least $i(G,\beta)$, since each $r_{s_j}$ represents the same marked foliation. The convex combination has the same lower bound. At time $t$ the conformal line element
\[
        e^{-t}\sqrt{a^2+e^{4t}b^2}\,|\dd\zeta_t|
\]
dominates $a|\dd x|+b|\dd y|$ by Cauchy--Schwarz and has area $e^{-2t}J_t(a,b)$. The functions $a_s,b_s$ are continuous in regular charts, and only finite convex combinations are used here. Take an embedded annulus $\mathcal A$ of core $\beta$ on the time-$t$ source, with cylinder coordinates $(u_1,u_2)\in(\mathbb R/\mathbb Z)\times(0,M)$. All but countably many circular cores avoid the discrete exceptional set and are smooth in the natural charts. If $\rho$ is the test density above, the transverse-length bound and Cauchy--Schwarz give, at almost every height,
\[
 i(G,\beta)^2\le\left(\int_0^1\rho(u_1,u_2)\,\dd u_1\right)^2
               \le\int_0^1\rho(u_1,u_2)^2\,\dd u_1.
\]
Integrating over $0<u_2<M$ gives $M i(G,\beta)^2\le\operatorname{Area}_\rho(\mathcal A)\le e^{-2t}J_t(a,b)$. This uses the finite area of the test metric, even though $|Q|$ has infinite mass.
Keeping $\beta$ fixed, we take the supremum of $M$ over all embedded annuli with core class $\beta$. Since this supremum is $1/\Ext_{\sigma^k g(t)}\beta$, we obtain
\[
   \frac{i(G,\beta)^2}{\Ext_{\sigma^k g(t)}\beta}\le e^{-2t}J_t(a,b).
\]
We now take the supremum over $\beta$. The quotient is unchanged when $\beta$ is given a positive weight, so by density of weighted simple curves and continuity its supremum is the one in the squared duality formula~\eqref{eq:duality}, namely $\Ext_{\sigma^k g(t)}G$. Recall that $r_t$ is the Hubbard--Masur differential on $\sigma^k g(t)$ with vertical foliation $G$, so this extremal length equals its area $\|r_t\|_1=e^{-2t}E_t$ by \eqref{eq:energy}. Multiplication by $e^{2t}$ now gives $E_t\le J_t(a,b)$, completing the proof of the claimed inequality~\eqref{eq:Jt}.

By \eqref{eq:Es}, \eqref{eq:energy} and \eqref{eq:Jt}, the pairs $(a_s,b_s)$ in $L^2(Y,|Q|)$ satisfy the hypotheses of Lemma~\ref{l:hilbert}. Hence $e^{2s}\|b_s\|_2\to0$, and $a_s$ converges in $L^2$ to a limit $v$ with $\|v\|_2^2=E$; moreover $v\ge0$, since each $a_s\ge0$.

On the initial source write $F_t=(e^{2t}r_t/Q_t)\circ\psi_t$ in regular charts. Then $\Re F_t=a_t^2-e^{4t}b_t^2$ and $|\Im F_t|=2a_t e^{2t}b_t$. In particular, with $B_t=e^{2t}b_t$ the global error is bounded by
\[
 \|F_t-v^2\|_1\le (\|a_t\|_2+\|v\|_2)\|a_t-v\|_2
           +\|B_t\|_2^2+2\|a_t\|_2\|B_t\|_2\longrightarrow0 .
\]
This estimate holds on the entire cover, including its infinite-area end when present. Indeed $\|e^{2t}r_t-v_t^2Q_t\|_1$ on the time-$t$ source equals $\|F_t-v^2\|_{L^1(Y,|Q|)}$, since $\psi_t^*|Q_t|=|Q|$.

In these markings the pushforward can be written on a simply connected regular target disk $D$ using its inverse branches $\xi_j$. For an integrable measurable differential $u$ its coefficient on $D$ is $\sum_j(u\circ\xi_j)(\xi_j')^2$, writing $u$ also for its coefficient in the source chart. The measurable case of Lemma~\ref{l:pushforward} gives $\|(H_t)_*u\|_1\le\|u\|_1$, also when there are infinitely many sheets. For $u=v_t^2Q_t$ every summand is a nonnegative multiple of the same target differential $q_t$, since $\xi_j^*Q_t=q_t$. Equality of the masses follows in this case. The removed sets and their inverse images have zero area.

Tonelli gives $\int S|q|=\int v^2|Q|=E<\infty$, which proves the asserted almost-everywhere finiteness of $S$. Under these Teichmüller markings, $(H_t)_*(v_t^2Q_t)=S_tq_t$. This gives \eqref{eq:pushforward}, including convergence of the masses in total variation. In particular the ratio required in the inverse-branch argument is
\[
 \frac{\|(H_t)_*r_t\|_1}{\|r_t\|_1}=\frac{\|P_t\|_1}{E_t}\longrightarrow1,
\]
since the numerator on the right tends to $\int S|q|=E$ and $E_t\to E>0$. The same depth $k$ occurs in both quantities.
\end{proof}

\begin{remark}
Proposition~\ref{l:transport} has two consequences: near-equality $\|(H_t)_*r_t\|_1/\|r_t\|_1\to1$, used on an inverse branch in Proposition~\ref{p:degree}, and one limiting measure $\mu$ at fixed depth. Total variation lets us select a cylinder on a late surface and compare its mass at different times using this same measure (Sections~\ref{sec:mass}--\ref{sec:persistence}). Section~\ref{sec:toHSS} compares this conclusion with the near-equality method of HSS.
\end{remark}

Write $\mu_t=(\phi_t)_*\mu$, and recall $|q_t|=(\phi_t)_*|q|$. The bound in \eqref{eq:pushforward} is uniform for all measurable sets, including sets that vary with $t$. A second uniform bound is
\begin{equation}\label{eq:omega}
 \mu_t(U)\le\omega\bigl(|q_t|(U)\bigr),\qquad
 \omega(\delta)=\sup_{|q|(V)\le\delta}\mu(V)\ \longrightarrow0
                 \quad(\delta\downarrow0).
\end{equation}
Indeed $\mu=S|q|$ with $S\in L^1(|q|)$, so $\omega(\delta)\le M\delta+\int_{\{S>M\}}S|q|$ for each $M>0$. Both bounds remain valid after any affine change of target coordinate when all measures are transported together.

\emph{Flat length and intersection at different depths.} A straight segment in a natural coordinate has Euclidean length $\sqrt{(\Delta x)^2+(\Delta y)^2}$. Its two coordinate variations measure transverse lengths for the vertical and horizontal foliations. At a zero or puncture a limiting shortest curve may be a broken geodesic made of \emph{saddle connections}, straight segments joining singularities with regular interiors. A regular core of a flat cylinder has a smooth closed geodesic representative instead. For a puncture at finite flat distance the marked homotopy data specify how a representative approaches and passes it. The comparison below applies in all of these cases. It converts small extremal lengths on several different pullback structures into small flat lengths on one surface. This common surface will let us compare intersections.

Let $\ell_{q_t}(\alpha)$ be the infimum of the flat lengths of representatives of the essential simple class $\alpha$ relative to $A^+$. There is a constant $D_{\rm flat}$ depending only on the topological type such that, for each $j\ge0$,
\begin{equation}\label{eq:flat}
 \ell_{q_t}(\alpha)\le C_j\sqrt{\Ext_{\sigma^j g(t)}\alpha},\qquad C_j=D_{\rm flat}\bigl(\sqrt C+e^{d_T(o,\sigma^jo)}\bigr),
\end{equation}
uniformly in $t$ and $\alpha$; in fact one may take $D_{\rm flat}=1$. Write $\mathcal H_q$ for the horizontal measured foliation at time zero. Choose the flat geodesic representative of $\alpha$ using the CAT(0) completion of the universal cover, with every mark retained as a puncture, including regular points of $q$ \cite[\S\S2.1, 2.4, Remark~8]{DLR}. Its lift and a regular horizontal or vertical leaf cannot form a bigon, since geodesics in this completion are unique. Thus the same representative minimizes transverse variation in both directions and realizes both intersection numbers, as in the proof of \cite[Lemma~9]{DLR}. Small detours at punctures preserve the marked class and add arbitrarily little flat and transverse length. Using $|\dd\zeta_t|\le|\dd\Re\zeta_t|+|\dd\Im\zeta_t|$, we obtain the two-line estimate
\begin{align*}
 \ell_{q_t}(\alpha)
 &\le e^ti(G,\alpha)+e^{-t}i(\mathcal H_q,\alpha)\\
 &\le\bigl(\sqrt C+e^{d_T(o,\sigma^jo)}\bigr)
                   \sqrt{\Ext_{\sigma^jg(t)}\alpha}.
\end{align*}
For the second line, Minsky's inequality and \eqref{eq:ray} bound the first term, while $\Ext_o\mathcal H_q=1$ and $d_T(o,\sigma^jg(t))\le d_T(o,\sigma^jo)+t$ bound the second by extremal-length distortion. This proves \eqref{eq:flat} directly.
In Section~\ref{sec:mass}, $n$ indexes the times $t_n\to\infty$, and the curve may vary with $n$. Since the pullback depth $k$ is fixed, the bound \eqref{eq:flat} with constant $\max_{0\le i<k}C_i$ holds for every $n$ and every intermediate depth $0\le j<k$.

\section{Near equality on an inverse branch}\label{sec:sheet}
We analyze a nonzero limit of the pushforward in a cluster chart. Near-equality selects an inverse branch with small angular error relative to its own mass, which may tend to zero. After normalizing that branch, we continue it across a large annulus. Unbounded covering degree gives an impossible univalent limit on a half-plane. Bounded degree gives a proper disk that can be capped to a polynomial, whose ramification bounds the number of marked poles.

\subsection*{An inverse branch on an annulus}\label{bg:branches}
\begin{figure}[ht]\centering
\begin{tikzpicture}[scale=0.82]
\begin{scope}
\draw[thick,pine,fill=mist] plot[smooth cycle,tension=0.8] coordinates {(-1.8,0) (-1.2,1.3) (0.2,1.6) (1.7,0.9) (1.9,-0.6) (0.6,-1.5) (-1.1,-1.3)};
\foreach \p in {(-0.9,0.4),(-0.2,-0.5),(0.6,0.3),(1.0,-0.6),(0.0,0.9)}{\fill[rust] \p circle (1.7pt);}
\fill[sand] (-0.7,-0.75) circle (0.22); \node[font=\tiny] at (-0.7,-0.75) {$\psi_n$};
\node[font=\scriptsize,pine] at (1.75,1.45) {$D_n$};
\node[font=\scriptsize,align=center] at (0,-2.1) {source: $\ge d+2$ marked poles of $\eta_n$};
\draw[thick,densely dotted,ink] (-2.6,-1.8)--(-2.6,1.8);
\node[font=\scriptsize,ink,rotate=90] at (-2.95,0) {other sheets};
\end{scope}
\draw[-{Stealth},thick] (2.35,0)--(3.55,0) node[midway,above,font=\scriptsize]{$H_n$} node[midway,below,font=\scriptsize]{degree $d$};
\begin{scope}[shift={(6.0,0)}]
\draw[pine,thick] (0,0) circle (1.8);
\draw[pine!50,dashed] (0,0) circle (0.5);
\foreach \p in {(-0.2,0.1),(0.15,-0.2),(0.25,0.2)}{\fill[rust] \p circle (1.4pt);}
\node[font=\scriptsize,pine] at (0,2.08) {$|z|=\rho_n=\sqrt{R_n}$};
\node[font=\scriptsize] at (0,-0.78) {$|z|\le1$};
\fill[sand] (1.1,0.6) circle (0.28); \node[font=\scriptsize] at (1.1,0.6) {$e^{V}$};
\node[font=\scriptsize,align=center] at (0,-2.1) {target chart of the cluster};
\end{scope}
\draw[-{Stealth},thick,plum] (8.25,0.2)--(9.45,0.2) node[midway,above,font=\scriptsize,plum]{cap};
\begin{scope}[shift={(11.6,0)}]
\fill[plum!8] (0,0) circle (1.9);
\draw[pine,thick,fill=mist] plot[smooth cycle,tension=0.8] coordinates {(-1.1,0) (-0.7,0.8) (0.1,1.0) (1.0,0.55) (1.15,-0.35) (0.35,-0.9) (-0.65,-0.8)};
\node[font=\scriptsize,align=center,plum] at (0,-2.25) {glue $\{|\xi|\ge\rho_n^{1/d}\}\cup\{\infty\}$ along $\partial D_n$:\\ $H_n|_{D_n}$ becomes a polynomial $P_n$};
\end{scope}
\end{tikzpicture}
\caption{The selected inverse branch determines an annular component. At finite degree $d$, a lifted circle bounds the proper source disk $D_n$. Attaching an exterior disk gives a polynomial model, whose finite ramification $d-1$ yields at least $d+2$ marked poles in $D_n$ (Steps~4--6). Other sheets of the entire map remain outside this disk.}
\label{fig:sheet}
\end{figure}
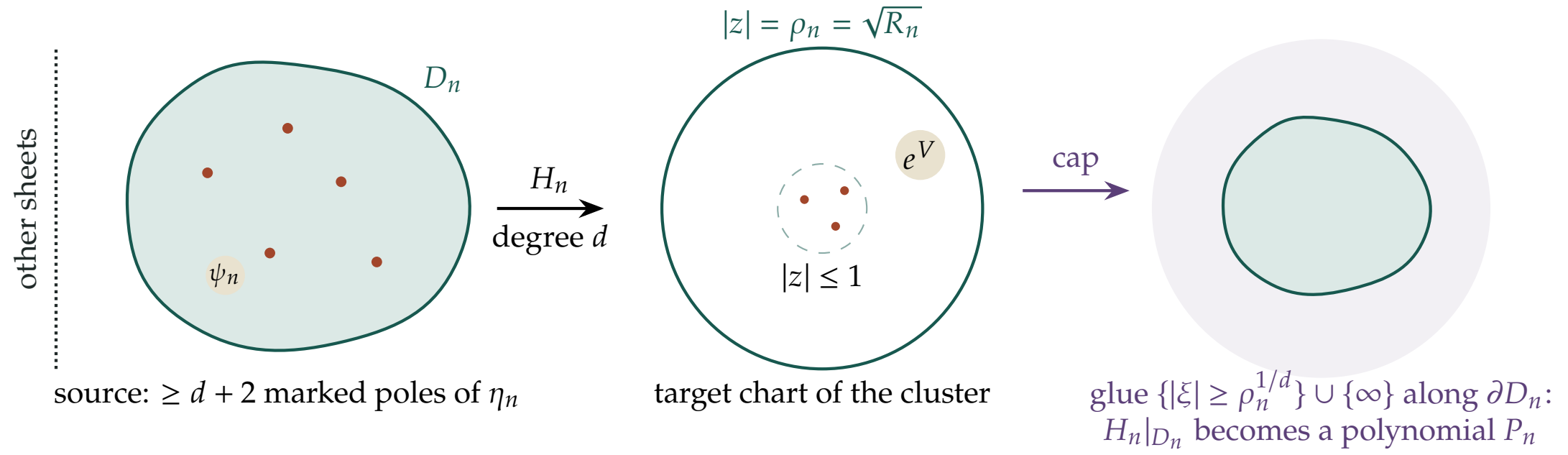

\addcontentsline{toc}{subsection}{An inverse branch on an annulus}
Continuation around an unmarked annulus permutes the inverse branches. On a connected lift, the orbit is either a finite cycle of length $d$ or a copy of $\mathbb Z$, corresponding to a degree-$d$ cover or the universal cover \cite[\S1.3, Theorem~1.38]{Hat}. In the coordinate $w=\log z$, a finite return gives period $2\pi i d$. For annuli $1<|z|<R_n$ with covering degrees $d_n$, the domains of injectivity exhaust the right half-plane when both $R_n\to\infty$ and $d_n\to\infty$. In the finite-degree case an exterior disk can be attached in the coordinate where the map is $\xi\mapsto\xi^d$; Figure~\ref{fig:sheet} shows the resulting polynomial model. The original entire map can have infinitely many other sheets outside the selected disk.

The pole count comes from the polynomial model: a nonzero integrable sphere differential has at least three finite simple poles, and degree $d$ allows only $d-1$ units of finite ramification. Step~6 proves that at least $d+2$ unramified poles survive and transfers them back to the original source disk.

\begin{samepage}
\begin{prop}\label{p:degree}
Let $H_n$ be transcendental entire maps. Each source has a set of $s$ finite marks. All finite singular values of $H_n$ and all images of source marks belong to finite target marking sets of uniformly bounded cardinality. In affine target coordinates all these target marks lie in $|z|\le1$ or $|z|>R_n$, where $1<R_n\to\infty$. Pass to a subsequence on which the labelled marks in $|z|\le1$ converge. Let $\eta_n$ be integrable rational quadratic differentials of norm one with finite poles only at source marks. Suppose $p_n=(H_n)_*\eta_n$ has $\|p_n\|_1\to1$ and converges locally uniformly off the limiting target marks to an integrable differential $p$ on the component sphere.

Put $\rho_n=\sqrt{R_n}$. If $p\ne0$, then after passage to a subsequence there are a fixed $d\ge1$ and Jordan disks $D_n$, each bounded by a lift of $|z|=\rho_n$, for which $H_n:D_n\to\{|z|<\rho_n\}$ is proper of degree $d$. Each $D_n$ contains at least $d+2$ distinct finite poles of $\eta_n$, hence
\begin{equation}\label{eq:degree}
                            d\le s-2 .
\end{equation}
If all finite target marks lie in $|z|\le1$ for every $n$ (the root chart), then $p=0$.
\end{prop}
\end{samepage}
\begin{proof}
Assume $p\ne0$. Steps~1--3 exclude the root case, proving the last conclusion by contradiction; the remaining steps construct the disks and prove the degree bound in the nonzero case.

\textbf{1. Select an inverse branch by its relative error.}
Write $z=e^w$ and $\mathcal P=\exp^*p$. Since $p$ is nonzero, we can fix disks $U\Subset V\Subset\{\Re w>0\}$ about a point $w_0$ such that $\exp$ is injective near $\overline V$ and $\mathcal P$ has neither a zero nor a pole there. For large $n$, $\exp(\overline V)$ lies in the unmarked annulus $1<|z|<R_n$, and $\mathcal P_n=\exp^*p_n$ converges uniformly to $\mathcal P$ on $\overline V$. In particular $\mathcal P_n$ is nonzero there. Since $\exp(V)$ is simply connected and contains no singular value, each preimage of $e^{w_0}$ determines a branch $\psi_{j,n}$ on $V$ satisfying $H_n\circ\psi_{j,n}(w)=e^w$. The branch images are disjoint and contain no poles of $\eta_n$, since the image of each such pole is marked. We will also use in Step~3 that $p$ is rational with at most simple poles by integrability on the sphere, so its coefficient is $O(z^{-3})$ at infinity.

On $V$ the pushforward is the sum of the branch contributions $b_{j,n}=\psi_{j,n}^*\eta_n$, so $\mathcal P_n=\sum_j b_{j,n}$. Disjointness and change of variables give $\sum_j\int_V|b_{j,n}|\le\|\eta_n\|_1=1$. Lemma~\ref{l:pushforward} gives absolute locally uniform convergence, including for infinitely many branches.

On $U$ define the branch mass $A_{j,n}$ and angular error $\Delta_{j,n}$ by
\[
 A_{j,n}=\int_U |b_{j,n}|,\qquad
 \Delta_{j,n}=\int_U\left(|b_{j,n}|-
       \Re\left(b_{j,n}\frac{\overline{\mathcal P_n}}{|\mathcal P_n|}\right)\right).
\]
Here we use the coefficients in the $w$-coordinate and integrate against Euclidean area. Where $b_{j,n}\ne0$, the integrand for $\Delta_{j,n}$ is $|b_{j,n}|(1-\cos\vartheta)$, where $\vartheta$ is the angle between the branch coefficient and the coefficient of the sum $\mathcal P_n$. Thus each error is nonnegative. Summing the projections onto the direction of $\mathcal P_n$ gives $|\mathcal P_n|$, so $\sum_j\Delta_{j,n}$ is exactly the mass lost under pushforward over $\exp(U)$. The loss over the complement is also nonnegative by the triangle inequality and the area formula. Hence
\[
   0\le\sum_j\Delta_{j,n}\le\|\eta_n\|_1-\|p_n\|_1=1-\|p_n\|_1\longrightarrow0 .
\]
Here $\|p_n\|_1\to1$ is a hypothesis of the proposition, supplied by Proposition~\ref{l:transport} in the application. The total branch mass on $U$ stays bounded below because $\sum_j A_{j,n}\ge\int_U|\mathcal P_n|\to\int_U|\mathcal P|>0$.

Each $A_{j,n}$ is positive since $b_{j,n}\not\equiv0$. The quotient $\sum_j\Delta_{j,n}/\sum_j A_{j,n}$ is therefore a weighted average of the individual ratios $\Delta_{j,n}/A_{j,n}$, with weights $A_{j,n}/\sum_\ell A_{\ell,n}$. For each large $n$ choose a branch whose ratio is at most this average, and denote the branch, its mass and its angular error by $\psi_n,A_n,\Delta_n$. Then
\[
  0\le\frac{\Delta_n}{A_n}
     \le\frac{\sum_j\Delta_{j,n}}{\sum_j A_{j,n}}
     \le\frac{1-\|p_n\|_1}{\int_U|\mathcal P_n|}\longrightarrow0 .
\]

\textbf{2. What a single converging branch determines.}
Suppose that, in affine source coordinates $\lambda_n$, the maps $\Psi_n=\lambda_n\circ\psi_n$ converge locally uniformly on $V$ to a nonconstant map $\Psi$. We will show that, after passing to a subsequence, there is a nonzero rational quadratic differential $T$ on the source sphere such that
\begin{equation}\label{eq:identity}
                 \Psi^*T=c\,\exp^*p\quad\text{on }U,\qquad c>0 .
\end{equation}

We first normalize the source differential so that its mass on the selected branch is one, even when $A_n\to0$. In the new source coordinate this gives
\[
 T_n=A_n^{-1}(\lambda_n^{-1})^*\eta_n,\qquad
 \Psi_n^*T_n=A_n^{-1}\psi_n^*\eta_n,\qquad
 \int_U|\Psi_n^*T_n|=1 .
\]
To obtain a limit, write $T_n=c_n(N_n/W_n)\,\dd z^2$ in reduced form, scaling the coefficient vectors of $N_n$ and $W_n$ separately to have norm one. Pass to a subsequence on which the number $m\le s$ of finite simple poles of $\eta_n$ is constant. Then $\deg W_n=m\le s$ and $\deg N_n\le m-3$ by integrability at infinity. These bounds are preserved by affine coordinate changes. After passing to a subsequence, $N_n,W_n$ converge coefficientwise to nonzero polynomials $N,W$. Since the branch avoids the poles, $W_n\circ\Psi_n$ is zero-free on $V$. Its limit $W\circ\Psi$ is not identically zero because $\Psi$ is nonconstant, so Hurwitz's theorem makes it zero-free on $V$ as well.

The nonconstant limit $\Psi$ of the univalent maps $\Psi_n$ is itself univalent by Hurwitz's theorem. The maps and their derivatives converge uniformly on $\overline U$, so $\Psi_n^*((N_n/W_n)\,\dd z^2)$ converges uniformly there to $\Psi^*((N/W)\,\dd z^2)$. The latter has finite, positive mass on $U$, since $N$ is not identically zero. The unit mass normalization above now shows that $|c_n|$ tends to the reciprocal of this mass. Passing to a subsequence on which the phases converge gives $c_n\to c_0\ne0$. Set $T=c_0(N/W)\,\dd z^2$, cancelling any common factors. We have $\Psi_n^*T_n\to\Psi^*T$ uniformly on $\overline U$ and $\int_U|\Psi^*T|=1$.

We can now use the angular estimate from Step~1. Divide by the target differential and write
\[
        \chi_n=\frac{\Psi_n^*T_n}{\mathcal P_n},\qquad
        \chi=\frac{\Psi^*T}{\mathcal P},\qquad \mathcal P=\exp^*p .
\]
Since $\mathcal P_n\to\mathcal P$ uniformly on $\overline U$ and $\mathcal P$ is nowhere zero there, $\chi_n\to \chi$ uniformly on $\overline U$. In terms of these quotients, the relative angular error from Step~1 is exactly
\[
  \frac{\Delta_n}{A_n}=\int_U (|\chi_n|-\Re \chi_n)\,|\mathcal P_n| .
\]
The left side tends to zero. Uniform convergence therefore gives $\int_U (|\chi|-\Re \chi)\,|\mathcal P|=0$. If the continuous nonnegative function $|\chi|-\Re \chi$ were positive anywhere on $U$, a small disk around that point would give a positive integral. Thus $|\chi|=\Re \chi$ everywhere on $U$, which says that $\chi$ is real and nonnegative. Since $\chi$ is holomorphic on the disk $U$, the open mapping theorem forces it to be constant. Its constant value is positive because $\int_U|\Psi^*T|=1$. This proves \eqref{eq:identity}.

We also need to know that each finite pole of $T$ is approached by poles of $T_n$. Such a pole is a zero of $W$. Choose a small disk around it whose boundary contains no zero of $W$. By coefficientwise convergence and Rouché's theorem, $W_n$ has a zero in this disk for all large $n$. Each such zero is an actual pole of $T_n$ because $N_n/W_n$ is reduced. This argument applies to every pole that remains after cancelling common factors of $N$ and $W$, even when other poles have merged or escaped to infinity.

\textbf{3. Rule out unbounded covering degree and universal covers.}
Let $d_n$ be the covering degree of the selected component over $1<|z|<R_n$, with $d_n=\infty$ for the universal cover. The branch continues to a lift of $w\mapsto e^w$ on the simply connected strip $0<\Re w<\log R_n$. If two points have the same lifted image, their difference is $2\pi i$ times an element of the stabilizer in $\mathbb Z$ of the chosen point in the fiber. The stabilizer is $d_n\mathbb Z$ for finite degree, giving additive period $2\pi i d_n$ in the logarithmic coordinate. The stabilizer of a universal cover is trivial. Thus the lift is univalent after restriction to $|\Im(w-w_0)|<\pi d_n$, with the whole strip allowed when $d_n=\infty$.

Suppose a subsequence has $d_n\to\infty$ in this extended sense. These restricted domains eventually contain every compact subset of the right half-plane. Affinely normalize the lifted maps to have value zero and derivative one at $w_0$. They have a subsequence converging on compact subsets to a univalent map $\Psi$ of the half-plane. Indeed, on any fixed simply connected relatively compact subdomain containing $w_0$ this is the normalized univalent normal-family theorem after a Riemann map \cite[Chapter~2]{Dur}. A diagonal exhaustion and the derivative at $w_0$ give the assertion. Step 2 now applies on $V$, and the meromorphic identity theorem extends \eqref{eq:identity} to the half-plane.

To contradict that identity we compare areas in the flat metrics $|\exp^*p|$ on the half-plane and $|T|$ on the sphere. Since $p(z)=O(|z|^{-3})$, there are $A_0,t_0>0$ with $|p(e^w)|e^{2\Re w}\le A_0e^{-\Re w}$ for $\Re w\ge t_0$. Hence the horizontal ray $s\mapsto w+s$, $s\ge0$, from a point with $\Re w\ge t_0$ has $|\exp^*p|$-length at most $2\sqrt{A_0}\,e^{-\Re w/2}$, and by \eqref{eq:identity} its image under $\Psi$ has $|T|$-length at most $2\sqrt{cA_0}\,e^{-\Re w/2}$. A path of finite $|T|$-length in $\widehat\CC$ converges: it stays away from the poles of $T$ of order at least two, which lie at infinite distance, and elsewhere the $|T|$-metric induces the usual topology. Rays from two points $w,w'$ have the same limit, because the segment from $w+s$ to $w'+s$ has $|\exp^*p|$-length at most $|w-w'|\sqrt{A_0}\,e^{-(\min(\Re w,\Re w')+s)/2}\to0$. Call the common limit $z_*$. Then $\Psi(\{\Re w>t\})$ lies in the $|T|$-ball of radius $2\sqrt{cA_0}\,e^{-t/2}$ about $z_*$ for every $t\ge t_0$. Since $z_*$ is at finite distance it is not a pole of order at least two, so $T$ is integrable near $z_*$, and small balls about $z_*$ have small finite $|T|$-area. On the other hand $\Psi$ is univalent, so the $|T|$-area of $\Psi(\{\Re w>t\})$ equals $c\int_t^\infty\!\int_{-\infty}^{\infty}|p(e^{x+iy})|e^{2x}\,\dd y\,\dd x$. This is infinite, because $y\mapsto|p(e^{x+iy})|$ is $2\pi$-periodic and not identically zero. The contradiction excludes $d_n\to\infty$.

When all finite marks lie in the unit disk, the exterior-cover argument of Section~\ref{sec:closed} applies to $|z|>1$: a finite-degree component would force $H_n$ to be a polynomial, so every component is a universal cover. Restriction to $1<|z|<R_n$ keeps the same stabilizer because the annulus carries the same fundamental group. Thus $d_n=\infty$, contradicting the preceding half-plane area argument and proving $p=0$ in the root case. We continue under $p\ne0$ in the remaining case.

\textbf{4. Identify the selected finite disk and cap its cover.}
Step~3 lets us pass to one fixed finite degree $d$. Over $|z|=\rho_n$, the selected annular component contains an analytic Jordan curve mapping as a $d$-fold cover. By the compact-lift disk argument in Section~\ref{sec:closed}, its bounded side $D_n$ maps properly onto $\{|z|<\rho_n\}$ with degree $d$. The part of the selected annular component below the circle is connected and locally lies on the interior side, so it lies in $D_n$. In particular it contains the chosen branch on $V$ for all large $n$.

On a collar of the circle, classification of annulus covers gives a holomorphic coordinate $\xi$ with $H_n=\xi^d$. This follows by lifting the logarithm, for which the subgroup of additive periods is $2\pi i d\mathbb Z$, and exponentiating one $d$th of the lifted logarithm. Glue the exterior in the $\xi$-plane, together with infinity, along this open collar. The resulting compact Riemann surface is topologically a sphere. Its holomorphic map to the sphere has one pole of order $d$. Uniformization with that point as infinity gives a polynomial $P_n$ and a conformal isomorphism
\[
 h_n:\Omega_n=P_n^{-1}(\{|z|<\rho_n\})\longrightarrow D_n,\qquad H_n h_n=P_n .
\]
Choose the input affine coordinate to make $P_n$ monic and centered, taking $P_n(z)=z$ for $d=1$. The target coordinate stays fixed. Every finite critical value of $P_n$ has modulus at most one. A critical point in the original disk maps to a marked value there, and the added exterior is unbranched at finite points.

\textbf{5. Compactness in the capped coordinate.}
Let $K_n=P_n^{-1}(\overline{\mathbb D})$. Each component of $P_n^{-1}(\widehat\CC\setminus\overline{\mathbb D})$ maps openly and properly onto the exterior disk, hence contains a preimage of infinity. There is just one such point, so the preimage is connected. Removing infinity leaves a connected finite unbranched cover of a punctured disk, since every finite critical value is in $\overline{\mathbb D}$. Filling the puncture gives a disk. Its complement $K_n$ is consequently connected. Logarithmic capacity is defined by the constant term in the exterior Green function, $g(z)=\log|z|-\log\operatorname{cap}(K)+o(1)$ at infinity. Here $g=d^{-1}\log|P_n|$, so monicity gives capacity one, and a continuum of capacity one has diameter at most four \cite[Ch.~5]{Ran}. All roots lie in $K_n$. Their mean is zero since $P_n$ is centered, so each root is within four of zero by averaging its differences from the other roots. This bounds the coefficients.

Take $P_n\to P_*$ coefficientwise, still of degree $d$. Every closed disk in $\CC$ lies in $\Omega_n$ for all large $n$, since $P_n$ is bounded on that disk and $\rho_n\to\infty$. Define $\lambda_n(z)=(z-h_n(0))/h_n'(0)$ for these indices. The injections $\widehat h_n=\lambda_n h_n$ converge locally uniformly on the plane to the identity. Indeed every subsequence has, by univalent normality on an exhaustion, a further limit that is an entire injection with value zero and derivative one at zero. Such an injection is the identity. Moreover each compact set of the normalized \emph{actual} source plane is in $\lambda_n(D_n)$ eventually. For a disk $|y|\le M$ use Rouché's theorem on $|z|=M+1$ for $\widehat h_n(z)-y$ compared with $z-y$.

\textbf{6. Continue the identity and count poles in the original source disk.}
On $V$ let $g_n=h_n^{-1}\psi_n$. Then $P_n(g_n(w))=e^w$. The roots of these equations are uniformly bounded on compact subsets of $V$ by the monic polynomial root bound. A subsequential limit $g$ satisfies $P_*(g(w))=e^w$. Thus $g$ is nonconstant, its derivative is nonzero, and it is univalent since $\exp$ is injective on $V$. The normalized \emph{same} branches $\lambda_n\circ\psi_n=\widehat h_n\circ g_n$ converge to $g$ on compact subsets. Applying step 2 in these affine source coordinates yields $g^*T=c\,\exp^*p=c\,g^*(P_*^*p)$ on $U$. Invert $g$ there and use meromorphic continuation of the two rational differentials on the sphere. This proves $T=cP_*^*p$ globally.

A nonzero rational quadratic differential has divisor degree $-4$. The poles of $p$ are at most simple, so at least three distinct poles $v_1,v_2,v_3$ are finite. Their preimages under $P_*$ have total multiplicity $3d$. A ramified point of local degree $e\ge2$ accounts for $e\le2(e-1)$ of these multiplicities. Summing over the three fibers uses at most the total finite ramification $d-1$, so at least $3d-2(d-1)=d+2$ distinct unramified points remain. In a local coordinate where the map is $z\mapsto z^e$, pullback takes an order-$m$ differential to one of order $em+2(e-1)$. Over a simple pole this is $e-2$, which is $-1$ for each of the unramified preimages just counted. These are finite poles of $T$. Choose disjoint closed disks about them with boundaries avoiding the zeros of $W$ from step 2. For all sufficiently large $n$, each disk lies in $\lambda_n(D_n)$ by step 5 and contains a pole of $T_n$ by step 2. Multiplication by $A_n^{-1}$ leaves poles in place, so applying $\lambda_n^{-1}$ produces at least $d+2$ distinct poles of $\eta_n$ in the selected $D_n$. All are marked. This proves \eqref{eq:degree}. The count permits all finite ramification to occur over the selected three poles. For example, for $P(z)=z^2$ and $p(w)=\dd w^2/[w(w-1)(w+1)]$, one gets $P^*p=4\,\dd z^2/(z^4-1)$. The ramified points over $0$ and $\infty$ are regular for the pullback and the four unramified poles attain $d+2$. This illustrates the local pole count.
\end{proof}

\section{Accounting for the whole measure}\label{sec:mass}
Assume that $f$ has no Levy cycle relative to $A^+$ and its pullback has no fixed point. Put $s=|A|\ge3$ and $r=s-2$, the largest possible size of a system of disjoint, mutually nonisotopic essential simple curves on this marked sphere. Fix once and for all
\[
       b=\min\{m\ge1 : m\text{ is an integer and }2^m>r\},\qquad k=rb+1.
\]
This depth makes a chain of lifts of total degree at most $r$ contain a block of $r$ degree-one lifts; once their classes are known to be essential and compatible, repetition gives a Levy cycle. The last lift ensures essentiality at depth $k-1$.

Use the positive measure $\mu$ and the number $E$ of Section \ref{sec:transport} for exactly this depth. Thus the order is to fix $o,G,C$, then the prescribed $k$ and its representatives, then $E$ and $\mu$, and only afterwards an arbitrary sequence $t_n\to\infty$. The contradiction below uses this choice of $k$; the measure construction in Section~\ref{sec:transport} was available at every fixed depth. Neither $G$ nor $\mu$ is reselected after the sequence or the cluster charts are chosen. Let $q_n=q_{t_n}$ and $P_n=e^{2t_n}(H_{t_n})_*r_{t_n}$, and let $\mu_n=(\phi_{t_n})_*\mu$. All are on the time-$t_n$ target. The maps $H_{t_n}$ are transcendental entire, as required by Proposition~\ref{p:degree}. We use these same symbols in affine target coordinates. In particular
\begin{equation}\label{eq:tv}
  \|\,|P_n|-\mu_n\|_{\rm TV}\to0,\qquad
  \|P_n\|_1\to E,\qquad \|q_n\|_1=1,\qquad
  \mu_n(U)\le\omega\bigl(|q_n|(U)\bigr).
\end{equation}
For densities $p_n,s_n$ relative to one area element, the variation measure of their difference is $|p_n-s_n|\,\dd A$. Its mass bounds the error on every measurable set, proving the first assertion from \eqref{eq:pushforward}. Thus positive limiting $|P_n|$-mass forces positive $\mu_n$-mass; the bound by $\omega$, with $\omega(\delta)\to0$, then forces positive $|q_n|$-area. We use this on component disks and on moving annular bands.

\subsection{Components in affine coordinates}\label{sec:components}
\subsection*{An example with three scales}\label{bg:cluster}
\addcontentsline{toc}{subsection}{An example with three scales}
For marks at $0,\epsilon_n^2,\epsilon_n,1$ with $\epsilon_n\downarrow0$, the charts $z$, $z/\epsilon_n$, and $z/\epsilon_n^2$ successively separate the nested clusters
\[
 \{0,\epsilon_n^2\}\subset\{0,\epsilon_n^2,\epsilon_n\}
 \subset\{0,\epsilon_n^2,\epsilon_n,1\}.
\]
The unmarked annuli $2\epsilon_n^2<|z|<\epsilon_n/2$ and $2\epsilon_n<|z|<1/2$ have moduli tending to infinity. For $p_n=u_n(z)\dd z^2$, the middle-chart coefficient is $\epsilon_n^2u_n(\epsilon_n Z)$; area is unchanged by this coordinate change. Limits are taken locally uniformly off the limiting marks. The annuli escape compact sets in both adjacent charts, so they can retain mass even when all component limits vanish \cite[\S5]{HSS}.

Start with the cluster $A$ of all finite target labels as the root of a tree. For each cluster $B$ of at least two labels, choose the position $c_{B,n}$ of a fixed label as center and its Euclidean diameter $\rho_{B,n}>0$ as scale. In the chart $z_B=(z-c_{B,n})/\rho_{B,n}$ the labels lie in the closed unit disk. Pass to a subsequence on which they converge. Their limiting diameter is one. The labels with each distinct limiting position form one child of $B$. An \emph{edge} $e=(\hat B,B)$ joins a parent $\hat B$ to an immediate child $B$, one step in this partition. Repeat the construction in each child containing at least two labels. Each split strictly decreases size, so the tree is finite and its leaves are singletons. If all limits in the root chart are distinct, every child of $A$ is already a leaf.

For a nonsingleton child $B$, $\rho_{B,n}/\rho_{\hat B,n}\to0$, giving a long annulus on its incoming edge. Singleton leaves need no further annuli: a small fixed disk about each singleton limit contains only that mark, and the simple-pole estimate in Point~1 of Proposition~\ref{l:cylinder-mass} makes its mass vanish. Thus the annulus entering a two-label cluster can carry mass even though its two leaf disks do not.

For a nonroot cluster $B$ of at least two labels and a label $a\notin B$, use the first ancestor where $a$ and $B$ occupy different children. Their separation is bounded below in that ancestor's chart while $B$ has vanishing relative diameter. Thus every outside label escapes to infinity in $B$'s chart. At a nonroot chart one may therefore choose $R_n$ smaller than the least modulus of an outside label and tending to infinity. The resulting round annulus separates precisely the child from the other finite labels. Its core is well defined at each index $n$, even if the marked isotopy classes of these cores vary with $n$.

\begin{figure}[ht]\centering
\begin{tikzpicture}[scale=0.88]
\begin{scope}
\draw[pine!40] (0,0) circle (2.6);
\node[font=\scriptsize,pine] at (0,2.85) {parent chart $z_{\hat B}$ (here the root)};
\fill (1.5,-1.2) circle (1.6pt); \fill (-1.6,1.0) circle (1.6pt);
\fill[mist,even odd rule] (-0.4,-0.2) circle (1.1) (-0.4,-0.2) circle (0.25);
\draw[pine] (-0.4,-0.2) circle (1.1) (-0.4,-0.2) circle (0.25);
\begin{scope}[shift={(-0.4,-0.2)},scale=0.28]
\foreach \p in {(-0.45,0.4),(0.55,0.4)}{\node[circle,fill,inner sep=0pt,minimum size=2.2pt] at \p {};}
\node[diamond,fill=plum,inner sep=0pt,minimum size=3.8pt] at (0,0) {};
\end{scope}
\node[font=\scriptsize,align=center] at (-0.4,-1.7) {edge annulus $e=(\hat B,B)$:\\ $16\rho_{B,n}<|z-c_{B,n}|<\vartheta_{\hat B}\rho_{\hat B,n}$};
\draw[dashed,-{Stealth}] (-0.15,-0.05)--(3.35,0.6);
\end{scope}
\begin{scope}[shift={(5.6,0.3)}]
\draw[pine!40] (0,0) circle (2.0);
\node[font=\scriptsize,pine] at (0,2.25) {child chart $z_B=(z-c_{B,n})/\rho_{B,n}$};
\fill (-0.45,0.4) circle (1.6pt); \fill (0.55,0.4) circle (1.6pt);
\node[diamond,fill=plum,inner sep=0pt,minimum size=4.6pt] at (0,0) {};
\node[below=2pt,font=\scriptsize,plum] at (0,0) {$z_B=0$};
\draw[pine!70,dashed] (0,0) circle (1.0);
\node[font=\scriptsize] at (0.0,-1.25) {$|z_B|\le1$};
\end{scope}
\begin{scope}[shift={(10.2,0)}]
\draw[thick] (0,0) rectangle (3.0,1.4);
\fill[sand] (0.3,0) rectangle (2.7,1.4);
\node[font=\scriptsize,align=center] at (1.5,0.7) {$\hat P_n\approx a(\hat P_n)$\\ $\hat q_n\approx a(\hat q_n)$};
\node[font=\scriptsize,below] at (0,0) {$0$}; \node[font=\scriptsize,below] at (3.0,0) {$L$};
\node[font=\scriptsize,above,align=center] at (1.5,1.45) {$w=\log\frac{z-c_{B,n}}{\rho^{\rm in}_n}$};
\draw[-{Stealth}] (3.2,0)--(3.2,1.4) node[midway,right,font=\scriptsize]{$2\pi$};
\end{scope}
\end{tikzpicture}
\caption{Rescaling clusters of target marks. The diamond $c_{B,n}$ becomes the origin in the child chart. Each nonsingleton child gives a long edge annulus. Lemma~\ref{l:components} eliminates component mass; Proposition~\ref{l:cylinder-mass} locates the remaining mass on these annuli and constructs $q_n$-flat subcylinders of diverging modulus on every positive-mass edge.}
\label{fig:clusters}
\end{figure}
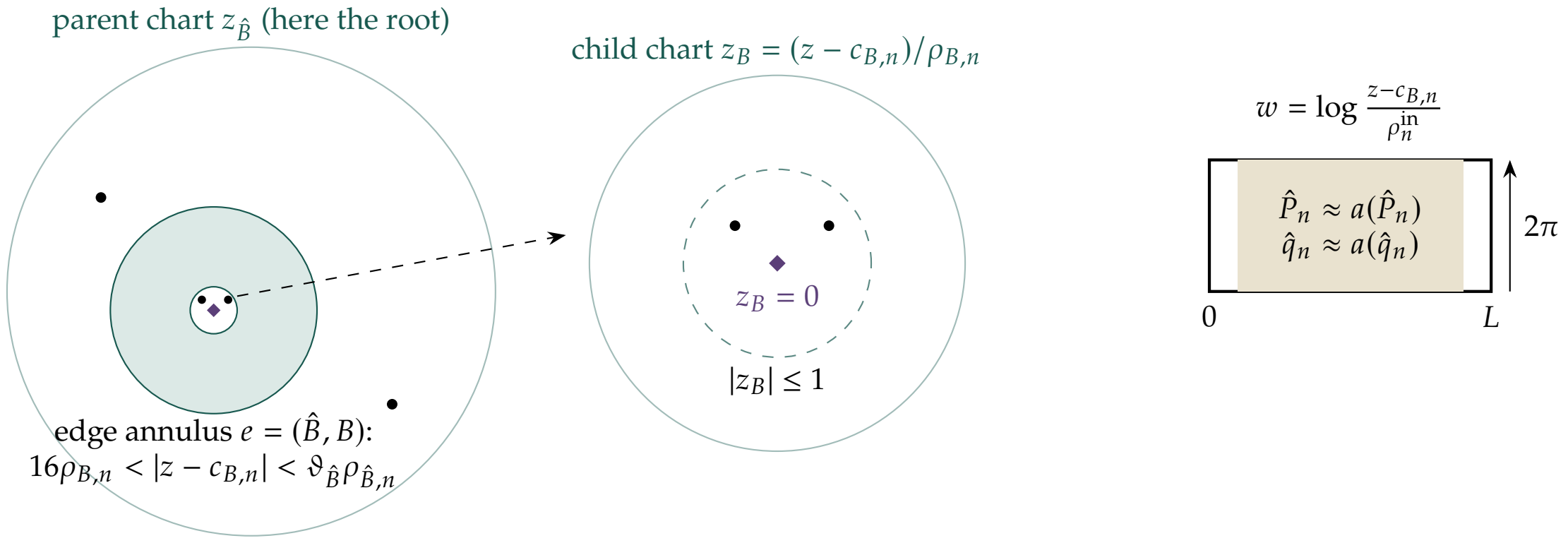

Recall that $P_n=P_{t_n}$, where $P_t$ is defined just above~\eqref{eq:pushforward} and the pullback depth $k$ stays fixed as $t_n\to\infty$.
Lemma~\ref{l:pushforward} shows that the differentials $P_n=e^{2t_n}(H_{t_n})_*r_{t_n}$ are holomorphic off target marks, including when the cover has infinitely many sheets, and have at most simple poles at the marks and infinity. In each of the finitely many cluster charts, the area bounds for $P_n,q_n$ and the area mean-value inequality on disks avoiding the limiting labels give locally uniform bounds on their coefficients. Montel's theorem and simultaneous extraction give locally uniform coefficient limits off those labels and infinity. Fatou's lemma shows that the area of each component limit is at most the liminf of the areas of the approximating differentials. These finite-area limits are rational differentials with at most simple poles at the omitted points, by the local integrability argument in the proof of the same lemma, applied also in the coordinate $1/z$ at infinity.

On this common subsequence either some component limit of $P_n$ is nonzero or all such limits vanish. Lemma~\ref{l:components} excludes the first alternative by applying Proposition~\ref{p:degree} with the normalization specified in its proof and then using the no-Levy hypothesis at the prescribed depth. In the second alternative Proposition~\ref{l:cylinder-mass} locates the mass on the annuli between clusters and supplies the cylinders used in Section~\ref{sec:persistence}. Every further extraction retains the same original measure $\mu$.

\emph{How classes on different surfaces are compared.} Suppose $\gamma_0$ is a curve on $X_0$ and $\gamma_j$ is a curve on $X_j$. The markings $m_0:S\to X_0$ and $m_j:S\to X_j$ are homeomorphisms from the same reference sphere $S$ with marked set $A^+$. Thus $m_0^{-1}(\gamma_0)$ and $m_j^{-1}(\gamma_j)$ are curves on one sphere. Their geometric intersection number is the minimum number of crossings among transverse representatives of the two classes on $S$, with isotopies fixing $A^+$.

For example, take four marks $1,2,3,\infty$. Draw a simple loop enclosing precisely $1,2$ on the first sphere and one enclosing precisely $1,3$ on the second. Bring both loops back to $S$ and consider the disks on the sides away from $\infty$. Both contain $1$, and each contains a mark that the other excludes. Disjoint boundaries would force the disks to be disjoint or nested. The same argument applies after any isotopy fixing the marks, so the two classes have positive intersection number.

The proof below first uses \eqref{eq:flat} to put short lifted classes on the common target. A crossing estimate then gives zero intersection; lifting disjoint representatives supplies the remaining pairwise disjointness needed for the bound $s-2$.

\begin{samepage}
\begin{lemma}\label{l:components}
Let $(f,A)$ be a transcendental marked entire Thurston map with $s=|A|\ge3$, no fixed point of its pullback, and no Levy cycle relative to $A^+$. Put $r=s-2$, $b=\min\{m\in\mathbb N:2^m>r\}$ and $k=rb+1$. Fix $o$ and a nonzero $G$ with $\Ext_oG=1$ satisfying \eqref{eq:displacement} with a single $C$, and the corresponding depth-$k$ objects from Section~\ref{sec:transport}. For any sequence $t_n\to\infty$, set $P_n=e^{2t_n}(H_{t_n})_*r_{t_n}$. On every subsequence on which the finite cluster tree and the locally uniform component limits described above exist, the component limit of $P_n$ vanishes in every cluster chart.
\end{lemma}
\end{samepage}
\begin{proof}
\textbf{1. A nonzero limit differential gives a proper map of bounded degree on a disk.}
Recall that ``in a chart $B$'' means in the rescaled coordinate $z_B=(z-c_{B,n})/\rho_{B,n}$ for one of the clusters constructed in the paragraph beginning ``Start with the cluster $A$''.
Suppose the component limit of $P_n$ is nonzero in a chart $B$. Choose a small closed disk away from the limiting marks on which the coefficient of this component limit in the $z_B$ coordinate is bounded away from zero. Uniform convergence on the disk then gives it positive limiting $|P_n|$-mass. By \eqref{eq:tv} its $\mu_n$-mass stays bounded below, and then so does its $|q_n|$-area by \eqref{eq:omega}. Therefore the locally uniform limit $q_B$ in this chart is nonzero.

Apply Proposition \ref{p:degree} with the normalized source differential
\[
 \eta_n=\frac{r_{t_n}}{\|r_{t_n}\|_1},\qquad
 (H_{t_n})_*\eta_n=\frac{P_n}{E_{t_n}},\qquad
 E_{t_n}=e^{2t_n}\|r_{t_n}\|_1\longrightarrow E>0.
\]
Thus $\|\eta_n\|_1=1$, and \eqref{eq:tv} gives $\|(H_{t_n})_*\eta_n\|_1=\|P_n\|_1/E_{t_n}\to1$. Division by $E_{t_n}\to E>0$ also preserves nonvanishing of the component limit. In the root chart all finite target marks lie in $|z_A|\le1$, so the last sentence of the proposition forces the component limit of $P_n/E_{t_n}$, and hence that of $P_n$, to vanish there. If all children of $A$ are singletons, the root is the only chart and this proves the lemma.

It remains to consider a proper cluster $B$. Its labels lie in $|z_B|\le1$ and all other labels are outside $|z_B|\le R_n$ for some $R_n\to\infty$. Set $R'_n=\sqrt{R_n}$. Along a further subsequence we have a disk over $|z_B|<R'_n$ of fixed degree $d\le r$ under $H_{t_n}$ with at least $d+2$ marked source poles. The target annulus $2<|z_B|<R'_n$ has no marked point. Its core represents a class $\alpha_{0,n}$ which may change with $n$.

\textbf{2. Follow the actual lifted components.}
Recall from the paragraph immediately before ``Reading a Thurston matrix'' in Section~\ref{sec:closed} that $F_j=m_{j-1}\circ f\circ m_j^{-1}:X_j\to X_{j-1}$ represents one application of $f$ in the markings at depths $j$ and $j-1$. Take $H=F_1\circ\cdots\circ F_k$ in Section~\ref{sec:transport}. For $t=t_n$ put $\chi_{0,t}=\phi_t$. On $X_j$ uniformize the structure pulled back from $\phi_t$ by $F_1\circ\cdots\circ F_j$, with marking $\chi_{j,t}$. Use the same final normalization $\chi_{k,t}=\psi_t$ as in Section \ref{sec:transport}. Then $F_{j,t}=\chi_{j-1,t}F_j\chi_{j,t}^{-1}$ is entire by construction, and cancellation gives $F_{1,t}\circ\cdots\circ F_{k,t}=\phi_t H\psi_t^{-1}=H_t$. The affine target chart for $B$ changes only the coordinate on the range of $F_{1,t}$. Thus we can follow the actual boundary of the selected disk and its marked poles through these factors.

Each successive image of the boundary is compact and open in its component over the target circle, hence is the full circle. Each bounds a Jordan disk on which the corresponding map is proper by the argument principle. Each disk maps onto its preceding one, with degree equal to the boundary degree. Write $D_{j,n}$ for the disk at depth $j$, $\alpha_{j,n}$ for its marked boundary class, and $d_j$ for its degree over depth $j-1$, suppressing $n$ in the latter notation. Their product is $d$. For $0\le j<k$, the proper map from the final disk $D_{k,n}$ to $D_{j,n}$ has degree
\[
 e_j=d_{j+1}\cdots d_k\le d.
\]
Each fiber has at most $e_j$ distinct points, since its total local-degree multiplicity is $e_j$. The images of a chosen set of $d+2$ final marked poles therefore comprise at least $\lceil(d+2)/e_j\rceil\ge2$ distinct points. They remain marked by the forward-image condition on the representatives. Thus the bounded side of $\alpha_{j,n}$ contains at least two finite marks.

For its other side, use the further compact lift from depth $j+1$. Section~\ref{sec:closed} shows that a curve enclosing all finite marks has no compact lift under a transcendental entire representative: a finite-degree exterior component would force a pole at infinity and hence a polynomial. Consequently $D_{j,n}$ omits a finite mark, which together with $\infty$ supplies two marks on the unbounded side. Thus $\alpha_{j,n}$ is essential for every $j<k$. The pigeonhole argument uses these depths $0,\ldots,k-1$; the final disk at depth $k$ may contain all finite source marks. After the subsequence and disk construction, the conclusions hold for all sufficiently large $n$, on one common tail for the finitely many depths.

At depth $j$, continue a parallel copy of $\partial D_{j,n}$ inward from its unbranched collar over the target annulus $2<|z_B|<R'_n$. This selects the component of the preimage of that annulus whose closure meets $\partial D_{j,n}$. It lies in $D_{j,n}$, contains no source mark, has core class $\alpha_{j,n}$, and maps onto the annulus with degree $d_1\cdots d_j\le d$ (empty product one). Indeed the map on $D_{j,n}$ is proper, and all singular values and images of source marks are outside the target annulus. The modulus is therefore the target modulus divided by this degree, giving
\[
 \Ext_{\sigma^j g(t_n)}\alpha_{j,n}
       \le\frac{2\pi d}{\log(R'_n/2)}\longrightarrow0
       \qquad(0\le j<k).
\]
With $C_*=\max_{0\le i<k}C_i$, estimate~\eqref{eq:flat} gives simultaneously
\[
 \ell_{q_n}(\alpha_{j,n})
       \le C_*\sqrt{\frac{2\pi d}{\log(R'_n/2)}}\longrightarrow0
       \qquad(0\le j<k).
\]

\textbf{3. Compare curve classes at a single index.}
Choose $R_0>2$ so that the closed band $R_0\le|z_B|\le2R_0$ avoids all finite limiting marks and zeros of the nonzero limit $q_B$. For all large $n$ the coefficient in $q_n=u_n(z_B)\,\dd z_B^2$ has modulus at least some $c>0$ throughout this band. Its boundary circles represent $\alpha_{0,n}$. Every representative of a class intersecting $\alpha_{0,n}$ essentially has a subarc crossing the band from one boundary to the other, and hence has flat length at least
\[
 \int\sqrt{|u_n|}\,|\dd z_B|
       \ge\sqrt c\int\left|\dd(|z_B|)\right|\ge\sqrt c\,R_0.
\]
This applies to the infimum over representatives. Comparing with the preceding length bound gives $i(\alpha_{j,n},\alpha_{0,n})=0$ for every $j<k$ on one common tail of the sequence.

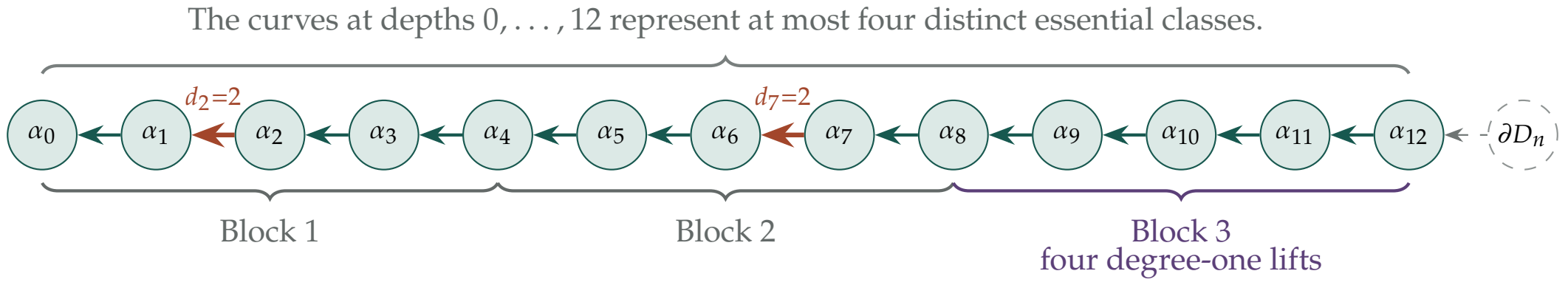
\begin{figure}[!ht]\centering
\begin{tikzpicture}[x=10.2mm]
\foreach \i in {0,...,12}{\node[circle,draw=pine,fill=mist,inner sep=0.8pt,font=\tiny,minimum size=6.2mm] (v\i) at (\i,0) {$\alpha_{\i}$};}
\node[circle,draw=ink!50,dashed,inner sep=0.8pt,font=\tiny,minimum size=6.2mm] (v13) at (13,0) {$\partial D_n$};
\foreach \i/\d in {0/1,1/2,2/1,3/1,4/1,5/1,6/2,7/1,8/1,9/1,10/1,11/1}{
  \pgfmathtruncatemacro{\j}{\i+1}
  \ifnum\d=1 \draw[pine,thick,{Stealth}-] (v\i)--(v\j); \else \draw[rust,very thick,{Stealth}-] (v\i)--(v\j) node[midway,above=2pt,font=\tiny,rust]{$d_{\j}{=}2$};\fi}
\draw[ink!65,dashed,{Stealth}-] (v12)--(v13);
\foreach \start/\finish/\block/\colour in {0/4/1/ink!70,4/8/2/ink!70,8/12/3/plum}{
  \draw[decorate,decoration={brace,amplitude=4pt,mirror},thick,\colour] (\start,-0.43)--(\finish,-0.43)
    node[midway,below=5pt,font=\scriptsize,\colour]{Block \block};}
\node[font=\scriptsize,plum] at (10,-1.14) {four degree-one lifts};
\draw[decorate,decoration={brace,amplitude=4pt},thick,ink!60] (0,0.55)--(12,0.55)
  node[midway,above=5pt,font=\scriptsize,ink!70]{The curves at depths $0,\ldots,12$ represent at most four distinct essential classes.};
\end{tikzpicture}
\caption{The block argument for $r=4$, $b=3$, $k=13$. Arrows point from lifts to images. Since $d\le4<2^3$, one of the three four-edge blocks has only degree-one lifts (here the third). Its five vertices repeat a class; a shortest return gives a Levy cycle. The dashed final lift ensures that $\alpha_{12}$ is essential.}
\label{fig:chain}
\end{figure}

Now fix one such $n$. For $0<i<j<k$, compare $\alpha_{0,n}$ and $\alpha_{j-i,n}$. They admit disjoint representatives relative to all marks, with parallel copies permitted for equal classes. Under the topological map $f^i$ these have components in the prescribed classes $\alpha_{i,n}$ and $\alpha_{j,n}$ by the marking convention of Section~\ref{sec:closed}, which applies equally to the time-$t$ factors $F_{j,t}$: transported to the reference surface, $\alpha_{i,n}$ is a component of $f^{-i}(\alpha_{0,n})$ and $\alpha_{j,n}$ a component of $f^{-i}(\alpha_{j-i,n})$, for the same topological $f^i$. The components are disjoint. Pairwise zero intersection gives simultaneously disjoint representatives, as in Section~\ref{sec:excluding}, so among $\alpha_{0,n},\ldots,\alpha_{k-1,n}$ there are at most $r$ distinct classes. The argument is applied separately at each such index $n$.

Divide the first $k-1=rb$ edges into $b$ consecutive blocks of $r$ edges each, as in Figure~\ref{fig:chain}. If every block contained an edge of degree at least two, the product of the degrees over these blocks would be at least $2^b>r$. Yet this product is $d_1\cdots d_{k-1}\le d\le r$ by Step~1. Hence one block consists entirely of degree-one edges. Its $r+1$ vertices represent at most $r$ distinct classes. Within this block choose two vertices of the same class with the fewest edges between them. The classes along this segment are distinct up to the return, they have disjoint representatives, and each edge gives a degree-one preimage under $f$. The segment therefore gives a Levy cycle, contradicting the hypothesis. This excludes the hypothesized nonzero component limit.
\end{proof}

\begin{remark}
For $s=3$ ($r=1$, $b=1$, $k=2$) the pigeonhole count is not needed. By \eqref{eq:degree}, $d=1$, so the selected disk contains at least $d+2=3$ marks, that is, all finite marks. The map from depth two to depth one has degree one, so the depth-one disk also contains all finite marks while having a compact lift, contradicting Section~\ref{sec:closed}. For $s=4$ the block argument uses $r=2$, $b=2$ and $k=5$.
\end{remark}

\subsection{Mass and flat cylinders on the annuli between clusters}\label{sec:annuli}
\emph{A differential whose mass lies between scales.}
On $1<|z|<e^L$ consider the local model
\[
       p_L=\frac{1}{2\pi L}\frac{\dd z^2}{z^2}
           =\frac{1}{2\pi L}\,\dd w^2,\qquad w=\log z .
\]
The cylinder has period $2\pi i$, area one and modulus $L/(2\pi)$. An end band of fixed width $b_0$ has mass $b_0/L\to0$, whereas $uL<\Re w<vL$ has mass $v-u$ for $0<u<v<1$. Thus coefficients may vanish in every fixed end chart while the middle retains mass.

In the application $P_n$ determines the limiting mass, whereas $q_n$ defines the flat cylinders. Total variation and the bound by $\omega$ transfer positive mass to positive $q_n$-area. The Laurent estimate then yields a nonzero dominant constant on the middle band, where Point~3 constructs the actual cylinder.

\emph{Edge annuli.} Use one edge $e=(\hat B,B)$ for each nonsingleton child $B$ of a parent $\hat B$. For each parent choose $0<\vartheta_{\hat B}<1/100$ smaller than one hundredth of the minimum separation of its distinct child limits. Denote by $\mathcal A_{e,n}$ the round annulus centered at $c_{B,n}$ with inner radius $\rho^{\rm in}_n=16\rho_{B,n}$ and outer radius $\rho^{\rm out}_n=\vartheta_{\hat B}\rho_{\hat B,n}$. Its logarithmic width $L=L_{e,n}=\log(\rho^{\rm out}_n/\rho^{\rm in}_n)$ tends to infinity. For large $n$ the closed enlarged annulus with radii $\rho^{\rm in}_n/\exp(1)$ and $\exp(1)\rho^{\rm out}_n$ misses all marked points. Sibling outer disks are disjoint by separation in the parent chart. The outer disks for the children of $B$ all lie strictly within the inner disk of the edge entering $B$, by the choice of the factor $16$. These facts give disjoint annuli, even for edges at different depths.

\begin{prop}[Mass carried by flat cylinders]\label{l:cylinder-mass}
Under the assumptions of Section~\ref{sec:mass}, let $t_n\to\infty$. After passing to a subsequence, each edge annulus $\mathcal A_{e,n}$ has limiting $\mu_n$-mass $m_e\ge0$, and $\sum_e m_e=E$. For each fixed $0<\epsilon<1/4$ and all sufficiently large $n$, the annulus for each edge with $m_e>0$ contains an embedded $q_n$-flat cylinder $C_{e,n}$ such that
\[
    \Mod(C_{e,n})\longrightarrow\infty,
    \qquad \mu_n(C_{e,n})\longrightarrow(1-4\epsilon)m_e.
\]
The cylinders belonging to different edges are disjoint.
\end{prop}
\begin{proof}
\textbf{1. Partition into annuli and complementary pieces.}
Use the disjoint edge annuli $\mathcal A_{e,n}$ defined above.

Up to boundary circles, the complementary piece at a nonroot cluster $B$ is its radius-$16\rho_{B,n}$ disk minus the outer disks at its nonsingleton children. At the root use the whole sphere minus those outer disks. We show that the $|P_n|$-mass of each complementary piece tends to zero. In the component chart, each removed outer disk contains a fixed smaller disk about the corresponding nonsingleton child limit. Set aside a small fixed disk around each singleton limit and, at the root, a neighborhood of infinity. The rest of the piece lies in a fixed compact set away from all limiting marks. By Lemma~\ref{l:components} the coefficient of $P_n$ tends uniformly to zero there, so its integral in absolute value tends to zero as well.
On one of the singleton disks $D$, write $P_n=p_n(z)\,\dd z^2$ and let $\xi_n$ be its unique marked point for large $n$. Since $p_n$ has at most a simple pole at $\xi_n$, the function $\chi_n(z)=(z-\xi_n)p_n(z)$ extends holomorphically across $\xi_n$. On the fixed boundary $\partial D$, Lemma~\ref{l:components} gives $p_n\to0$ uniformly, and $|z-\xi_n|$ is uniformly bounded there. Hence $\chi_n\to0$ uniformly on $\partial D$, and the maximum modulus principle gives $\sup_D|\chi_n|\to0$. If $D$ has radius $r_D$, then $D\subset\{|z-\xi_n|<2r_D\}$, and polar integration about $\xi_n$ gives
\[
   \int_D |p_n(z)|\,\dd A(z)
   \le \sup_D |\chi_n| \int_D \frac{\dd A(z)}{|z-\xi_n|}
   \le 4\pi r_D\,\sup_D |\chi_n| \longrightarrow 0 .
\]
Thus the simple-pole factor has a uniformly bounded integral even as $\xi_n$ moves, with $r_D$ fixed. At the root the coordinate $1/z$ gives the same estimate around infinity. There are only finitely many complementary pieces, so their total $|P_n|$-mass tends to zero. Since $\|P_n\|_1\to E$ by \eqref{eq:tv}, the long annuli together carry mass tending to $E$. The decomposition is related to Hubbard--Schleicher--Shishikura \cite[Section~4]{HSS}.

\textbf{2. Identify all the mass on each edge.}
Use $w=x+iy=\log((z-c_{B,n})/\rho^{\rm in}_n)$ with imaginary period $2\pi i$, so $0<x<L$. Write $P_n=\hat P_n(w)\,\dd w^2$ and $q_n=\hat q_n(w)\,\dd w^2$ on this edge. The coefficients are holomorphic on a neighborhood of the closed enlarged annulus. The buffered circles $x=-1,L+1$ lie in compact regular sets of the child and parent charts, respectively. In either chart with coordinate $Z$, the change to $w$ multiplies the differential coefficient by $(Z-Z(c_{B,n}))^2$, uniformly bounded on its buffered end. Hence $|\hat P_n|\le M_n=o(1)$ on both circles, and $|\hat q_n|\le M=O(1)$ there.

For a periodic coefficient $\kappa(w)=\sum_{\ell\in\mathbb Z}a_\ell e^{\ell w}$ bounded by $M'$ on those circles, Cauchy's formula gives $|a_\ell|\le M'e^{-\ell(L+1)}$ and $|a_{-\ell}|\le M'e^{-\ell}$ for $\ell>0$. Thus on $0\le x\le L$,
\begin{equation}\label{eq:laurent}
        |\kappa(w)-a_0|\le C_0 M'\bigl(e^{-x}+e^{-(L-x)}\bigr)
\end{equation}
with an absolute geometric-series constant $C_0$. For example the positive modes have sum of absolute values bounded by $M'e^{-(L+1-x)}/(1-e^{-(L+1-x)})$, since each mode is bounded on the outer buffered circle. The negative modes give the analogous estimate from $x=-1$. Write $a(\kappa)=a_0$ for the constant Laurent coefficient. It is the angular mean of $\kappa$ on every circle of fixed $x$, so $2\pi L|a(\kappa)|\le\int_{\rm edge}|\kappa|\,\dd x\,\dd y$.

Fix an edge $e$ and put $a_{e,n}=a(\hat P_n)$ on this edge. The constant differential $a_{e,n}\,\dd w^2$ has mass $2\pi L|a_{e,n}|$, since a fundamental rectangle has width $L=L_{e,n}$ and height $2\pi$. To estimate the error for $P_n$, integrate \eqref{eq:laurent}. Each end exponential has integral $1-e^{-L}$ in $x$, so
\[
  \int_0^L\!\int_0^{2\pi}|\hat P_n(x+iy)-a_{e,n}|\,\dd y\,\dd x
  \le 4\pi C_0 M_n(1-e^{-L})\longrightarrow 0 .
\]
The bound is uniform in the growing width $L$ because the nonconstant modes decay away from their respective ends. The reverse triangle inequality gives $\bigl||\hat P_n|-|a_{e,n}|\bigr|\le|\hat P_n-a_{e,n}|$. Thus the actual $|P_n|$-mass of the annulus differs from $2\pi L|a_{e,n}|$ by at most this same error.

By angular averaging, $2\pi L|a_{e,n}|\le\|P_n\|_1$, and $\|P_n\|_1\to E$. Since there are finitely many edges, we may take one subsequence on which $2\pi L|a_{e,n}|\to m_e\ge0$ for every edge $e$. The preceding error estimate then gives limiting mass $m_e$ on the whole annulus. For any interval $I\subset[0,L]$, the constant differential has mass $2\pi |I|\,|a_{e,n}|$ on the band $x\in I$. The integrated error on the band is bounded by the error on the whole annulus. In particular, taking $I=(uL,vL)$ gives limiting mass $(v-u)m_e$ for every fixed $0<u<v<1$ on that same subsequence. The total variation estimate \eqref{eq:tv} transfers both the annulus and the band limits to $\mu_n$. Point~1 shows that the complementary pieces have vanishing total mass. Consequently
\begin{equation}\label{eq:masses}
                             \sum_e m_e=E.
\end{equation}
If there are no edges this already contradicts $E>0$.

\textbf{3. Produce a genuine flat cylinder carrying the prescribed mass.}
Fix an edge with $m_e>0$ and write $L=L_{e,n}$. Point~2 gives the limiting $\mu_n$-mass of bands in the logarithmic coordinate $w$. We first construct a coordinate in which $q_n$ has constant coefficient, and then compare the resulting cylinder with those bands.

The middle third $L/3<\Re w<2L/3$ has $\mu_n$-mass tending to $m_e/3>0$. Its $|q_n|$-area is therefore bounded below by some $\delta_e>0$, since \eqref{eq:tv} bounds its mass by a function of its area tending to zero with the area. Put $b_{e,n}=a(\hat q_n)$, the constant Laurent coefficient of $\hat q_n$. By \eqref{eq:laurent}, $|\hat q_n-b_{e,n}|=O(e^{-L/3})$ on the middle third. Integrating over its coordinate area $2\pi L/3$ shows that its $|q_n|$-area differs from $(2\pi L/3)|b_{e,n}|$ by $O(Le^{-L/3})=o(1)$. Thus the constant term must account for a positive amount of area. Angular averaging on the full edge and the unit area of $q_n$ also give $2\pi L|b_{e,n}|\le1$. Together these bounds give, for some $c_e>0$ and all large $n$,
\[
                      c_e/L\le |b_{e,n}|\le 1/(2\pi L) .
\]
Fix $0<\epsilon<1/4$. On the wider band $\epsilon L<\Re w<(1-\epsilon)L$, division of \eqref{eq:laurent} by $|b_{e,n}|$ gives the relative estimate $\hat q_n/b_{e,n}=1+O(L e^{-\epsilon L})$. For large $n$ the ratio takes values in a disk about one missing zero. Its square root $\varpi_n$, chosen near one, is therefore single-valued and periodic. Its mean over one period,
\[
     \bar\varpi_n=(2\pi i)^{-1}\int_w^{w+2\pi i}\varpi_n(s)\,\dd s,
\]
is independent of $w$ because the derivative of the integral is $\varpi_n(w+2\pi i)-\varpi_n(w)=0$. The mean is $1+O(L e^{-\epsilon L})$. On the lifted strip define $\zeta$ by $\zeta'=\varpi_n/\bar\varpi_n$ and fix $\zeta(L/2)=L/2$. Dividing by the mean ensures that one turn around the annulus gives exactly the translation $2\pi i$. Since $\hat q_n=b_{e,n}\varpi_n^2$, we obtain
\[
 \zeta(w+2\pi i)=\zeta(w)+2\pi i,\qquad
 q_n=b_{e,n}\bar\varpi_n^2\,\dd\zeta^2 .
\]
Thus lengths in the $q_n$-metric are Euclidean lengths in the $\zeta$-coordinate multiplied by a single positive factor, even when the argument of $b_{e,n}$ varies with $n$.

Since $\zeta'-1=O(Le^{-\epsilon L})$, the Lipschitz constant of $\zeta-\mathrm{id}$ on the convex strip is less than one for large $n$. Hence $\zeta$ is injective, and the period identity gives injectivity on the quotient. Integrating $\zeta'-1$ from the midpoint in one fundamental rectangle and using periodicity gives $\sup|\zeta(w)-w|=O((L+1)Le^{-\epsilon L})=o(1)$. The image contains the straight strip $2\epsilon L<\Re\zeta<(1-2\epsilon)L$. Indeed each point $z_0$ of this narrower strip is at distance at least $\epsilon L$ from the boundary of the original strip. On $|w-z_0|=1$ we have $|\zeta(w)-w|<1$ for large $n$, so Rouché's theorem gives a zero of $\zeta(w)-z_0$ inside the circle.

Quotient the straight strip by its period $2\pi i$ and take its preimage under $\zeta$. This gives an embedded $q_n$-flat cylinder $C_{e,n}$ in the original edge. The strip retains width $(1-4\epsilon)L$ after trimming $2\epsilon L$ from each end. Width and circumference acquire the same metric factor, so their ratio, the modulus of $C_{e,n}$, is $(1-4\epsilon)L/(2\pi)\to\infty$.

To compute $\mu_n(C_{e,n})$, return to the $w$-coordinate of Point~2. Fix a sufficiently small $\delta'>0$. Since $\sup|\zeta(w)-w|=o(1)$, for all large $n$ the cylinder contains the $w$-band with fractional endpoints $2\epsilon+\delta',1-2\epsilon-\delta'$ and lies inside the band with endpoints $2\epsilon-\delta',1-2\epsilon+\delta'$. Applying Point~2 to the two bands gives
\[
 (1-4\epsilon-2\delta')m_e
 \le \liminf_n\mu_n(C_{e,n})
 \le \limsup_n\mu_n(C_{e,n})
 \le (1-4\epsilon+2\delta')m_e .
\]
Letting $\delta'\downarrow0$ gives the limiting mass $(1-4\epsilon)m_e$. The cylinders lie in their respective disjoint edge annuli. For each fixed $\epsilon$ the construction uses the subsequence already chosen in Point~2.
\end{proof}

\section{Why the cylinders cannot keep changing and existence of a fixed point}\label{sec:persistence}
Recall that $\Sigma$ is the initial punctured target representing $o$, with unit-area differential $q$. For $K>1$ define $\mathcal U(t,K)$ on $\Sigma$ to be the inverse image under $\phi_t$ of the union of all \emph{open maximal} $q_t$-flat cylinders of modulus strictly greater than $K$. By Section~\ref{sec:flat} it is a union of at most $r$ pairwise disjoint maximal $q$-cylinders, each a fixed open subset of $\Sigma$ whose modulus at time $t$ is given by \eqref{eq:modulus}.

The measure $\mu$ was fixed before the times or cylinders were chosen. Section~\ref{sec:mass} gives concentration along a subsequence of every sequence tending to infinity. We first turn this into a bound at \emph{all} sufficiently late times, so it applies at the later time when a chosen cylinder's modulus falls to $K$.

\textbf{Mass in cylinders of modulus greater than $K$.} For the measure $\mu$ fixed at the start of Section \ref{sec:mass},
\begin{equation}\label{eq:coverage}
 \mu(\mathcal U(t,K))\longrightarrow E\quad(t\to\infty).
\end{equation}
To verify the quantifiers, if \eqref{eq:coverage} fails there are $\delta_0>0$ and times $t_n\to\infty$ with complement mass at least $\delta_0$ at every $n$. Apply Proposition~\ref{l:cylinder-mass} to this sequence. Choose $0<\epsilon<1/4$ with $4\epsilon E<\delta_0$. On the extracted subsequence, the cylinder constructed in each positive-mass edge has modulus tending to infinity and sits in a maximal cylinder of at least its modulus (enlargement in the same direction preserves circumference). Thus all these cylinders are in the time-$t_n$ union for large $n$. The constructed subcylinders are disjoint, including if some have the same maximal extension. Their total $\mu_n$-mass tends to $(1-4\epsilon)\sum_e m_e=(1-4\epsilon)E$. Consequently the limsup of the complement mass along the extracted sequence is at most $4\epsilon E<\delta_0$. Since $\mu_n=(\phi_{t_n})_*\mu$, this contradicts the lower bound on the mass outside the union. This establishes \eqref{eq:coverage} for the real time parameter and with the same measure at every time.

\begin{figure}[ht]\centering
\begin{tikzpicture}
\begin{axis}[width=12cm,height=5.4cm,xmin=0,xmax=3.6,ymin=0,ymax=8.5,
  xlabel={$t$},ylabel={$m_C(t)$},axis lines=left,tick label style={font=\scriptsize},label style={font=\small},
  restrict y to domain=0:40,unbounded coords=jump,
  legend style={font=\scriptsize,draw=none,fill=none,at={(0.5,-0.32)},anchor=north,legend columns=3,/tikz/every even column/.append style={column sep=4mm}},legend cell align=left]
\addplot[domain=0:3.6,samples=300,pine,very thick] {0.3*exp(2*x)};
\addlegendentry{vertical core ($h=0$): grows like $e^{2t}$}
\addplot[domain=0:3.6,samples=300,rust,thick] {0.5/(0.0049*exp(2*x)+0.25*exp(-2*x))};
\addlegendentry{$C_1$: $(h,v)=(0.07,0.5)$}
\addplot[domain=0:3.6,samples=300,plum,thick,dashed] {0.043/(0.000036*exp(2*x)+0.36*exp(-2*x))};
\addlegendentry{$C_2$: $(h,v)=(0.006,0.6)$}
\addplot[domain=0:3.6,black!55,thin] {3};
\node[font=\scriptsize,anchor=south west] at (axis cs:0.02,3.02) {threshold $K$};
\draw[densely dotted,rust] (axis cs:1.74,0)--(axis cs:1.74,3);
\node[font=\scriptsize,rust,anchor=north] at (axis cs:1.74,0) {$t_*$};
\draw[densely dotted,rust] (axis cs:0.9,0)--(axis cs:0.9,7.1);
\node[font=\scriptsize,rust,anchor=north] at (axis cs:0.9,0) {$t_0$};
\end{axis}
\end{tikzpicture}
\caption{The nonvertical cylinder $C_1$ passes from modulus $>K$ at $t_0$ to $K$ at $t_*$. Its fixed $\mu$-mass then lies outside every maximal cylinder of modulus $>K$. The increasing curve illustrates the vertical case excluded in Section~\ref{sec:closed}; all curves follow \eqref{eq:modulus}.}
\label{fig:moduli}
\end{figure}
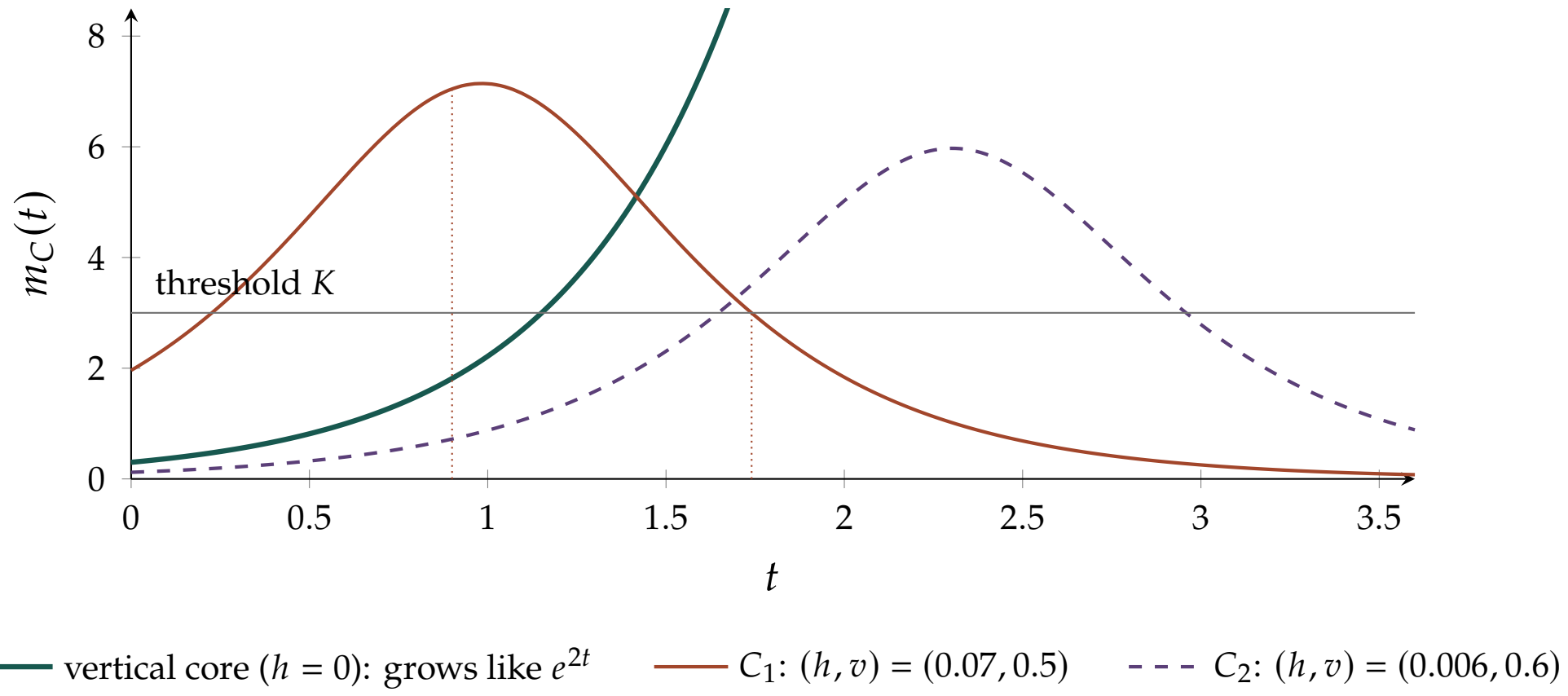

\textbf{The cylinder at modulus $K$.} By Proposition~\ref{l:cylinder-levy}, $q$ has no positive-width vertical cylinder under the joint hypotheses of Section~\ref{sec:mass}. Fix $K>1$ and $0<\delta<E/(r+1)$. By \eqref{eq:coverage} there is $t_0$ such that $\mu(\Sigma\setminus\mathcal U(t,K))<\delta$ at every time $t\ge t_0$. At time $t_0$ the at most $r$ cylinders forming $\mathcal U(t_0,K)$ carry $\mu$-mass more than $E-\delta$, so one of them, $C$, has $\mu(C)>(E-\delta)/r$. Since $C$ is not vertical, \eqref{eq:modulus} gives a time $t_*>t_0$ with $m_C(t_*)=K$. At time $t_*$ the cylinder $C$ is not part of $\mathcal U(t_*,K)$, and every cylinder of $\mathcal U(t_*,K)$ is disjoint from the open set $C$, because distinct maximal cylinders of modulus at least $K>1$ are disjoint (Section~\ref{sec:flat}). Hence $\mu(\Sigma\setminus\mathcal U(t_*,K))$ is at least $\mu(C)>(E-\delta)/r>\delta$, contradicting the choice of $t_0$. Therefore a transcendental map without a Levy cycle with $s\ge3$ has a fixed point.

\section{Necessity, uniqueness and convergence}\label{sec:necessity}
\subsection*{Two uses of contraction}\label{bg:rigidity}
\addcontentsline{toc}{subsection}{Two uses of contraction}
Necessity uses the hyperbolic metric on a punctured sphere: coverings are local isometries, while a proper holomorphic inclusion strictly shortens the metric by Schwarz--Pick. Convergence uses the Teichmüller metric on marked structures: compact balls about a fixed point are forward invariant, and a uniformly contracting iterate forces convergence there. For the general contraction framework, see \cite[\S3]{DH}, \cite[\S3]{HSS}, \cite[\S2.6]{MPR}. For polynomials with extra marks, we establish the required iterate below; its dependence on the marking is explained in Remark~\ref{r:poly}.

\emph{Necessity of the no-Levy-cycle condition.} Suppose that the transcendental map is realized by $h$. Transfer the marks and the curve classes through the equivalence and denote the finite marked set again by $A$. We may take $s\ge3$, since a sphere with at most three marks has no essential simple class. Set $X=\CC\setminus A$ and $Y=h^{-1}(X)$. Forward invariance gives $Y\subset X$. The inclusion is strict. In fact a transcendental entire function takes every finite value, with at most one exception, infinitely often near infinity by the great Picard theorem. Some $a\in A$ consequently has a preimage outside the finite set $A$. Both surfaces are hyperbolic and $h:Y\to X$ is a covering since all finite singular values are excluded. The deleted set $h^{-1}(A)$ is locally finite in the plane, so its complement is connected. The inclusion into the hyperbolic surface $X$ also makes $Y$ hyperbolic by uniformization: a plane universal cover of $Y$ would give a nonconstant map from $\CC$ to the disk by lifting the inclusion to the universal cover of $X$, contradicting Liouville's theorem.

For any essential target class take its closed hyperbolic geodesic in $X$. A degree-one preimage of the class is represented by a closed lift of this geodesic, by homotopy lifting. With the complete hyperbolic metric of $Y$ the covering is a local isometry, so this lifted representative has the original length. The inclusion $Y\hookrightarrow X$ strictly decreases the hyperbolic line element: lift it to a holomorphic map between the universal disks. If equality held at one tangent direction it would be a disk automorphism by the equality case of Schwarz--Pick. Its projection would then be onto $X$, contrary to the strict inclusion. If $\widetilde\gamma$ is the closed degree-one lift of $\gamma$, this gives
\[
  \ell_X([\widetilde\gamma])\le
  \operatorname{length}_X(\widetilde\gamma)<
  \operatorname{length}_Y(\widetilde\gamma)=\ell_X([\gamma]).
\]
Here the source class is taken in $X$, as required by the marked definition of a Levy cycle. Applying the strict inequality to the successive classes in a Levy cycle would return to the original geodesic length after finitely many strict decreases.

\emph{Existence and convergence.} As in Section~\ref{sec:records}, a fixed point is equivalent to a realization relative to $A$ (\mpr, Proposition~2.24). For $s\le2$ the marked Teichmüller space is a point and there are no essential curves. For $s\ge3$, Section~\ref{sec:persistence} gives existence for transcendental maps without a Levy cycle. For topological polynomials the marked polynomial criterion gives existence in the absence of a Levy cycle \cite[\S1.3 (the \"Ubertheorem) and Proposition~B.1]{BMW}. For $s\ge3$, necessity for a realized polynomial of degree $d$ also follows by the length argument above. The number of distinct points above $A$ is at least $ds-(d-1)>s$, since the total finite ramification is $d-1$. Thus $h^{-1}(X)\subsetneq X$ in this case as well.

Let $x_*$ be a fixed point. Every closed metric ball about $x_*$ is compact and forward invariant. Suppose that for some fixed $N\ge1$ the cotangent map of $\sigma^N$ strictly decreases the norm of every nonzero differential at every point. Compactness of the unit cotangent spheres and continuity of the derivative give a uniform operator bound $\lambda<1$ on each such ball. Integrating along a geodesic from $x_*$ to $x$ within the ball gives $d_T(\sigma^Nx,x_*)\le\lambda d_T(x,x_*)$. Iteration, followed by nonexpansion for the remaining fewer than $N$ steps, proves convergence to $x_*$ from every starting point and uniqueness. For transcendental maps, Lemma~\ref{l:pushforward} and the cotangent identification give this strict decrease with $N=1$ (\mpr, Proposition~2.22(2)).

For a polynomial of degree $d\ge2$, assume now $s\ge3$ and choose $N=\lfloor\log_d(s-2)\rfloor+1$. Put $D=d^N>s-2$. At any point of $\T$, compose the $N$ successive polynomial representatives to obtain $H$ of degree $D$. Through the successive markings, $H$ represents the topological iterate $f^N$. The cotangent identity proved after Lemma~\ref{l:pushforward} identifies the cotangent map of $\sigma^N$ with $H_*$ on the full marked spaces; see also \cite[Proposition~3.2]{DH}.

Take a nonzero source cotangent differential $q$, so that all its finite poles lie among the $s$ source marks, and put $p=H_*q$. The lemma gives $p$ integrable and holomorphic off the target marks, hence meromorphic with at most simple marked poles, including at infinity. Suppose now that $\|p\|_1=\|q\|_1$. Then $p\ne0$, and the equality case gives $q=cH^*p$; applying $H_*$ and using $H_*H^*p=Dp$ yields $c=1/D$.

The divisor of the nonzero quadratic differential $p$ on the sphere has degree $-4$. Since its poles are at most simple, it has at least four of them, at least three of them finite. Above three such poles there are at least $3D-2(D-1)=D+2$ distinct unramified preimages, as in step~6 of Proposition~\ref{p:degree}. The identity $q=H^*p/D$ makes each a finite pole of $q$. Since the finite poles of $q$ lie among the $s$ source marks, this gives the contradiction $D+2>s$.

Thus $H_*$ strictly decreases every nonzero cotangent norm, and the preceding compact-ball argument applies to $\sigma^N$.

\emph{Affine uniqueness.} For $s\ge2$ let $\phi_i f=h_i\psi_i$ be two realizations, with $\phi_i$ isotopic to $\psi_i$ relative to $A$. Their target markings represent the unique fixed point. After conjugating $h_1$ affinely and postcomposing both its markings by the same affine map, choose an isotopy $\Phi_u$, $0\le u\le1$, relative to $A$ with $\Phi_0=\phi_1$ and $\Phi_1=\phi_2$. The isotopy $\Phi_u\phi_1^{-1}$ starts at the identity and fixes $\phi_1(A)$ pointwise. Lift it through $h_1$ starting at the identity on the source, first on the covering of the complement of the target marks and then across the discrete fibers by isotopy lifting (\mpr, Proposition~2.13). Each point above a marked value stays fixed. The terminal lift composed with $\psi_1$ is a marking $\widetilde\psi_1$ agreeing with $\psi_2$ on $A$ and satisfying $h_1\widetilde\psi_1=\phi_2 f=h_2\psi_2$. Hence the plane homeomorphism $B=\widetilde\psi_1\psi_2^{-1}$ satisfies $h_1B=h_2$. On the complement of the inverse images of the marked values it is holomorphic by local inversion, and continuity removes its isolated singularities. It is an affine automorphism fixing the common images of at least two finite marks, so $B$ is the identity. The normalized functions therefore coincide. For $s=1$, \cite[Proposition~2.17]{MPR} gives the power-map class. Any realization is a polynomial with a totally ramified fixed point and hence is affinely conjugate to $z^d$. This proves Theorems~A and~B.

\begin{remark}\label{r:poly}
For the postcritical marking and hyperbolic orbifold, strict contraction of the second iterate is classical \cite[Proposition~3.3(b)]{DH}. The extension in \cite[Proposition~2.22(3)]{MPR} to arbitrary extra marks fails; our proof instead uses an iterate depending on $s$. For example, take $f(z)=z^2+i$, with critical orbit $0\mapsto i\mapsto-1+i\mapsto-i\mapsto-1+i$, and the forward-invariant set
\[
 A=\{0,\ i,\ -i,\ -1+i,\ 1-i,\ \pm\sqrt{1-2i}\},\qquad s=7.
\]
Put $p=\dd z^2/\big((z-i)(z+i)(z+1-i)\big)$, whose poles are $i,-i,-1+i,\infty$. The fibers of $f^2$ over these finite poles are, respectively, $\{\pm\sqrt{-i}\}$, $\{\pm i,\pm\sqrt{1-2i}\}$, and $\{0,\pm(1-i)\}$. The points $\pm\sqrt{-i}$ and $0$ have local degree two; infinity has local degree four. The pullback order rule from Proposition~\ref{p:degree}, Step~6 therefore makes $(f^2)^*p$ regular at those finite ramified points, with a double zero at infinity and exactly six simple poles $\pm i,\pm\sqrt{1-2i},\pm(1-i)$, all in $A$. Consequently $q=(f^2)^*p/4$ is integrable relative to $A^+$ and
$\|(f^2)_*q\|_1=\|p\|_1=\|q\|_1$.
Thus the cotangent map of $\sigma^2$ has norm one at the fixed point.
\end{remark}

\clearpage
\section{Relation to the real proof and to Hubbard--Schleicher--Shishikura}\label{sec:comparison}
Table~\ref{tab:comparison} compares the real proof \cite{vS}, the present argument and the exponential-map proof \cite{HSS}. Result numbers in the outer columns refer to the respective papers.

\begin{table}[!ht]\centering\footnotesize
\renewcommand{\arraystretch}{1.18}
\begin{tabularx}{\textwidth}{@{}>{\raggedright\arraybackslash}p{1.75cm}>{\raggedright\arraybackslash}X>{\raggedright\arraybackslash}X>{\raggedright\arraybackslash}X@{}}
\toprule
& \textbf{Real entire maps} \cite{vS} & \textbf{This paper} & \textbf{Exponential maps} \cite{HSS}\\
\midrule
Space & open simplex $W$ of ordered real configurations & $\T=\mathrm{Teich}(\widehat\CC,A^+)$ & Teichmüller space of the marked topological exponential map\\
Nearness to the boundary & smallest gap $m(x)\asymp\operatorname{dist}(x,\partial W)$ (Lemma~4.9) & $\frac12\log\Ext_xG$ for one foliation $G$, a horofunction up to a bounded error (Theorem~\ref{l:record}, Remark~\ref{r:wolff}) & not used\\
No fixed point gives & $x$ near $\partial W$ with $\|T(x)-x\|_\infty\le Km(x)$ (Theorem~5.1) & $\Ext_{\sigma^kx}G\le C\Ext_xG$ for all $x$ and all $k$ (Theorem~\ref{l:record}) & no analogue; only $d(\tau,\sigma\tau)\le d_0$ along the images of a path, by nonexpansion\\
Route to existence & Brouwer-type criterion; no contraction estimate & Wolff--Denjoy-type theorem, then contradiction; no uniform contraction & $\|d\sigma^k\|\le r<1$ where $d(\tau,\sigma\tau)<d_0$ (Proposition~3.3); Cauchy orbit\\
Short object & a shortest edge and its images under the combinatorial map & cores of vertical cylinders (\S\ref{sec:closed}); cores of annuli about clusters and their lifts (Lemma~\ref{l:components}) & cores of fat annuli on which the composition is injective (Proposition~3.2)\\
Use of no Levy cycle & every edge meets a critical label within $N$ steps (Definition~3.6, Lemma~3.7) & $\rho(M_\Lambda)<1$; in a chain of compatible essential classes, each block of $r$ lifts includes a degree $\ge2$ & a chain of short injective annuli would repeat a class (Proposition~3.3)\\
Effect of local degree & an edge of length $\ell$ at a critical endpoint has image of length $\le\frac12\|F''\|_\infty\ell^2$ (Proposition~12.2) & a closed lift of degree $d$ of an annulus of modulus $M$ has modulus $M/d$ & no critical points\\
Analytic estimate & the double-logarithmic bound on the hyperbolic diameter of the preimage interval gives polylogarithmic bounds on $\|F'\|_\infty$ and $\|F''\|_\infty$ after rescaling (Theorem~10.1, Proposition~11.1) & at fixed depth, the coefficient of $e^{2t}$ in reciprocal extremal length on the four-marked sphere is the sum of $M/d$ over finite cylinder lifts in the selected class (Proposition~\ref{p:quartet}) & near equality for the pushforward gives a disk on which the composition is injective (Proposition~3.2)\\
Final step & $1\le C_K\,m(1+\log(1/m))^B$ excludes small $m=m(x)$ at bounded relative displacement (Theorem~12.3) & for a cylinder core $\gamma$, $0<c_\gamma\le(M_\Lambda^k m)_\gamma$ for all $k$, whereas $M_\Lambda^k\to0$ without a Levy cycle (Proposition~\ref{l:cylinder-levy}) & the contraction gives a Cauchy orbit\\
Mass of quadratic differentials & --- & if $G$ has no cylinder, one inverse branch and the cluster estimates locate the mass of a fixed measure. The same cylinder retains its mass as its modulus falls (\S\S\ref{sec:transport}--\ref{sec:persistence}) & thick--thin decomposition, limit models and exp-efficient pushforwards (Sections~4--6)\\
Uniqueness and convergence & strict decrease; for polynomials $D^{n_0}+2>s$ (Proposition~14.3) & compact balls and \mpr; for polynomials $d^N+2>s$, the same count (\S\ref{sec:necessity}) & a contraction has at most one fixed point\\
Necessity & interval argument and subhyperbolicity; Schwarz--Pick form (Proposition~3.3, Remark~3.4) & hyperbolic lengths under $Y\subset X$ (\S\ref{sec:necessity}) & annulus moduli; the same hyperbolic-length argument (remark in \S3.4)\\
\bottomrule
\end{tabularx}
\caption{Comparison of the three proofs. The real derivative bounds keep distortion along a short interval from offsetting the quadratic gain at a critical point. In the cylinder case, the four-point upper bound shows that only finite lifts of vertical cylinders contribute to the leading modulus in the selected class.}
\label{tab:comparison}
\end{table}
\FloatBarrier

\subsection{What the proof takes from the real case}\label{sec:fromreal}
In the real proof, absence of a fixed point gives configurations near the boundary with $\|T(x)-x\|_\infty\le K m(x)$, where $m(x)$ is the smallest marked gap \cite[Lemma~4.9 and Theorem~5.1]{vS}. This is the \mbox{de Melo--van Strien} displacement criterion \cite[Section~II.4]{dMvS}, \cite[Proposition~8.1]{Sma}. With the enlarged marking, absence of Levy cycles forces a short interval eventually to meet a critical label. Taylor's theorem then gives a quadratic gain, while logarithmic derivative bounds control the intervening distortion and yield a contradiction \cite[Lemma~3.7, Theorem~10.1, Proposition~11.1 and Theorem~12.3]{vS}.

Here Theorem~\ref{l:record} gives one foliation $G$ whose extremal length is bounded by $Ce^{-2t}$ at every pullback depth. For a cylinder core $\gamma$, Proposition~\ref{p:quartet} converts this into the matrix bound $0<c_\gamma\le(M_\Lambda^k m)_\gamma$ of Proposition~\ref{l:cylinder-levy}. Without a Levy cycle, matrix decay contradicts this bound. Thus control uniform under iteration replaces the real proof's logarithmic derivative estimates.

For a real-symmetric map with $A\subset\RR$, this cylinder case suffices. Let $\rho$ be reflection and choose a reflection-fixed base point. Commutation $\rho\sigma=\sigma\rho$ makes the limiting set $K_0$ in Theorem~\ref{l:record} reflection invariant. If $G_0$ is the basis sum constructed there, then $G_0$ and $\rho G_0$ are compatible, since their summands belong to $K_0$, and both satisfy \eqref{eq:displacement} with the original constant $C_0$. For $G=G_0+\rho G_0$, additivity of intersection and duality give
\[
 \sqrt{\Ext_{\sigma^k x}G}
 \le\sqrt{C_0}\bigl(\sqrt{\Ext_xG_0}+\sqrt{\Ext_x\rho G_0}\bigr)
 \le2\sqrt{C_0}\sqrt{\Ext_xG}.
\]
Thus $G$ satisfies \eqref{eq:displacement} with $C=4C_0$ and may be normalized at $o$. A nonzero reflection-invariant measured foliation on a sphere with all marks real is a weighted multicurve, so it has a positive-width cylinder.

In the complex case, markings can diverge by twisting while the unmarked configuration stays fixed. Theorem~\ref{l:record} supplies $G$ without relying on real order; Sections~\ref{sec:transport}--\ref{sec:persistence} handle the case where it has no cylinder.

\subsection{Relation to Hubbard--Schleicher--Shishikura}\label{sec:toHSS}
Both proofs use pushforward as the cotangent map, strict norm decrease, and near-equality to locate the mass of quadratic differentials \cite[\S3.2]{HSS}. The integrable component limits and annular models $\dd z^2/z^2$ of \cite[Theorem~5.2; cf.\ Proposition~4.2]{HSS} correspond to our cluster charts and Laurent estimate \eqref{eq:laurent}.

HSS obtain a uniform contraction factor for an iterate on the region of bounded displacement, making the pullback orbit Cauchy \cite[Proposition~3.3]{HSS}. Here Theorem~\ref{l:record} supplies a fixed foliation and an auxiliary ray. Proposition~\ref{l:transport} yields near-equality and a single limiting measure on the initial surface. This additional conclusion lets us follow the same cylinder at two times: its mass stays fixed while its modulus falls, giving the contradiction in Section~\ref{sec:persistence}.

For exponential maps, the disk supplied by \cite[Proposition~3.2]{HSS} is mapped injectively by the relevant composition, so every boundary lift has degree one. Proposition~\ref{p:degree} instead gives degree $1\le d\le s-2$. The block argument in Lemma~\ref{l:components} uses this degree bound to force a Levy cycle. The inverse-branch series, univalent normal families and polynomial cap apply to general entire maps; the logarithmic half-plane retains the infinite-area argument at the start of the proof of \cite[Proposition~6.1]{HSS}.

\section*{Use of artificial intelligence}
\addcontentsline{toc}{section}{Use of artificial intelligence}
After completing the real-case manuscript \cite{vS}, I used OpenAI's ChatGPT to explore a complex analogue of its displacement condition $\|T(x)-x\|_\infty\le Km(x)$, where the smallest marked gap $m(x)$ is comparable to the distance to the boundary. This led, through several iterations, to Theorem~\ref{l:record}. I next investigated whether it would allow a different strategy for the real case, replacing analytic results about real entire maps by geometric methods, leading to Section~\ref{sec:closed}. I then explored whether it would improve the near-equality of the pushforward used by Hubbard, Schleicher and Shishikura~\cite{HSS}, building on Douady--Hubbard's exact equality case~\cite{DH}. ChatGPT assisted in developing these arguments.

I also used Anthropic's Claude (Claude Opus~5.5, through Claude Code) to check an earlier draft line by line, identify gaps and errors, and develop corrections. It helped simplify proofs, reorganize the exposition, prepare figures and LaTeX, check cited results against accessible sources, draft the comparison section and background notes, and review rewritten arguments in separate passes.

I take full responsibility for the mathematical statements, proofs, references and any remaining errors. The arguments are presented in the paper, with external results attributed to the cited literature.


\begin{thebibliography}{CFPP12}
\bibitem[AIM]{AIM}
K.~Astala, T.~Iwaniec and G.~Martin.
\emph{Elliptic Partial Differential Equations and Quasiconformal Mappings in the Plane}.
Princeton Mathematical Series, vol.~48. Princeton University Press, 2009.

\bibitem[Ant]{Ant}
S.~M.~Antonakoudis. \emph{The complex geometry of Teichmüller spaces and bounded symmetric domains}.
In \emph{In the Tradition of Ahlfors--Bers, VII}, Contemporary Mathematics \textbf{696} (2017), 1--16.
\href{https://arxiv.org/abs/1510.07340v1}{arXiv version}; page pinpoints refer to this version.

\bibitem[Bea]{Bea}
A.~F.~Beardon. \emph{The dynamics of contractions}.
Ergodic Theory Dynam. Systems \textbf{17} (1997), 1257--1266.

\bibitem[BMW]{BMW}
J.~Belk, D.~Margalit and R.~R.~Winarski.
\emph{Thurston's theorem and the Nielsen--Thurston classification via Teichmüller's theorem}.
arXiv:2309.06993 (2023). \href{https://arxiv.org/abs/2309.06993}{Preprint}.

\bibitem[BCT]{BCT}
X.~Buff, G.~Cui and L.~Tan. \emph{Teichmüller spaces and holomorphic dynamics}.
In A.~Papadopoulos (ed.), \emph{Handbook of Teichmüller Theory}, Vol.~IV, IRMA Lectures in Mathematics and Theoretical Physics, vol.~19. European Mathematical Society, 2014, 717--756.
\href{https://www.math.univ-toulouse.fr/~buff/Preprints/Teichmuller/Teichmuller.pdf}{Author version}.

\bibitem[CFPP]{CFPP}
J.~W.~Cannon, W.~J.~Floyd, W.~R.~Parry and K.~M.~Pilgrim.
\emph{Nearly Euclidean Thurston maps}.
Conform. Geom. Dyn. \textbf{16} (2012), 209--255.
\href{https://personal.math.vt.edu/floyd/research/papers/netmaps.pdf}{Article}.

\bibitem[dMvS]{dMvS}
W.~de Melo and S.~van Strien. \emph{One-Dimensional Dynamics}.
Ergebnisse der Mathematik und ihrer Grenzgebiete (3), vol.~25. Springer, 1993.

\bibitem[DH]{DH}
A.~Douady and J.~H.~Hubbard.
\emph{A proof of Thurston's topological characterization of rational functions}.
Acta Math. \textbf{171} (1993), 263--297.
\href{https://pi.math.cornell.edu/~hubbard/ThurstonRatMaps.pdf}{Article}.

\bibitem[DLR]{DLR}
M.~Duchin, C.~J.~Leininger and K.~Rafi. \emph{Length spectra and degeneration of flat metrics}.
Invent. Math. \textbf{182} (2010), 231--277.
\href{https://www.math.toronto.edu/~rafi/Papers/Flat-Metrics.pdf}{Article}.

\bibitem[DS]{DSou}
S.~D'Souza. \emph{A partial characterization of cosine Thurston maps}.
arXiv:2504.01273 (2025). \href{https://arxiv.org/abs/2504.01273}{Preprint}.

\bibitem[Dur]{Dur}
P.~L.~Duren. \emph{Univalent Functions}.
Grundlehren der mathematischen Wissenschaften, vol.~259. Springer, 1983.
\href{https://link.springer.com/book/9780387907956}{Publisher}.

\bibitem[FMarg]{FMar}
B.~Farb and D.~Margalit. \emph{A Primer on Mapping Class Groups}.
Princeton Mathematical Series, vol.~49. Princeton University Press, 2012.
\href{https://press.princeton.edu/books/ebook/9781400839049/a-primer-on-mapping-class-groups-pdf}{Publisher}.

\bibitem[FMas]{FM}
B.~Farb and H.~Masur.
\emph{Teichmüller geometry of moduli space, I: Distance minimizing rays and the Deligne--Mumford compactification}.
J. Differential Geom. \textbf{85} (2010), 187--227. arXiv:0803.0141.

\bibitem[FLP]{FLP}
A.~Fathi, F.~Laudenbach and V.~Poénaru. \emph{Travaux de Thurston sur les surfaces}.
Astérisque 66--67, Soc.\ Math.\ France, 1979.

\bibitem[GM]{GM}
F.~P.~Gardiner and H.~Masur. \emph{Extremal length geometry of Teichmüller space}.
Complex Variables Theory Appl. \textbf{16} (1991), 209--237.

\bibitem[Hat]{Hat}
A.~Hatcher. \emph{Algebraic Topology}.
Cambridge University Press, 2002.
\href{https://pi.math.cornell.edu/~hatcher/AT/AT.pdf}{Author version}.

\bibitem[Hub]{Hub}
J.~H.~Hubbard. \emph{Teichmüller Theory and Applications to Geometry, Topology, and Dynamics}, Vol.~2: Surface Homeomorphisms and Rational Functions. Matrix Editions, 2016.

\bibitem[HM]{HM}
J.~Hubbard and H.~Masur. \emph{Quadratic differentials and foliations}.
Acta Math. \textbf{142} (1979), 221--274.
\href{https://pi.math.cornell.edu/~hubbard/HubbardMasurFol.pdf}{Article}.

\bibitem[HSS]{HSS}
J.~Hubbard, D.~Schleicher and M.~Shishikura.
\emph{Exponential Thurston maps and limits of quadratic differentials}.
J. Amer. Math. Soc. \textbf{22} (2009), 77--117.
\href{https://pi.math.cornell.edu/~hubbard/ExpoSpidersJAMS.pdf}{Article}.

\bibitem[Kar1]{Kar}
A.~Karlsson. \emph{Non-expanding maps and Busemann functions}.
Ergodic Theory Dynam. Systems \textbf{21} (2001), 1447--1457.

\bibitem[Kar2]{Kar14}
A.~Karlsson. \emph{Two extensions of Thurston's spectral theorem for surface diffeomorphisms}.
Bull. Lond. Math. Soc. \textbf{46} (2014), 217--226. arXiv:1208.3735.

\bibitem[Ker]{Ker}
S.~P.~Kerckhoff. \emph{The asymptotic geometry of Teichmüller space}.
Topology \textbf{19} (1980), 23--41.
\href{https://doi.org/10.1016/0040-9383(80)90029-4}{Article}.

\bibitem[Kin]{Kin}
J.~Kinnunen. \emph{Sobolev spaces}.
Aalto University lecture notes, 2026 version (accessed 19 November 2026).
\href{https://math.aalto.fi/~jkkinnun/files/sobolev_spaces.pdf}{Notes}.

\bibitem[LS]{LS}
L.~Liu and W.~Su. \emph{The horofunction compactification of the Teichmüller metric}.
In A.~Papadopoulos (ed.), \emph{Handbook of Teichmüller Theory}, Vol.~IV, IRMA Lectures in Mathematics and Theoretical Physics, vol.~19. European Mathematical Society, 2014, 355--374.
\href{https://doi.org/10.4171/117-1/9}{Published chapter}.

\bibitem[Min]{Min}
Y.~N.~Minsky. \emph{Teichmüller geodesics and ends of hyperbolic 3-manifolds}.
Topology \textbf{32} (1993), 625--647.

\bibitem[Min99]{Min99}
Y.~N.~Minsky. \emph{The classification of punctured-torus groups}.
Ann. of Math. (2) \textbf{149} (1999), 559--626.
\href{https://arxiv.org/abs/math/9807001}{Preprint}.

\bibitem[Miy]{Miy}
H.~Miyachi. \emph{Teichmüller rays and the Gardiner--Masur boundary of Teichmüller space}.
Geom. Dedicata \textbf{137} (2008), 113--141.

\bibitem[MPR]{MPR}
M.~Mukundan, N.~Prochorov and B.~Reinke.
\emph{Dynamical approximations of postsingularly finite entire maps}.
To appear in \emph{Transactions of the American Mathematical Society}. arXiv:2305.17793, version~2 (2024).
\href{https://arxiv.org/abs/2305.17793v2}{Preprint}.

\bibitem[OP]{OP}
M.~Oliva and M.~Prats. \emph{Sharp bounds for composition with quasiconformal mappings in Sobolev spaces}.
J. Math. Anal. Appl. \textbf{451} (2017), 1026--1044.
\href{https://arxiv.org/abs/1612.00689}{Author version}.

\bibitem[Pr1]{ProSmall}
N.~Prochorov. \emph{Finite and infinite degree Thurston maps with a small postsingular set}.
To appear in the \emph{Journal of the London Mathematical Society}. arXiv:2410.01146 (2024). \href{https://arxiv.org/abs/2410.01146}{Preprint}.

\bibitem[Pr2]{ProExtra}
N.~Prochorov. \emph{Finite and infinite degree Thurston maps with extra marked points}.
arXiv:2410.06206 (2024). \href{https://arxiv.org/abs/2410.06206}{Preprint}.

\bibitem[Ran]{Ran}
T.~Ransford. \emph{Potential Theory in the Complex Plane}.
London Mathematical Society Student Texts, vol.~28. Cambridge University Press, 1995.

\bibitem[Rees]{Rees}
M.~Rees. \emph{An alternative approach to the ergodic theory of measured foliations on surfaces}.
Ergodic Theory Dynam. Systems \textbf{1} (1981), 461--488.
\href{https://doi.org/10.1017/S0143385700001383}{Article}.

\bibitem[Sch]{SchTalk}
D.~Schleicher. \emph{Transcendental Thurston Theory in Toulouse. A topological characterization of postsingularly finite meromorphic functions}.
Talk, Toulouse, 13 June 2025, reporting joint work with S.~Shemyakov.
\href{https://indico.math.cnrs.fr/event/12318/attachments/5832/10134/JHH80_Dierk\%20Schleicher.pdf}{Conference slides}.

\bibitem[Sh]{Shem}
S.~Shemyakov. \emph{A topological characterization of certain postsingularly finite entire functions: transcendental dynamics and Thurston theory}.
PhD thesis, Aix-Marseille Université, 2022. \href{https://theses.fr/2022AIXM0017.pdf}{theses.fr/2022AIXM0017}.

\bibitem[ST]{ShiTan}
M.~Shishikura and Tan Lei. \emph{A family of cubic rational maps and matings of cubic polynomials}.
Experiment. Math. \textbf{9} (2000), 29--53.
\href{https://doi.org/10.1080/10586458.2000.10504634}{Article}.

\bibitem[Sma]{Sma}
D.~Smania. \emph{Phase space universality for multimodal maps}.
Bull. Braz. Math. Soc. (N.S.) \textbf{36} (2005), 225--274.
\href{https://doi.org/10.1007/s00574-005-0038-y}{Article}.

\bibitem[Str]{Str}
K.~Strebel. \emph{Quadratic Differentials}.
Ergebnisse der Mathematik und ihrer Grenzgebiete (3), vol.~5. Springer, 1984.

\bibitem[Tho]{Tho}
C.~Thomassen. \emph{The Jordan--Sch\"onflies theorem and the classification of surfaces}.
Amer. Math. Monthly \textbf{99} (1992), 116--130.
\href{https://webhomes.maths.ed.ac.uk/~v1ranick/jordan/thomass.pdf}{Article}.

\bibitem[vS]{vS}
S.~van Strien. \emph{Thurston algorithm for real entire maps: existence and uniqueness}.
arXiv:2611.07818 (2026). \href{https://arxiv.org/abs/2611.07818}{Preprint}.

\bibitem[Wal]{Wal}
C.~Walsh. \emph{The asymptotic geometry of the Teichmüller metric}.
Geom. Dedicata \textbf{200} (2019), 115--152. \href{https://arxiv.org/abs/1210.5565}{arXiv version}; page pinpoints refer to this version.
\end{thebibliography}
\end{document}